\documentclass{amsart}
\usepackage{amssymb,latexsym,amsxtra,amscd,ifthen}
\usepackage{amsmath}
\usepackage{amsfonts}
\usepackage{verbatim}
\usepackage[arrow,matrix,cmtip]{xy}

\usepackage{fullpage}

\numberwithin{equation}{section}
\theoremstyle{plain}
\newtheorem{theorem}[equation]{Theorem}

\newtheorem{lemma}[equation]{Lemma}

\newtheorem{corollary}[equation]{Corollary}
\newtheorem{proposition}[equation]{Proposition}

\theoremstyle{definition}
\newtheorem{definition}[equation]{Definition}
\newtheorem{notation}[equation]{Notation}

\newtheorem{assumptions}[equation]{Assumptions}

\newtheorem{example}[equation]{Example}
\newtheorem{examples}[equation]{Examples}

\newtheorem{remark}[equation]{Remark}
\newtheorem{remarks}[equation]{Remarks}
\newcommand{\HB}{\operatorname{H}}

\newcommand{\Morph}{\operatorname{Morph}}

\newcommand{\Mod}{\operatorname{-Mod}}

\newcommand{\GT}{\operatorname{GT}}

\newcommand{\bpr}{\begin{proof}}
\newcommand{\epr}{\end{proof}}

\newcommand{\reg}{\operatorname{reg}}

\newcommand{\spec}{\operatorname{Spec}}
\newcommand{\ext}{\operatorname{Ext}}

\newcommand{\lra}{\longrightarrow}

\newcommand{\ra}{\rightarrow}

\newcommand{\mc}{\mathcal}
\newcommand{\mf}{\mathfrak}

\newcommand{\mb}{\mathbb}

\renewcommand{\hom}{\operatorname{Hom}}

\newcommand{\wt}{\widetilde}
\newcommand{\ptfn}{\operatorname{\mathcal{P}}}
\newcommand{\qptfn}{\operatorname{\mathcal{P}_{\mathrm{qgr}}}}

\newcommand{\coH}{\operatorname{H}}

\newcommand{\tor}{\operatorname{Tor}}

\newcommand{\cd}{\operatorname{cd}}
\newcommand{\gld}{\operatorname{gld}}

\newcommand{\rproj}{\operatorname{proj-\negmedspace}}

\newcommand{\rcoh}{\operatorname{coh}}

\newcommand{\rmod}{\operatorname{mod-\negmedspace}}

\newcommand{\catmod}{\operatorname{-mod}}
\newcommand{\catMod}{\operatorname{-Mod}}

\newcommand{\rcatMod}{\operatorname{Mod-\negmedspace}}
\newcommand{\Gr}{\operatorname{-Gr}}

\newcommand{\qgr}{\operatorname{-qgr}}
\newcommand{\rGr}{\operatorname{Gr-\negmedspace}}
\newcommand{\rgr}{\operatorname{gr-\negmedspace}}
\newcommand{\rQgr}{\operatorname{Qgr-\negmedspace}}
\newcommand{\rqgr}{\operatorname{qgr-\negmedspace}}

\newcommand{\rTors}{\operatorname{Tors-\negmedspace}}

\newcommand{\aut}{\operatorname{Aut}}

\newcommand{\supp}{\operatorname{Supp}}
\newcommand{\cosupp}{\operatorname{coSupp}}

\newcommand{\nsr}{na{\"\i}ve blowup algebra}

\newcommand{\gnbba}{generalized \naive\ blowup bimodule algebra}
\newcommand{\gnba}{generalized \naive\ blowup algebra}
\newcommand{\gns}{generalized na{\"\i}ve sequence}
\newcommand{\mcII}{\mc{I}\hskip -4pt \mc{I}}

\providecommand{\qedhere}{\qed}
\newboolean{ams2.0}
\ifthenelse{\equal{\qedhere}{\qed}}%
    {\setboolean{ams2.0}{false}}%
    {\setboolean{ams2.0}{true}}%

\DeclareMathOperator{\calHom}{\mathcal{H}\mathit{om}}

\DeclareMathOperator{\calTor}{\mathcal{T}\mathit{or}}
\DeclareMathOperator{\Hom}{Hom}
\DeclareMathOperator{\Ext}{Ext}

\DeclareMathOperator{\Aut}{Aut}

\DeclareMathOperator{\Supp}{Supp}

\newcommand{\op}{\text{op}}

\newcommand{\NN}{{\mathbb N}}

\newcommand{\LL}{{\mathcal L}}

\newcommand{\calR}{{\mathcal R}}

\newcommand{\OO}{{\mathcal O}}

\newcommand{\bound}{bounding}
\newcommand{\naive}{na{\"\i}ve}
\newcommand{\Naive}{Na{\"\i}ve}

\newcommand{\lcl}[2]{#1_{{\scriptscriptstyle [}#2{\scriptscriptstyle ]}}}

\begin{document}
\title[Naive Noncommutative Blowing
 Up]{Na{\"\i}ve Noncommutative Blowups  at  Zero-Dimensional  Schemes}

\author{D. Rogalski and J. T. Stafford}

\address{(Rogalski)
Department of Mathematics, UCSD, La Jolla, CA 92093-0112, USA. }
\email{drogalsk@math.ucsd.edu}
\thanks{The first author was  partially supported by the NSF  through grants
DMS-0202479  and  DMS-0600834   while the second author was
  partially supported by the NSF through grants DMS-0245320 and
  DMS-0555750  and also  by the  Leverhulme Research Interchange Grant
F/00158/X.  Part of this work was written up while the second author
was visiting and partially supported by the Newton Institute, Cambridge. We
would like to thank all three institutions for their financial
support.}

\address{(Stafford) Department of Mathematics,
University  of Michigan, Ann Arbor, MI 48109-1043, USA.}
\email{jts@umich.edu}

\begin{abstract}  In an earlier paper \cite{KRS}
%  (D.\ S.\ Keeler, D.\ Rogalski, and J.\ T.\ Stafford,
% ``Na{\"\i}ve  noncommutative blowing up,'' \emph{Duke Math.\ J.},
 %  \textbf{126} (2005), 491-546),
 we defined and investigated the properties
  of the  \naive\ blowup of an integral projective scheme $X$ at a single closed point.
   In this paper we extend those results to
   the case when one \naive ly blows up $X$ at any
   suitably generic zero-dimensional subscheme $Z$.
 The resulting algebra $A$ has a number of curious properties; for example it is
  noetherian but never strongly noetherian and the point modules are never
parametrized by a projective scheme. This is despite the fact that the
category of torsion modules in $\rqgr A$ is equivalent to the category of torsion
coherent sheaves over $X$.

   These results are used in the companion paper \cite{RS}
%   ``A class of noncommutative projective surfaces''
   to prove that a large class of noncommutative surfaces can be written as \naive\ blowups.
\end{abstract}
 \maketitle
 \tableofcontents

  \clearpage

%%%%%%%%%%%%%%%%%%%%%%%%%%%%%
\section{Introduction}\label{intro}

The concept of a \naive\ blowup of a scheme $X$ was introduced in
\cite{KRS}, where it was  also shown that these objects had
properties quite unlike their commutative counterparts. \Naive\
blowups  are also  used  in the companion paper to this one,
 \cite{RS}, in order to classify a
large class of noncommutative algebras. Unfortunately, the algebras
considered in \cite{KRS} were obtained by \naive ly blowing up a
single closed point whereas the applications in  \cite{RS} require
one to \naive ly blow up any suitably general zero-dimensional
subscheme. The aim of this paper is therefore to study this more
general case.
 Before describing the results we need some definitions.

 Throughout, $k$ will be an algebraically closed base field. A  $k$-algebra $A$ is  called
 \emph{connected graded}
(\emph{cg})\label{cg-defn}   if $A=\bigoplus_{n\geq 0} A_n$, where
$A_0=k$ and  $\dim_kA_n<\infty$ for each $n$.  Given a cg
$k$-algebra $A$, the category of
noetherian graded $A$-modules modulo those of finite length is
written $\rqgr A$ (see page~\pageref{qgr-defn}  for more
details).  A \emph{point module}   is a
 cyclic graded $A$-module $M=\bigoplus_{n\geq 0} M_n$ such that $\dim_kM_n=1$
  for all $n\geq 0$.
A \emph{point module in $\rqgr A$} is defined to be the
   image in $\rqgr A$ of a cyclic graded $A$-module $M=\bigoplus_{n\geq 0} M_n$,
generated in degree zero,  such that $\dim_k M_n=1$
  for all $n\gg 0$.

 The underlying data for a \naive\ blowup is as follows.
 Fix an integral projective scheme  $X$, with automorphism  $\sigma$
and  $\sigma$-ample sheaf $\mc{L}$, as   defined
in~\eqref{sigma-ample}.
 Let $Z=Z_{\mc{I}}\subset X$ be  a  zero-dimensional subscheme,
  with defining ideal  $\mc{I}\subseteq \mc{O}_X$.   In a
manner reminiscent of the  Rees ring   construction of the blowup of
a commutative scheme,
 we form the \emph{bimodule algebra}
 $\mc{R}=\mc{R}(X, Z, \mc{L},\sigma) =\mc{O}_X\oplus
\mc{R}_1\oplus\mc{R}_2\oplus\cdots,\ $ where $\mc{R}_n =
\mc{L}_n\otimes_{\mc{O}_X}\mc{I}_n$, for $\mc{L}_n = \mc{L} \otimes
\sigma^* \mc{L} \otimes \dots \otimes (\sigma^{n-1})^*\mc{L}$ and
$\mc{I}_n = \mc{I} \cdot \sigma^* \mc{I} \cdots (\sigma^{n-1})^*
\mc{I}$.  %As is described in  Section~\ref{definitions}, this
 This bimodule algebra has a natural
multiplication and the  {\it \nsr\ of $X$ at $Z$} is then  the algebra of sections
  $$R \ = \ R(X,Z, \LL,\sigma)
   \ = \ \mathrm{H}^0(X, \, \mc{R})  \ = \  k\oplus \mathrm{H}^0(X, \, \mc{R}_1)
  \oplus\mathrm{H}^0(X, \, \mc{R}_2)\oplus\cdots $$

One can also form $R=R(X,Z,\mc{L},\sigma)$
 when  $Z$ is the empty set,  in which case
$R$ is simply \emph{the  twisted homogeneous coordinate ring}
$B(X,\mc{L},\sigma)$ from \cite{AV} that is so important in noncommutative projective geometry
(see \cite{SV}, for example).
% and its properties are extremely nice. For example
If $B=B(X,\mc{L},\sigma)$ for  a $\sigma$-ample invertible sheaf $\mc{L}$, then
$B$ has extremely pleasant properties, among which we mention:
\begin{enumerate}
 \item[(a)] \cite[Proposition~4.13]{ASZ}  $B$ is  \emph{strongly noetherian};
that is, for all commutative noetherian $k$-algebras $C$, the ring $B\otimes_kC$ is noetherian;
\item[(b)]   \cite[Theorem~1.3]{AV}   $\rqgr B \simeq \rcoh X$, the category of
coherent sheaves on $X$;
 \item[(c)]  at least when $B$ is generated in degree one,
  the set of point modules for $B$, both in $\rgr B$
 and in $\rqgr B$, is  parametrized by the scheme $X$
 (use   \cite[Theorems~1.1 and 1.2]{RZ}  and \cite[Proposition~10.2]{KRS});
  \item[(d)]    \cite[Theorem~7.3]{Ye}  $B$ has a balanced dualizing complex
in the sense of that paper.
\end{enumerate}
In contrast,  the \nsr\
$R(X,Z,\mc{L},\sigma)$ for $Z\not=\emptyset$ has properties quite unlike 
those just mentioned.  %of a  commutative cg algebra.
Some of these properties  are given by the next theorem and will be
discussed in more detail later in this introduction. 

In order to state the theorem, we need one more definition.
A set of closed points  $\mc{C}\subset X$
  is   \emph{critically dense} if  $\mc{C}$ is infinite and every infinite
subset of $\mc{C}$ has Zariski  closure equal to $X$.
A zero-dimensional subscheme
  $Z\subset X$ is \emph{right saturating},  respectively \emph{saturating}  if, for each point
$c\in S=\supp Z$, the   set $\{c_i=\sigma^{-i}(c) : i\geq 0\}$,
respectively  $\{c_i=\sigma^{-i}(c): i \in \mb{Z} \}$ is critically
dense. By convention, the empty subscheme is  (right)
 saturating.

\begin{theorem}\label{mainthm}
Let $X$ be an integral projective scheme  with $\dim X\geq 2$ and
$\sigma\in \Aut(X)$. Assume that $\mathcal L$ is a $\sigma$-ample
invertible sheaf on $X$ and that  $Z$ is a {\bf (}non-empty{\bf )}
zero-dimensional, saturating subscheme of $X$.
 If $\mc{R}=\mc{R}(X,Z,\mc{L},\sigma)$ with global
sections $R=R(X,Z,\mc{L},\sigma)$, then:
\begin{enumerate}

\item  $\rqgr R \simeq \rqgr {\mc{R}}$ and $\rqgr \mc{R}$ is independent
of the choice of $\mathcal L$.

\item $R$ is always a noetherian domain.

\item $R$ is never strongly noetherian.
% that is,
% there exists a commutative noetherian $k$-algebra $C$ such that
% $R\otimes_kC$ is not noetherian.

\item The category
 of (Goldie) torsion objects in $\rqgr R$, as defined in Section~\ref{R-modules},
  is equivalent to the category of torsion coherent
 $\mc{O}_X$-modules.

\item  In particular, the isomorphism classes of simple objects in $\rqgr R$,
which are also the point modules in $\rqgr R$, are
in 1-1 correspondence with the closed points of $X$.

\item
However,  the point modules in $\rqgr R$ are \emph{not} parametrized by
any scheme of locally finite type.
 When $R_1 \neq 0$, \
the point modules  in $\rgr R$  are not
parametrized by any scheme of locally finite type.

\item   $\rqgr R$ has finite cohomological dimension. If $X$
is smooth, $\rqgr R$ has finite
homological dimension.

\item If ${\coH}^1(R) ={\Ext}^1_{\rqgr R}(R,R)$,
then $\dim_k {\coH}^1(R)= \infty$. Consequently,  $R$ does
not have  a balanced dualizing complex.

\item  If $U$  is \emph{any}  open affine subset of $X$,
 then generic flatness,  as defined in Section~\ref{further-section},
 fails for the finitely generated $R\otimes {\mathcal
O}_X(U)$-module  ${\mathcal R}(U) =\bigoplus {\mathcal R}_n(U)$.
\end{enumerate}
\end{theorem}

\begin{remark}\label{intro-remark}
All the properties described by Theorem~\ref{mainthm}
pass easily to subrings of finite index; that is,  the theorem
 still holds if one replaces  $R$ by any cg $k$-algebra $R'\subset R$ such 
that $\dim_k R/R'<\infty$. See Remark~\ref{intro-remark-proof} for the details.
\end{remark}

Theorem~\ref{mainthm} summarizes many of the results of this paper
and so its proof takes up much of  the paper.  Specifically, parts~(1)
and  (2) of the theorem are proved  by combining 
Proposition~\ref{L no matter} and
 Theorem~\ref{general-ampleness} and   their proof takes up most of
Sections~\ref{definitions}--\ref{sect-ample-new}. The rest of the
paper is then concerned with applying this theorem to get a deeper
understanding of the properties of $\mc{R}$ and $R$. In particular,
parts~(3) and (9)  of Theorem~\ref{mainthm} are proved in
Theorem~\ref{not strong noeth}; part~(4)  in
 Theorem~\ref{GT equiv}; part~(5) by combining Corollary~\ref{GT equiv1} and
 Proposition~\ref{non-represent2};
 part~(6) in Theorem~\ref{non-represent};
part~(7) in  Theorem~\ref{fin-cohom-dim}; and part~(8)  in
Theorem~\ref{chi_1, not chi_2} and Remark~\ref{dualizing-remark}.

The proofs of many of the results intermediary to
Theorem~\ref{mainthm} are very similar to those of the analogous
results from \cite{KRS}, but there are some significant differences
which we now discuss. Part~(1) of the theorem is the fundamental
tool for studying the \nsr\ $R=R(X,Z,\mc{L},\sigma)$ since it
provides the avenue for introducing geometry into  that study. In
turn, the key step in its proof is the following result.

\begin{proposition}\label{ample-intro} {\rm (Theorem~\ref{general-ampleness})}
Keep the hypotheses of Theorem~\ref{mainthm}. Then the  sequence $\{\mc{R}_n=\mc{L}_n\otimes\mc{I}_n\}$ is
 ample in the sense   that, for every $\mc{F}\in \rcoh X$ and all  $n\gg 0$, the sheaf
 $\mc{F}\otimes \mc{R}_n$ is generated by its global  sections and satisfies
 $\mathrm{H}^i(X,\mc{F}\otimes \mc{R}_n)=0$ for all $i>0$.
 \end{proposition}

The proof of Proposition~\ref{ample-intro} is considerably more subtle
 than that of its counterpart  \cite[Proposition~4.6]{KRS}.
Even leaving aside the issue of the
number of points at which one is \naive ly blowing up,
  the proof in \cite{KRS}  only works when $\mc{L}$ is both
$\sigma$-ample and very ample. In contrast, in the application in
Theorem~\ref{RS-mainthm},  below, one is only allowed to assume that
$\mc{L}$ is $\sigma$-ample. The proof of
Proposition~\ref{ample-intro} in this more general case requires
a delicate analysis of the number of points separated by $\mc{L}_n$.

 As well as  its applications in this paper,
Proposition~\ref{ample-intro} is also the starting point for the
classification in \cite{RS} of a large class of noncommutative algebras:

\begin{theorem}\label{RS-mainthm} {\rm (}\cite[Theorem~1.1]{RS}{\rm )}
 Let $A$ be a cg noetherian
domain that is generated in degree one and assume that the graded quotient ring $Q(A)$
has the form  $Q(A) \cong k(Y)[t, t^{-1}; \sigma]$, where $\sigma$ is induced from
 an automorphism of the integral projective surface $Y$.
Then, up to a finite dimensional vector space, $A$ is isomorphic to either
$B(X,\mc{L},\sigma)$ or to $
R(X,Z,\LL,\sigma)$, for some projective surface  $X$ birational to
$Y,$ with an induced action of $\sigma$, a $\sigma$-ample invertible
sheaf $\mc{L}$ and a  zero-dimensional saturating
subscheme $Z$.
 \end{theorem}

One striking consequence of this theorem is the
fact that the properties  described
by  Theorem~\ref{mainthm} are not exceptional:  as soon as    $Y$
has at least one  critically dense $\sigma$-orbit then  Theorem~\ref{RS-mainthm}
and Remark~\ref{intro-remark}
imply  that,  generically,
each noetherian cg subalgebra of $k(Y)[t,t^{-1};\sigma]$ has these properties.

A surprising feature of \naive\ blowups at multiple points concerns
torsion extensions. Given a \nsr\ $R=R(X,Z,\mc{L},\sigma)$
satisfying the hypotheses of Theorem~\ref{mainthm}, the
\emph{maximal right  torsion extension}   of $R$ is defined to be
the ring $T=T(R)=\{x\in Q(R) : xR_{\geq n} \subseteq R\ \text{for
some } n\geq 0\}$.   If $R$ is the \nsr\ at a
single point (or a twisted homogeneous coordinate ring)  then $T(R)/R$ is always
finite dimensional \cite[Theorem~1.1(8)]{KRS}. In contrast, if one
blows up at multiple points then $T(R)/R$ can be infinite
dimensional (see Example~\ref{easy-eg1}) and, indeed, the ring $R$
from that example is even an \emph{idealizer subring} $R=\{\theta\in
T : M\theta\subseteq M\}$ for some left ideal $M$ of $T$ (see the
discussion after Lemma~\ref{ideal-lem}). This implies that \nsr s
can have very non-symmetric properties. It also implies that the
$\chi_1$ condition, as defined in Section~\ref{chi}, will fail for
such a ring $R$. In contrast, \cite[Theorem~1.1(8)]{KRS} shows that
the $\chi_1$ condition always holds when one \naive ly blows up a
single point.   The details behind these assertions are given in
Sections~\ref{gen naive blowups} and~\ref{chi}.

In order to describe the maximal right torsion extension of a \nsr,
and even to prove  parts of Theorem~\ref{mainthm},
one needs  to pass to a slightly larger class of algebras, called \gnba s.
These algebras are discussed in detail in Section~\ref{gen naive blowups} and,
because it takes little extra work,
 much of Theorem~\ref{mainthm} is actually proved at their level of
 generality.  

We end the introduction by describing the significance of some of
the other parts of Theorem~\ref{mainthm} and we keep the notation
and hypotheses from that result. Although part~(1)   justifies the
idea that $\rqgr R\simeq \rqgr \mc{R}$
 is a kind of a noncommutative
blowup of $X$, the category  $\rqgr R $ is in fact much closer to
  $\mathrm{coh}\, X$  than it is to  $\mathrm{coh}\, \wt{X}$
for the (classical) blowup $\wt{X}$ of $X$ at $Z$. For example,
suppose that $k(x)$ is  the skyscraper sheaf corresponding to
 a closed point $x\in \supp Z$. If one tensors $k(x)$ with the
 sheaf of Rees rings
$\mc{R}(X,Z,\mc{L}, \mathrm{Id})$ corresponding to $\wt{X}$ then, of
course, one  obtains the $\mc{O}_{\wt{X}}$-module corresponding to
an exceptional divisor on $\wt{X}$. In contrast,
$k(x)\otimes_{\mc{O}_X} \mc{R}$  
 is a finite direct sum of simple objects  from $ \rqgr\mc{R}$
(use the computation from \cite[Proposition~5.3]{KRS}).
 Although we take a different approach,
this argument can be used to prove much of  part~(5) of the theorem.
This result, together with its generalization in part~(4),
 shows   that the differences between the categories
 $\rqgr R$ and $\mathrm{coh}\, X$ are really quite subtle.

The idea of considering the strongly noetherian condition arose in the work of
Artin, Small and Zhang \cite{ASZ,AZ2},  who showed that many algebras have this
property and that this
 has a number of important consequences for the  algebras in question.
Notably,   a strongly noetherian graded $k$-algebra
  $A$  will always satisfy generic flatness \cite[Theorem~0.1]{ASZ}.
 If  $A$ is also generated in degree one then
 its point modules, both in $\rgr A$ and  in $\rqgr A$,
will be parametrized by a projective scheme (see
\cite[Corollary~E4.5]{AZ2}, respectively
\cite[Proposition~10.2]{KRS}). Thus part~(9) and both assertions from
part~(6) of Theorem~\ref{mainthm} all imply that $R$ is not strongly
noetherian.   Finally,
parts~(7) and (8) of the theorem show contrasting homological properties
of $\rqgr R$; in particular, part~(8) means that  the homological
machinery developed by Yekutieli and Zhang in their papers on
dualizing complexes cannot easily be applied to the study of $R$.
For an illustration of the complications this causes and the ways in
which one can circumvent some of these problems, see
Remark~\ref{non-rep-remark}.

Finally, we mention that several peripheral
 results that are stated but not proved in this paper
are proved in full generality in an appendix \cite{RS2}
 that will be available on the web but not published.

%%%%%%%%%%%%%%%%%%%%%%%%%%%%
\section{Definitions and background material  }\label{definitions}

 In this section we set up the appropriate
notation relating to  the  bimodule algebras
$\mc{R}= \mc{R}(X,Z,\mc{L},\sigma)$ and their   section rings $R=\mathrm{H}^0(X,\mc{R})$
and determine, among other things, when
$\mc{R}$ is noetherian (see Proposition~\ref{bimod alg noeth}).
With the exception of the proof of that proposition, most of this section is similar to the material in
 \cite[Sections~2 and~3]{KRS}, to which the reader is referred  for further
 details.

Fix throughout an integral projective scheme $X$ over an algebraically closed
field $k$.   The category of quasi-coherent, respectively coherent, sheaves on
$X$ will be written $\mc{O}_X\Mod$, respectively $\mc{O}_X\catmod$. We use the
following notation  for pullbacks: if $\sigma\in \Aut(X)$
is a $k$-automorphism of   $X$,  and $\mc{F}\in \mc{O}_X\catmod$, then
  $\mc{F}^{\sigma}=\sigma^{*}(\mc{F})$.
We adopt the usual convention that   $\sigma$
acts on functions by $f^\sigma(x) = f(\sigma(x))$, for $x\in X$.

\begin{definition}
\label{bimod def}
A \emph{coherent $\mc{O}_X$-bimodule} is a coherent sheaf $\mc{F}$
on  $X \times X$ such that  $Z = \supp{\mc{F}}$ has the property that both
projections  $\rho_1, \rho_2: Z \ra X$ are finite morphisms.
An \emph{$\mc{O}_X$-bimodule} is a  quasi-coherent sheaf $\mc{F}$
on  $X \times X$ such that every coherent $X \times X$-subsheaf is a
coherent $\mc{O}_X$-bimodule.  The left and right $\mc{O}_X$-module
structures associated to $\mc{F}$ are defined
to be $_{\mc{O}_X} \mc{F} = (\rho_1)_{*} \mc{F}$ and
$\mc{F}_{\mc{O}_X} = (\rho_2)_{*} \mc{F}$ respectively.

Given $\mc{F}\in \mc{O}_X\catmod$ and $\tau,\sigma\in \Aut(X)$,
 define an $\mc{O}_X$-bimodule ${}_\tau
\mc{F}_{\sigma}$ by $(\tau, \sigma)_*\mc{F}$ where
 $(\tau, \sigma): X \to X \times X$.
We usually  write $\mc{F}_{\sigma}$
for ${}_1 \mc{F}_{\sigma}$, where $1$ is the identity automorphism.
The reader may check that ${}_1
\mc{F}_{\sigma}$ has left $\mc{O}_X$-module structure
$\mc{F}$ but right $\mc{O}_X$-module structure $\mc{F}^{\sigma^{-1}}$.
\end{definition}

When no other bimodule structure  is
given, a   sheaf $\mc{F} \in \mc{O}_X\catmod$ will be  assumed to have the bimodule
structure $_1 \mc{F} _1$.  Thus all sheaves become bimodules, and all tensor
products can be thought of as tensor products of bimodules.
 Unless otherwise stated,  when thinking of a bimodule $\mc{G}$ as a sheaf, we mean the left
$\mc{O}_X$-module structure of $\mc{G}$. Thus,
 when we write $\coH^i(X,\mc{G})$ or say that $\mc{G}$ is generated by
its global sections we  are referring to the left structure of $\mc{G}$.
Working on the left will have notational advantages, but, as in \cite{KRS},
it  is otherwise not  significant.

The following  special case of Van den Bergh's bimodule algebras
  will form the main objects of interest in this
paper.

\begin{definition}\label{bimod-alg}  Let $\sigma\in\Aut(X)$.
A   \emph{graded $(\mc{O}_X, \sigma)$-bimodule  algebra} is an
$\mc{O}_X$-bimodule $\mc{B}=\bigoplus_{n\geq 0} \mc{B}_n$
with a unit map $1: \mc{O}_X \ra \mc{B}$
and a product map  $\mu: \mc{B} \otimes \mc{B} \ra \mc{B}$ satisfying the usual
axioms as well as:  \begin{enumerate}
  \item For each $n$,   $\mc{B}_n\cong {}_1(\mc{E}_n)_{\sigma^n}$, for some
$\mc{E}_n\in\ \mc{O}_X\catmod$ with
$\mc{B}_0={}_1(\mc{O}_X)_1$.
\item The multiplication map satisfies $\mu
(\mc{B}_m \otimes \mc{B}_n) \subseteq \mc{B}_{m+n}$ for all $m,n$
 and $1(\mc{O}_X) \subseteq
\mc{B}_0$. Equivalently $\mu$ is defined by   $\OO_X$-module maps
$\mc{E}_n\otimes \mc{E}_m^{\sigma^n}\to \mc{E}_{m+n}$ satisfying the
appropriate associativity conditions.
\end{enumerate}
We will write $\mc{B} =
\bigoplus{}_1(\mc{E}_n)_{\sigma^n}$ throughout the section.
\end{definition}

\begin{definition}
\label{graded def}
Let $\mc{B}$ be a graded $(\mc{O}_X, \sigma)$-algebra.  A
\emph{graded right $\mc{B}$-module}   is a quasi-coherent
right $\mc{O}_X$-module $\mc{M} = \bigoplus_{n \in \mb{Z}} \mc{M}_n$
 together with a right $\mc{O}_X$-module  map $\mu:
\mc{M} \otimes \mc{B} \ra \mc{M}$  satisfying the usual axioms.  The \emph{shift of}
$\mc{M} $ is defined by $\mc{M}[n]=\bigoplus \mc{M}[n]_i$ with
$\mc{M}[n]_i=\mc{M}_{i+n}$.

The $\mc{B}$-module $\mc{M}$ is \emph{coherent} (as a  $\mc{B}$-module) if
there is a coherent $\mc{O}_X$-module $\mc{M}_0$  and a surjective map
 $\mc{M}_0 \otimes \mc{B} \ra \mc{M}$ of ungraded
$\mc{B}$-modules.  Left $\mc{B}$-modules are
defined similarly and  the bimodule algebra $\mc{B}$ is \emph{right (left)
noetherian} if every   right (left) ideal of $\mc{B}$ is coherent.
For the algebras that interest us, a more natural definition of
coherence will be  given in Lemma~\ref{coherent}.
\end{definition}

One can   give  various different bimodule structures to a graded
right $\mc{B}$-module $\mc{M}=\bigoplus \mc{M}_i$ and
   it will cause no loss of generality to
assume that all right $\mc{B}$-modules  have the form
\begin{equation}\label{nice-modules}
\mc{M} = \bigoplus_{n
\in \mb{Z}} {}_1 (\mc{G}_n)_{\sigma^n}
\text{ for some (left) sheaves }\mc{G}_n\in\mc{O}_X\Mod.
\end{equation}
The advantage of this choice is that the $\mc{B}$-module structure on $\mc{M}$
is given by a family of $\mc{O}_X$-module  maps
$\mc{G}_n \otimes \mc{E}_m^{\sigma^n} \to \mc{G}_{n+m}$,
 again satisfying the appropriate associativity conditions.

Graded right $\mc{B}$-modules form an abelian category $\rGr
\mc{B}$, with homomorphisms graded  of degree zero.  Its subcategory
of coherent modules is denoted $\rgr \mc{B}$. A graded
$\mc{B}$-module $\mc{M}=\bigoplus \mc{M}_i$ is \emph{right bounded}
if $\mc{M}_i = 0$ for all $i \gg 0$ and
\emph{bounded}\label{bounded-defn} if $\mc{M}_i = 0$ for all but
finitely many $i$.  A module $\mc{M}
\in \rGr\mc{B}$ is called \emph{torsion} if every coherent submodule
of $\mc{M}$ is bounded. Let $\rTors \mc{B}$ denote  the full
subcategory of $\rGr \mc{B}$ consisting of torsion modules, and
write $\rQgr \mc{B}$ for the quotient category $\rGr
\mc{B}/\negthinspace\rTors \mc{B}$.\label{qgr-defn}
 The analogous quotient category of $\rgr \mc{B}$
 will be denoted  $\rqgr \mc{B}$. The
 corresponding  categories of left modules will be  denoted by $\mc{B}\Gr$, etc.
Similar   definitions and notation  apply to modules over a graded ring $A=\bigoplus_{i\geq 0}A_i$
and we denote the natural quotient maps by $\pi_{\mc{B}} : \rGr \mc{B}\to \rQgr\mc{B}$ and
$\pi_A : \rGr A \to \rQgr A.$  Both maps are written  as $\pi$ if no confusion is possible.

 A basic technique for us will be to
 pass between a bimodule algebra and its ring of sections,
 and the next theorem, due to Van den Bergh,  gives one situation in which this is possible.

\begin{definition}
\label{ample def}
Suppose that  $\{ \mc{J}_n \}_{n \in \mb{N}}$  is a sequence of $\mc{O}_X$-bimodules.
  Then   the sequence is \emph{ample} (or, more formally, \emph{right ample})
   if the following conditions hold  for any $\mc{M}  \in \OO_X\catmod$:
\begin{enumerate}
\item $\mc{M} \otimes_{\mc{O}_X} \mc{J}_n$
is generated by global sections for $n \gg 0$.
\item $\coH^i(X, \mc{M} \otimes \mc{J}_n) = 0$ for all $i > 0$ and $n \gg 0$.
\end{enumerate}
\end{definition}

\begin{theorem}   \label{VdB main theorem}
Let $\mc{B}=\bigoplus \mc{B}_i$ be a right noetherian graded  $(\mc{O}_X,\sigma)$-algebra.
Assume that     $\{ \mc{B}_n \}_{n \in \mb{N}}$ is an
  ample sequence of $\mc{O}_X$-bimodules such that each $\mc{B}_n$
is contained in a locally free left $\mc{O}_X$-module. Then:
\begin{enumerate}
\item  The  section algebra $B = \mathrm{H}^0(X,\, \mc{B})$ is right noetherian, and
there is an equivalence of categories $\xi : \rqgr \mc{B} \simeq \rqgr B$
via the   inverse equivalences $ \mathrm{H}^0(X,\, -)$ and $- \otimes_B \mc{B}$.

\item  If $\mc{M}\in \rgr \mc{B}$ then $\mathrm{H}^0(X, \mc{M})$ is a noetherian $B$-module.
\end{enumerate}
\end{theorem}

\begin{proof}    (1) This is essentially  \cite[Theorem~5.2]{VB1};
see   \cite[Theorem~2.12]{KRS} for more details.

 (2) This is   Step~3 in the proof of  \cite[Theorem~5.2]{VB1}.
 \end{proof}

An important special case of Definition~\ref{ample def} and Theorem~\ref{VdB
main theorem} occurs when $\mc{J}_n = \mc{B}_n = (_1\mc{L}_\sigma)^{\otimes n}$
for an invertible sheaf $\mc{L}$ on $X$. We will usually write
$\mc{L}_\sigma^{\otimes n}$ for
$(_1\mc{L}_\sigma)^{\otimes n}$.
It is customary to say that
\begin{equation}\label{sigma-ample}
\mc{L} \text{ is } \sigma\text{\emph{-ample if }}\ \{
\mc{L}_{\sigma}^{\otimes n} \}_{n \geq 0}\ \text{ is an   ample
sequence of bimodules.}
\end{equation}
We will always write the corresponding bimodule algebra $ \bigoplus \mc{L}_\sigma^{\otimes n}$
as   $\mc{B}=\mc{B}(X,\mc{L},\sigma) $
with    section algebra $B=B(X,{\mathcal L},\sigma) =
 \bigoplus_{n\geq 0}
\mathrm{H}^0(X,\, \mc{L}_{\sigma}^{\otimes n}).$ This is an
equivalent definition of  the  \emph{twisted homogeneous coordinate
ring} of $X$ from \cite{AV}.
  For more detailed   results about $B(X,{\mathcal L},\sigma)$, see \cite{AV,Ke1}.

We now turn to a second special case of bimodule algebras; that of \naive\ blowups.
For this we need the  following assumptions,  which will also be fixed for the rest of  the section.

\begin{assumptions}\label{global-convention}
Fix an integral projective scheme $X$. Fix $\sigma\in\Aut(X)$, an
invertible sheaf $\mc{L}$ on $X$ and let $\mc{I}=\mc{I}_Z$ denote
the sheaf of ideals defining a subscheme $Z$ that is either
zero-dimensional or empty.  Let $S = \supp Z$. We always assume that
each $p \in S$ has infinite order under $\sigma$. Our convention on
automorphisms from the beginning of the section means that
$\mc{I}^{\sigma^i} = \mc{I}_{\sigma^{-i}(Z)}$, and so $\supp
\mc{O}_X/ \mc{I}^{\sigma^i}=\sigma^{-i}(S)$.
\end{assumptions}

Mimicking classical blowing up we set
$$\mc{I}_n =\mc{I} \mc{I}^{\sigma} \dots
\mc{I}^{\sigma^{n-1}},\quad \mc{L}_n = \mc{L}\otimes \mc{L}^{\sigma}
\otimes\dots\otimes \mc{L}^{\sigma^{n-1}}\quad\text{and}\quad
\mc{R}_n={}_1(\mc{I}_n\otimes \mc{L}_n)_{\sigma^n},$$ where all
tensor products are over $\OO_X$.  From this data we define a
bimodule algebra
 $\mc{R} = \mc{R}(X,Z,\mc{L},\sigma) = \bigoplus_{n =
0}^{\infty} \mc{R}_n$ with  \emph{\nsr} $R = R(X,Z,\mc{L},\sigma) =
\mathrm{H}^0(X,\, \mc{R}) =\bigoplus_{n\geq 0} \mathrm{H}^0(X,\,
\mc{R}_n).$
These algebras $\mc{R}$ and $R$ are called  \emph{nontrivial}
 if $Z \neq \emptyset$.
By \cite[Lemma~2.9]{KRS},
 \begin{equation}\label{rees opposite}
\mc{R}(X,Z,\mc{L},\sigma)^\op \cong
 \mc{R}(X, \sigma(Z), \mc{L}^{\sigma^{-1}}, \sigma^{-1}),
 \end{equation}
where $\mc{R}^\op$ denotes the opposite bimodule algebra in the obvious
 sense (see  \cite[Definition~2.8]{KRS} for the formal definition).
Thus   any result proved  on the right can immediately be transferred to the left.

We note that, in distinction to the situation in \cite{KRS} where
$Z$ is a single reduced point,   $\mc{I}_n = \mc{I}
\mc{I}^{\sigma} \dots \mc{I}^{\sigma^{n-1}}$ need not equal
  $\mc{I} \otimes \mc{I}^{\sigma} \otimes \dots \otimes
\mc{I}^{\sigma^{n-1}}$; indeed, in general the latter sheaf need not even be
an ideal sheaf.  This also means that the multiplication map
$\mc{R}_m \otimes \mc{R}_n \to \mc{R}_{m+n}$ is not an isomorphism
of sheaves.  In generalizing many of the results in \cite{KRS} this
is merely an annoyance, but in Sections~\ref{gen naive blowups}
and~\ref{chi}    it will make a significant difference to the results themselves.

 \begin{proposition}
\label{L no matter} Given  invertible sheaves  $\mc{L}$ and
$\mc{L}'$, then  $\rGr\mc{R}(X,Z,\mc{L},\sigma)\sim \rGr\mc{R}(X,Z,\mc{L}',\sigma)$.
\end{proposition}

\begin{proof}  The proof of  \cite[Proposition~3.5]{KRS}   can be used without change.
\end{proof}

\begin{lemma}\label{coherent}
A  module
  $\mc{M}=\bigoplus_{n\in \mathbb Z} \mc{M}_n\in \rGr \mc{R}$ is coherent if and only if
the following conditions hold:
\begin{enumerate}
\item Each   $\mc{M}_n $ is a coherent $\OO_X$-module, with  $\mc{M}_n = 0$
 for $n \ll 0$.
\item     The natural map
 $\mu_n:\mc{M}_n\otimes \mc{R}_1\to \mc{M}_{n+1}$
is surjective for   $n\gg 0$.
\end{enumerate}\end{lemma}

\begin{proof}  The proof of \cite[Lemma~3.9]{KRS}  can be used essentially unchanged.
(The only difference is that the morphism $\phi_1$ in the
commutative diagram on \cite[p.~504]{KRS} will now  be a surjection
rather than an isomorphism, but this does not affect the proof.)
\end{proof}

The next result, which  is the main result of this section,
determines when the bimodule algebra $\mc{R}$ is noetherian.

\begin{proposition}
\label{bimod alg noeth} Keep the hypotheses of
\eqref{global-convention}.  Then:
\begin{enumerate}
\item  The bimodule algebra $\mc{R} =
\mc{R}(X,Z,\mc{L},\sigma)$ is right noetherian if and only if $Z$ is
right saturating
and noetherian if and only if $Z$ is saturating.
\item If $\mc{R}$ is not right noetherian, there exists an
infinite ascending chain of coherent right ideals of $\mc{R}$ with
non-torsion factors.
\end{enumerate}
\end{proposition}

\begin{proof}  (1)
By Proposition~\ref{L no matter}, the result is independent of the choice
of $\mc{L}$ and we choose  $\mc{L} = \mc{O}_X$.  We start by  assuming  that $Z$
is right saturating.
An arbitrary right ideal $\mc{G}$ of $\mc{R}$ is given by a sequence
of bimodules $\mc{G}_i = {}_1 (\mc{H}_i)_{\sigma^i} \subseteq
\mc{R}_i$, where $\mc{H}_i\subseteq \mc{I}_i$ is an ideal sheaf,
  such that the multiplication maps  $\mc{R}_i \otimes
\mc{R}_1 \to \mc{R}_{i+1}$ restrict   to maps $\mu_i: \mc{G}_i \otimes \mc{R}_1 \to
\mc{G}_{i+1}$ for all $i \geq 0$.  By Lemma~\ref{coherent}, $\mc{G}$
is  coherent if and only if $\mu_i$ is
surjective for all $i \gg 0$. Note that the image of $\mu_i$ is $_1
(\mc{H}_i \mc{I}^{\sigma^i})_{\sigma^{i+1}}$. Equivalently,
 we are given ideal sheaves
  \begin{equation}\label{key-prop}
\mc{H}_{i} \mc{I}^{\sigma^{i}} \subseteq
 \mc{H}_{i+1}\subseteq \mc{I}\mc{I}^{\sigma}\cdots\mc{I}^{\sigma^{i}}
\qquad \text{ for all } i \geq 0,
\end{equation}
 and $\mc{G}$ is  coherent if and only if
  $\mc{H}_{i} \mc{I}^{\sigma^{i}} =
\mc{H}_{i+1}$ for $i \gg 0$. To avoid trivialities, we assume that $\mc{G}\not=0$.

Pick $r$ such that $\mc{H}_r\not=0$.
 Since $Z$ is right
saturating, $\supp \mc{O}_X/\mc{H}_r$ contains at most finitely
many points from $\{ s_i=\sigma^{-i}(s) | s\in S,  i \geq 0 \}$.
  Let $d$ be the largest   integer such that
$\sigma^d(s) = t$ for some $s, t \in S$ and  pick   $m \geq§ r
+d$ such that $(\supp \mc{O}_X/\mc{H}_r) \bigcap \sigma^{-j}(S) =
\emptyset$ for all $j \geq m$.
Put $T = \bigcup_{i = m-d}^{m-1} \sigma^{-i}(S)$,
 and let $Y$ be the open subscheme $X \smallsetminus T$  of $X$.

Set $\mc{U} = \mc{H}_{m-d} \vert_Y$, and for $n \geq m$ set
$\mc{W}_n = \mc{H}_n \vert_Y$, as well as  $\mc{J}_n =
\mc{I}^{\sigma^n} \vert_Y$ and  $\mc{V}_n = \prod_{i = m}^{n-1}
\mc{I}^{\sigma^i} \vert_Y$ (with the convention that  the product
of an empty set of
ideal sheaves equals $\mc{O}_Y$). By \eqref{key-prop}  and  the
choice of $m$, the ideal sheaves $\mc{H}_{m-d}$ and $\prod_{i =
m}^{n-1} \mc{I}^{\sigma^i}$ are comaximal in $\mc{O}_X$, so
certainly  $\mc{U}$ and $\mc{V}_n$ are comaximal in $\mc{O}_Y$. By
induction, \eqref{key-prop} shows that $\mc{H}_{m-d}
(\prod_{i=m}^{n-1} \mc{I}^{\sigma^i})(\prod_{i = m-d}^{m-1}
\mc{I}^{\sigma^i}) \subseteq \mc{H}_n  \subseteq \prod_{i = m}^{n-1}
\mc{I}^{\sigma^i}.$ Restricting   to $Y$, we get $\mc{U} \cap
\mc{V}_n = \mc{U} \mc{V}_n \subseteq \mc{W}_n \subseteq \mc{V}_n$.
Thus  \cite[Lemma~3.8]{KRS}
 implies that
$\mc{Z}_n=\mc{U}+\mc{W}_n$ is maximal among ideal sheaves $\mc{Z}$ satisfying
 $\mc{Z}  \mc{V}_n \subseteq \mc{W}_n$.
 Since
$\mc{Z}_n\mc{V}_{n+1}=\mc{Z}_n\mc{V}_n\mc{J}_n \subseteq
\mc{W}_n\mc{J}_n \subseteq \mc{W}_{n+1},$ this implies that
$\mc{Z}_n\subseteq \mc{Z}_{n+1}$  for all $n\geq m$. Thus we may
pick $n_0\geq m$ such that $\mc{Z}_n=\mc{Z}_{n+1}$ for all $n\geq
n_0$. For all such $n$,  \cite[Lemma~3.8]{KRS}    implies that
$\mc{W}_{n+1}=\mc{Z}_n \mc{V}_{n+1}=\mc{Z}_n \mc{V}_n \mc{J}_n
=\mc{W}_{n}\mc{J}_n.$ In other words,   $\mc{H}_{n+1} \vert_Y
= \mc{H}_n \mc{I}^{\sigma^n} \vert_Y$ for all $n \geq n_0$.

We need to extend this last  equation  to all of $X$. If  $t \in T$, then  $t
\not \in \sigma^{-n}(S)$ for $n\gg 0$. Looking locally at $t$ this means that
$(\mc{H}_n)_t = (\mc{H}_n \mc{I}^{\sigma^n})_t    \subseteq
(\mc{H}_{n+1})_t$ and hence  that $(\mc{H}_n)_t = (\mc{H}_{n+1})_t$ for $n \gg 0$.  Since
$ T$ is finite, there exists  a single integer
$n_1 \geq n_0$ such that for $n \geq n_1$ we have $(\mc{H}_{n}
\mc{I}^{\sigma^n})_x = (\mc{H}_{n+1})_x$ locally at every $x \in T$. By the last paragraph
  $\mc{H}_{n} \mc{I}^{\sigma^n} = \mc{H}_{n+1}$ for such $n$ and so
  $\mc{G}$ is coherent and $\mc{R}$ is right noetherian.

Conversely,   suppose that $\mc{I}$ is not right
saturating.   Write $\mc{I} = \mc{J}^{(1)} \mc{J}^{(2)} \dots
\mc{J}^{(r)}$ with   $S_i = \supp
\mc{O}_X/\mc{J}^{(i)}$ so that $S=\bigcup S_i$ partitions $S$
 into elements of distinct $\sigma$-orbits.
   Pick $s \in S$ such that the set $\{s_i \}_{i
\geq 0}$ is not critically dense.  By renumbering we may assume that
$s \in S_1$ and we may further assume that
  $S_1\subset \{\sigma^d(s), d\geq 0\}$.  Now choose an infinite
set $A$ of nonnegative integers such that the Zariski closure of
$\{s_i\}_{i \in A}$ is a proper closed subset $W$ of $X$.

Set $\mc{J} = \mc{J}^{(1)}$, and put $\mc{J}_n = \mc{J}
\mc{J}^{\sigma} \dots \mc{J}^{\sigma^{n-1}}$ as usual.  We claim
that $\mc{K} = \bigcap_{i \in A} \mc{J}^{\sigma^i}$ is nonzero. To see this,
let  $\mc{M}_t$ be the the ideal sheaf defining a
  closed point $t\in S_1$; thus
$\prod_{t \in S_1} (\mc{M}_t)^e \subseteq \mc{J}$ for some $e\geq 1$.
 The  closure  of $\{ \sigma^{-i}(S_1) | i \in A \}$  equals $W' = \bigcup
  \sigma^d(W) $, where the union is over the finite set
  $ \{d\in \mathbb{Z} \, |\,  \sigma^d(s) \in S_1 \} $.
Hence  $W'\not=X$ and the ideal sheaf $\mc{M}_{W'}$ defining $W'$
satisfies $(\mc{M}_{W'})^e \subseteq \mc{K}$. Thus $\mc{K}\not=0$.

Set $\mc{H}_n = \mc{K} \cap \mc{I}_n$ for $n \geq 0$, and observe
that $\mc{G} = \bigoplus \mc{G}_n = \bigoplus
{}_1(\mc{H}_n)_{\sigma^n}$ is a right ideal of $\mc{R}$. Write
$\mf{m}_{s_n}$ for the maximal ideal in the local ring $\mc{O}_{X,
s_n}$. Pick $n \in A$ and note that $s_n \not \in \Supp
\mc{O}_X/\mc{J}_n = \bigcup_{i = 0}^{n-1} \sigma^{-i}(S_1)$
 but $s_n \in \supp \mc{O}_X/\mc{J}^{\sigma^n}$ by the
choice of $s$. Equivalently, $(\mc{J}_n)_{s_n} = (\mc{O}_X)_{s_n}$
but $(\mc{J}_{n+1})_{s_n} = (\mc{J}^{\sigma^n})_{s_n}   \subseteq
\mf{m}_{s_n}$. Thus $\mc{K}_{s_n} \cap (\mc{J}_n)_{s_n} =
\mc{K}_{s_n} $ and $\mc{K}_{s_n} \cap (\mc{J}_{n+1})_{s_n} =
\mc{K}_{s_n} \cap (\mc{J}^{\sigma^n})_{s_n}=\mc{K}_{s_n}$ as $n\in A$.
By Nakayama's  lemma
\[
(\mc{H}_n \mc{I}^{\sigma^n})_{s_n} = (\mc{H}_n
\mc{J}^{\sigma^n})_{s_n} \subseteq \mc{K}_{s_n} \mf{m}_{s_n}
\subsetneq \mc{K}_{s_n} = \mc{K}_{s_n} \cap (\mc{J}_{n+1})_{s_n} =
\mc{K}_{s_n} \cap (\mc{I}_{n+1})_{s_n} = (\mc{H}_{n+1})_{s_n}.
\]
Since this happens for infinitely many $n$, the right ideal $\mc{G}$
is not coherent and  $\mc{R}$ is not  right noetherian.

By \eqref{rees opposite}, $\mc{R}^\op \cong \mc{R}(X,
\sigma(Z),\mc{L}^{\sigma^{-1}},\sigma^{-1})$  and so the analogous
result on  the left (that $\mc{R}$ is left noetherian if and only if
$\{\sigma^i(c) :  i\geq 0\}$ is critically dense for each $c\in
\supp Z$) follows from the one on the right.
 The result for noetherian algebras is then obvious. This completes the proof
 of~(1).

(2) If $\mc{R}$ is not noetherian,   let $\mc{G} = \bigoplus
\mc{G}_n$ be the non-coherent right ideal defined above. Set
$\mc{M}^n = \sum_{0 \leq i \leq n} \mc{G}_i \mc{R}$; thus
$(\mc{M}^n)_j = ((\mc{K}\cap \mc{I}_n)\mc{I}^{\sigma^n }\cdots
 \mc{I}^{\sigma^{j-1}})_{\sigma^j}$   for $j\geq n$. This gives a chain of
coherent right ideals $\mc{M}^0 \subseteq \mc{M}^1 \subseteq \dots
\subseteq \mc{R}$. As before, looking locally at a point $s_n$ for
$n \in A$ shows that $(\mc{M}^n)_j \subsetneq (\mc{M}^{n+1})_j$ for
all $n\in A$ and all $j\geq n+1$.   Thus the
 subsequence $\{ \mc{M}^n : n \in A\}$ gives the desired   chain of right ideals.
\end{proof}

%%%%%%%%%%%%%%%%%%%%%%%%%%%%%
\section{Ampleness}\label{sect-ample-new}

We maintain the hypotheses from \eqref{global-convention}.   The
main aim of this section (Theorem~\ref{general-ampleness}) is to
prove, in considerable generality, that the sequence of bimodules
$\{ \mc{R}_n = (\mc{I}_n\otimes \mc{L}_n)_{\sigma^n} \}$ is ample in
the sense of Definition~\ref{ample def}. Combined with the results
of Section~\ref{definitions}
 this will prove parts~(1) and (2) of Theorem~\ref{mainthm}.  This section differs significantly from the
proof of ampleness in \cite{KRS}, because we are proving a much
stronger result. First,  we need much stronger  estimates of the number of
points in a  $\sigma$-orbit that can be separated by the sheaf $\mc{L}_n$.
Secondly, for applications in \cite{RS} we need to prove the result for
a $\sigma$-ample sheaf $\mc{L}$, whereas  \cite{KRS} only  proved the
result when $\mc{L}$ was also very ample.

Here is our goal:
\begin{theorem}
\label{general-ampleness}  Keep Assumptions~\ref{global-convention},
 and assume in
addition that $\dim X \geq 2$, that $\mc{L}$ is $\sigma$-ample, and
that $\mc{I}$ defines a subscheme $Z$ of $X$ such that each point of
$S = \supp Z$ lies on a dense $\sigma$-orbit.  Then:
\begin{enumerate}
\item  $\{ \mc{R}_n=(\mc{I}_n\otimes\mc{L}_n)_{\sigma^n} \}$ is an ample sequence.

\item  Assume in addition that $Z$ is saturating, and set
$\mc{R}=\mc{R}(X,Z,\mc{L},\sigma)$.  Then the \nsr\  $R =
\mathrm{H}^0(X,\, \mc{R})$ is  noetherian and there is an equivalence of
categories $\xi:\rqgr \mc{R} \simeq \rqgr R$ via the inverse
equivalences $\mathrm{H}^0(X,\, -)$ and $-\otimes_R\mc{R}$. Similarly,
$\mc{R}\qgr\simeq R\qgr$.
\end{enumerate}
\end{theorem}

 We begin with a useful combinatorial notion and some preliminary results related to it.

\begin{definition}\label{sparse-defn}
 For us, the  natural numbers $\mb{N}$
contains $0$  and we write $\mb{N}_n = \{0,1, \dots, n-1 \}$.
A subset $S \subset \mb{N}$ is called \emph{sparse} if for every
 $m \in \mb{N}_+=\mb{N}\smallsetminus \{0\}$ there exists $N(m)\in \mb{N}_+$ such that
$|S \cap \mb{N}_n  | \leq n/m$ for all $n \geq N(m)$.  If $S$ is a
sparse set, then any monotonically increasing
 function $N: \mb{N}_+ \to \mb{N}_+$
 satisfying this condition is called a \emph{\bound\ function}  for $S$.
\end{definition}

\begin{lemma}
\label{sparselemma1} Let $T_1, T_2, \dots T_d \subset \mb{N}$
be sparse sets with respective \bound\ functions $N_1, N_2, \dots,
N_d$. Then $T = \bigcup_{i = 1}^d T_i$ is sparse, with one  \bound\
function $N$ given by $N(m) = \underset{1 \leq i \leq d}{\max}
N_i(md)$.
\end{lemma}
\begin{proof}  If  $n \geq \underset{1 \leq i \leq d}{\max} N_i(md)$, then
$
|T \cap \mb{N}_n | \leq \sum_{i =1}^d | T_i \cap \mb{N}_n | \leq
\sum_{i = 1}^d n/(md) = n/m$.
\end{proof}

For any   $S \subset \mb{N}$ and   $i
\geq 1$, define $S + i = \{ s + i | s \in S \} \subset \mb{N}$ and
 set $S_d = \bigcup_{i = 1}^d ((S + i) \cap
S)$ for $d\geq 1$.  In other words, $S_d$ consists of   those
numbers $s\in S$  such that some integer $s'\in [s-d, s-1]$  also lies in  $S$.

\begin{lemma}
\label{sparselemma2} Let  $S \subset \mb{N}$ and suppose that, for all $d\geq 1$, the sets
$S_d \subset \mb{N}$ are sparse, with respective \bound\ functions $N_d$.  Then $S$ is also
sparse, with one \bound\ function  being
 $N(m) = \max(3m, N_{3m}(3m))$.
\end{lemma}

\begin{proof} Write $T^c =\mb{N}\smallsetminus T$ for the complement of
a set $T\subseteq \mb{N}$.   Fix   $d \geq 1$ and
suppose that $s \in S \smallsetminus S_d$. Then as long as $s \geq
d$, we know that $\{ s-1, s-2, \dots, s-d \} \subset S^{c}$. If
$t\not=s \in S \smallsetminus S_d$ is another natural number with $t
\geq d$,
 then $|s-t| > d$ and so  $\{ t -1, t-2, \dots, t-d \}
\subset S^c$, with
$$\{ s-1, s-2, \dots, s-d \} \cap \{ t -1, t-2,
\dots, t-d \} = \emptyset.$$
 Also, there is at most one   $u \in S
\smallsetminus S_d$ with $0 \leq u < d$. Combining these
observations, we conclude that
\begin{equation}
\label{main estimate} |S^c \cap \mb{N}_n | \geq d \big( \big| (S
 \smallsetminus S_d) \cap \mb{N}_n \big| - 1 \big)  \qquad {\rm for \ each \ } n\geq 1.
\end{equation}

Now use the formul\ae\ $ | S \cap \mb{N}_n | + | S^c \cap \mb{N}_n |
= n $ and $ | S \cap \mb{N}_n | = | (S \smallsetminus S_d) \cap
\mb{N}_n | + | S_d \cap \mb{N}_n| $  to transform
  \eqref{main estimate} into
 \begin{equation}
\label{main estimate 2} | S \cap \mb{N}_n | \leq \frac{n + d + d\, |
S_d \cap \mb{N}_n | }{d+1}\qquad {\rm for \ each \ } n\geq 1.
\end{equation}

For a given $m > 0$, take $d = 3m$ in this calculation.
If $n \geq \max(3m, N_{3m}(3m))$, then
\[
|S \cap \mb{N}_n | \leq \frac{n + 3m + 3m\, | S_{3m} \cap \mb{N}_n
|}{3m+1} \leq \frac{n + 3m + 3m(n/3m)}{3m+1}
\leq \frac{3n}{3m+1} \leq \frac{n}{m}.
\]
So $N(m) =\max(3m, N_{3m}(3m))$ defines
a \bound\ function for $S$ and  hence $S$ is sparse.
\end{proof}

Define a {\it reduced $d$-cycle $W$ on $X$}\label{cycle-defn} to be
a formal sum $W=\sum W_i$, where the $W_i$ are distinct integral
closed subschemes of $X$ of dimension $d$. Equivalently, we may
identify $W$ with the reduced and equidimensional closed subscheme
$\bigcup W_i$.   In the proof of the next proposition we
  will want to induct on the degrees of cycles on $X$. To do this, fix
some closed  immersion $\iota: X\hookrightarrow \mb{P}^s$ and define 
the degree of a  reduced $d$-cycle  $W$   in
  $X$ to be the intersection number  $i(H^{\cdot d }\cdot W; \mb{P}^s)
 = i(\overbrace{H\cdot H\cdots H}^d\cdot W; \mb{P}^s)$, where $H$ is a
hyperplane of $\mb{P}^s$. This is, of course, dependent on the given embedding.

We would like to thank Mark Gross for his helpful suggestions concerning the next result.

\begin{lemma} \label{int theory facts}
Let $X$ be an integral projective scheme of
dimension $\geq 2$, with $\sigma \in \mathrm{Aut}(X)$ and a fixed projective embedding $\iota: X \to
\mb{P}^s$ which will be used to measure degrees. Then
there exists  $M \geq 2$  such that $\deg
\sigma(Z) \leq M \deg Z$ for all integral closed subschemes  $Z\subseteq X$.
\end{lemma}

\begin{proof}
Let $\dim Z=d$.   It will  be convenient to interpret degrees  in terms of
 the intersection theory  on $X$,  as developed in
\cite{Kl}. Thus, write  $\mc{N} = \iota^* \mc{O}_{\mb{P}^s}(1)$
 and let $(\mc{N}^{\cdot d}\cdot
\sigma(Z))_X$ denote the intersection number from \cite{Kl}; by
\cite[Proposition~5, p.~298 and
 Corollary~3, p.~301]{Kl} this equals
 $i(H^{\cdot d}\cdot \sigma(Z); \mb{P}^s)$. Next, by Bertini's Theorem
  \cite[Th\'eor\`eme~I.6.10(2,3)]{Jou},
 we may choose generic hyperplanes $H_i$ on $\mb{P}^s$,  so that
   $\sigma(Z)\cap H_1\cap\cdots \cap H_d$ is a reduced set of points and each
    $E_i=H_i\cap X$ is a reduced and irreducible subscheme of $X$.

    By \cite[Remark~1, p.~301]{Kl}  applied to $\mb{P}^s$,
 $\deg\sigma(Z) = \#(\sigma(Z)\cap H_1\cap\cdots \cap H_d)$, the number of points in
 this intersection,  and hence $\deg\sigma(Z)=\#(\sigma(Z)\cap E_1\cap\cdots \cap E_d)$.
This is also the number of points in $(Z\cap
\sigma^{-1}(E_1)\cap\cdots \cap \sigma^{-1}(E_d))$ and, of course,
each $\sigma^{-1}(E_i)$ is a reduced and irreducible subscheme of
$X$ and hence of $\mb{P}^s$. Thus, by Bezout's Theorem  \cite[Example~8.4.6]{Fu},
$\deg\sigma(Z)\leq \deg Z \left(\prod_{i=1}^d \deg
\sigma^{-1}(E_i)\right)$. Finally, by \cite[Remark~2, p. 301]{Kl},
the degree of the divisor $\sigma^{-1}(E_i)$  regarded as a subscheme of $X$ equals its
degree regarded as an element of the complete linear system
$|\mc{O}_X( \sigma^{-1}(E_i)) | = |\sigma^*(\mc{N})|$. In other
words, $\deg\sigma(Z)\leq M\deg Z$, where
 $M=\left(\deg \sigma^*(\mc{N}) \right)^d$.
 \end{proof}

Let  $p_1, \dots p_n$ be closed points  on a projective scheme $X$. Then an invertible sheaf $\mc{L}$
\emph{separates}\label{separates-defn}   $\{ p_1, \dots, p_n \}$ if, for each
$1 \leq i \leq n$, there is a section  $s_i \in \HB^0(X,\mc{L})$ such
that $s_i(p_i) \neq 0$ but $s_i(p_j) = 0$ for  $j \neq i$.

 In Proposition~\ref{separating orbit points}   we will    use intersection theory
  to   estimate the number of  points in a
$\sigma$-orbit  on $X$ that can be separated by
  $\mc{L}_n$, in the notation of Assumptions~\ref{global-convention}.
  The main idea, which is provided by the next proposition,
 is to study hyperplane sections $W = H \cap X$ of a closed immersion
 $\iota: X\hookrightarrow \mb{P}^s$, and study the size of the
sets $W \cap \{ \sigma^i(x) | i\geq 0 \}$ for $x \in X$ lying on a
dense $\sigma$-orbit.  These sets need not be finite, so
we actually show that $\{ i \geq 0 | \sigma^i(x) \in W \}$ is
sparse, and estimate its \bound\ function.

  \begin{proposition} \label{universal bound functions}
Let $X$ be an integral projective scheme of dimension $\geq 2$ with
automorphism $\alpha$ and an immersion $\iota: X \to \mb{P}^s$ that
will be used to measure degrees.  Then, for any   $e \geq 1$ and
$0\leq d<\dim X$,  there exists a \bound\ function $N_{d,e}(m)$,
depending   on ($X, \alpha, \iota$),   with the following property:
\begin{enumerate}
\item[] For any $x \in X$ for which $P = \{ \alpha^i(x) | i \in \mb{Z} \}$
 is dense in $X$, and any   reduced
$d$-cycle $Z$ on $X$  with $\deg Z \leq e$, the set $S = \{i \geq 0 |
\alpha^i(x) \in Z \}$ is sparse with \bound\ function
$N_{d,e}(m)$.
\end{enumerate}
\end{proposition}

\begin{proof} If  $Z$ is a reduced $0$-cycle with $\deg Z\leq e$,
then $\deg Z$ is just the number of points in $Z$ and so $\#S\leq e$ for any $x\in X$.
 It therefore suffices to  take $N_{0,e} = e$.
Now, for some $0<d<\dim X$,  suppose by induction that   $N_{c,f}$
has been constructed  for every $0 \leq c<d$ and $f \geq 1$.  By
taking supremums, we can assume that $N_{c,f}=N_{b,f}$ for all
$0\leq b\leq c<d$ and that $N_{c,f}\leq N_{c,f+1}$ for all $f$.

Fix  $x \in X$ such that $P = \{ \alpha^i(x) | i \in \mb{Z} \}$  is dense, let $Y$ is an irreducible
$d$-dimensional subvariety of $X$ with
 $\deg Y \leq e$, and define $S =S^Y= \{i \geq 0 | \alpha^i(x) \in Y
\}$.  Fix some $a \geq 1$ and  define the set $S_a$ as before
Lemma~\ref{sparselemma2}. Since $S+j = \{\ell \geq j :
\alpha^{\ell}(x) \in \alpha^j(Y)\}$,  clearly
$ S_a   \; \subseteq\;   \left\{i \geq 0 \, | \, \alpha^i(x) \in \bigcup_{k =1}^a ( Y \cap
\alpha^k(Y) ) \right\}. $
For each $k$, write  $ ( Y \cap \alpha^k(Y) ) = C_{k1} \cup C_{k2}
\cup \dots \cup C_{k,r_k}$ where   the $C_{kj}$ are distinct irreducible
components.

Suppose first  that $Y = \alpha^k(Y)$ for some $k \geq 1$. In this case, if
 $S\not=\emptyset$, then certainly  $Y\cap P \not=\emptyset $ and so   $P$ is entirely
contained in $\bigcup_{i = 1}^k \alpha^i(Y)$. This contradicts the
density of $P$. Therefore, $S = \emptyset$, which is trivially
sparse   and any \bound\ function whatsoever will do.

  Thus   we may assume that $Y \neq \alpha^k(Y)$ for   $k \geq 1$.  In
particular, since $Y$ is irreducible, $\dim C_{ij} < \dim Y$ for all
$i,j$.  Let $M\geq 2$ denote  the constant  from  Lemma~\ref{int theory facts}.
 Then Bezout's Theorem \cite[Example~8.4.6]{Fu}   implies that, for any $i\leq r_k$,
 \begin{equation}\label{bezout}
\deg C_{ki} \leq \sum_{j = 1}^{r_k} \deg C_{kj} \leq
  (\deg Y)(\deg \alpha^k(Y)) \leq  (\deg Y)(M^k\deg Y)  \leq M^{a}
  e^2.
  \end{equation}
For each $k,j$  set $T_{kj} = \{i \geq 0 |
\alpha^i(x) \in C_{kj} \}$.  Now apply induction, recalling our
assumptions on the $N_{c,f}$ for $c<d$. This implies that  each
$T_{kj}$ is sparse and that we can use $N_{c,b(a)}(m)$, where
$c=d-1$ and
 $b(a)=M^ae^2$, as its   \bound\
function.
By definition,   $\deg C_{kj}\geq 1$ for each $k,j$ and so \eqref{bezout} also implies that
  $r_k\leq b(a)$. Therefore, as
   $S_a = \bigcup T_{kj}$,
Lemma~\ref{sparselemma1} implies  that $S_a$ is a sparse set with
\bound\ function $\widetilde{N}_a(m) = N_{c,b(a)}(ab(a)m)$.   By
Lemma~\ref{sparselemma2}, $S=S^Y$ is a sparse set with \bound\
function $\widetilde{N}(m) =\max\{3m, \widetilde{N}_{3m}(3m)\}
=\max\{3m, N_{c,b(3m)}(9m^2 b(3m)) \}$.   Notice that this \bound\ function is independent both
of $x$ and of $Y$ (and also works for an irreducible subscheme $Y$
satisfying $Y=\alpha^k(Y)$).

Finally,  let $Z$ be an arbitrary reduced $d$-cycle with $\deg Z
\leq e$ and set  $S' = \{ i \geq 0 | \alpha^i(x) \in Z \}$. Write $Z
= Y_1\cup \cdots \cup Y_r$, where the $Y_j$ are  distinct
irreducible components of dimension $d$.  By the multilinearity of
the intersection form,   $\deg Z = \sum_{i=1}^r \deg Y_i$ and so $1\leq \deg
Y_j\leq e$ for each $j$ and  hence $r\leq e$. For each $j$, set $U_j
=S^{Y_j}= \{i \geq 0 \, | \, \alpha^i(x) \in Y_j \}$. By the
conclusion of the last paragraph, each $U_j$ is sparse with \bound\
function $\widetilde{N}(m) $. Since $S' = \bigcup_{j=1}^r U_j$ with
$r\leq e$, Lemma~\ref{sparselemma1} implies that $S'$ is sparse with
\bound\ function $N(m) = \widetilde{N}(em).$ This depends only on
$e$ and the previously constructed \bound\ functions, and all
constructions are independent of the choice of $x$.  Thus we may
take $N_{d,e} = N$, completing the induction.   \end{proof}

The following easy  fact was   observed in the proof of \cite[Lemma~4.4]{KRS}.

\begin{lemma}\label{separate} Suppose that an invertible sheaf
$\mc{L}$  separates the closed points $\{ p_1, \dots, p_n \}$ on a
projective scheme $X$. If $\mc{N}$ is a very ample  invertible sheaf
 then $\mc{N}\otimes \mc{L}$ also separates
$\{ p_1, \dots, p_n \}$.\qed
\end{lemma}

We next  use the estimates from Proposition~\ref{universal bound functions}
 to get good bounds on the separation of points. The idea behind the proof is
 similar to that of \cite[Proposition~4.6]{KRS}, although for the applications in this paper
 we need  much
 more efficient bounds than those given in \cite{KRS}.

\begin{proposition}
\label{separating orbit points} Let $X$ be an integral
projective scheme with $\dim X \geq 2$, $\sigma \in \aut(X)$, and a
$\sigma$-ample invertible sheaf $\mc{L}$.
Then, for any real number $\delta >0$, there exists $M = M(\delta) \geq 0$
with the following property:
\begin{enumerate}
\item[]  Let  $x\in X$ be a closed point,
such that     $\{\sigma^i(x) : i \in \mb{Z}\}$ is dense and write
   $P_n =  \{ \sigma^{-i}(x): 0 \leq i \leq n-1\}$ for each $n \geq 1$. Then, for all $n \geq M$ and
all $p \in \mb{Z}$, the sheaf $\mc{L}_{\lfloor \delta n
\rfloor}^{\sigma^p }$ separates $P_n$.
\end{enumerate}
\end{proposition}

\begin{proof}
Fix $\delta > 0$. As observed at the beginning of the proof of \cite[Lemma~4.4]{KRS},
  $\mc{L}_{\lfloor \delta n
\rfloor}^{\sigma^p }$ separates $P_n$ if and only if
$\mc{L}_{\lfloor \delta n \rfloor}$ separates $\sigma^{p}(P_n)=
 \{ \sigma^{-i}(y): 0 \leq i \leq n-1\}$, where $y=\sigma^{p}(x)$.
Since  $M$ will be  chosen independently  of the point $x$, it
therefore suffices to prove the result when $p=0$. By
\cite[Proposition~3.2]{AV}, $\mc{L}_r$ is very ample for all $r \gg
0$.  Fix some such
  $\mc{L}_r$, set
 $d +1 = \dim_k \HB^0(X,\mc{L}_r)$   and  write
$\tau: X \hookrightarrow  \mb{P}^d$ for the closed immersion
associated to a  fixed  basis of $\HB^0(X,\mc{L}_r)$. More generally, for
any $j \in \mb{Z}$ write $\tau_j=\tau\sigma^j : X \hookrightarrow
\mb{P}^d$, which is the closed immersion associated to an appropriate
basis of $\HB^0(X,\mc{L}_r^{\sigma^j})$.

Given a closed point $y\in X$ such that
  $\{\sigma^i(y) : i \in \mb{Z}\}$ is dense, a hyperplane $H$ of $\mb{P}^d$  and $j\in \mb{Z}$, set
$S(y,H,j)=\{i\geq 0  : \sigma^{-i}(y)\in \tau^{-1}_j(H)\}$. For such
a hyperplane  $H$ it follows, for example from Bezout's Theorem
that $\tau^{-1}(H) = H \cap X$, thought of as a reduced $(\dim X -
1)$-cycle,  satisfies $\deg \tau^{-1}(H)\leq \deg X$. Thus, by
Proposition~\ref{universal bound functions} (applied with
$\alpha=\sigma^{-1}$), there is a fixed \bound\ function $N(m)$ such
that,
 independently of $H$ and $y$, the set
$S(y,H,0)$  is sparse with \bound\ function $N$.  Notice that $\sigma^{j-i}(y)\in \tau^{-1}(H) $
 if and only if $\sigma^{-i}(y)\in \sigma^{-j}(\tau^{-1}(H)) = \tau^{-1}_j(H)$ and so
 $S(y,H,j)=S(\sigma^{j}(y),H,0)$.
Since $N(m)$ was chosen independently of $y$, we conclude that  each
$S(y,H,j)$ is sparse with \bound\   function $N(m)$.
   To rephrase, fix $m > 0$ and set $M_1 = N(m)$. Then  for
all $n \geq M_1$ and $j \in \mb{Z}$, at most $n/m$ of the points in
the set $P_n$ lie on a single  hyperplane section $\tau_j^{-1}(H)$.
We emphasize that the choice of $M_1$  is independent of $y$.

A set of closed points $\{y_1, \dots, y_t\}\subset  \mb{P}^d$ is
called \emph{linearly general} if the smallest linear subspace
containing the points has dimension $t-1$.  Suppose $P$ is a set of
closed points in $\mb{P}^d$ such that $|P| \geq d+1$, and such that
$P$ is not entirely contained in any hyperplane of $\mb{P}^d$. Then
given $y \in P$, an easy inductive argument shows that we can find
$d$ other points $y_1, \dots, y_d \in P$ such that $\{ y, y_1,
\dots, y_d \}$ is linearly general in $\mb{P}^d$.  Now choose $n
\geq M_1$ and pick any $z \in P_n$. Suppose first  that $n
> n/m$.  By the previous paragraph  we can
pick a set of points $T_0 \subseteq P_n \smallsetminus \{ z\}$ with
$|T_0| = d$ such that   $\tau(T_0 \cup \{ z\})$ is
linearly general in $\mb{P}^d$.  If  $n - d > n/m$,   repeat
the process to
 find   $T_1 \subseteq P_n \smallsetminus (\{z\} \cup T_0)$
with $|T_1| = d$ such that   $\tau_r(\{z \} \cup T_1)$
is linearly general in $\mb{P}^d$.  Continue this process
inductively as long as possible.  This partitions $P_n$ into
disjoint subsets
\[
P_n = \{ z \} \cup T_0 \cup \dots \cup T_{k-1} \cup V
\]
where $V$ contains $q\leq n/m$ elements and,  for each $i$,  $|T_i| = d$
  and $\tau_{ir}(\{z\} \cup T_i)$  is linearly general in $\mb{P}^d$.

Fix $0\leq j\leq k-1$.  Since $\tau_{jr}(\{z \} \cup T_j)$ is
linearly general in $\mb{P}^d$,   we can find a hyperplane $H$ in
$\mb{P}^d$ such that $\tau_{jr}(T_j) \subseteq H$ but $\tau_{jr}(z)
\not \in H$. Since the
morphism $\tau_{jr} : X\hookrightarrow \mb{P}^d$ is defined via a
basis of $\HB^0(X,\mc{L}_r^{\sigma^{jr}})$ this is equivalent to the
existence of  a section  $s_j \in \HB^0(X,\mc{L}_r^{\sigma^{jr}})$ with
$s_j(z) \neq 0$ but $s_j(y) = 0$ for all $y \in T_j$.  Now
$\mc{L}_r^{\sigma^{(k+i)r}}$ is very ample for all $i\geq 0$ and so  it separates any pair of
points. Thus, if  $V = \{v_0, \dots, v_{q-1}\}$  then,  for each $i\leq q-1$,
we can also find   $t_i \in \HB^0(X,
\mc{L}_r^{\sigma^{(k+i)r}})$  such that
$t_i(z) \neq 0$ but $t_i(v_i) = 0$. Consequently,
\[
\mc{L}_{(k+q)r} = \mc{L}_r \otimes \mc{L}_r^{\sigma^r} \otimes \dots
\otimes \mc{L}_r^{\sigma^{(k+q-1)r}}
\]
has a section $ s = s_0 \otimes \dots \otimes s_{k-1} \otimes t_0
\otimes \dots \otimes t_{q-1}$ with $s(z) \neq 0$ but $s(y) = 0$ for
all $y \in P_n \smallsetminus \{z\}$. Since $z\in P_n$ was arbitrary, we
conclude that $\mc{L}_{(k+q)r}$ separates $P_n$ for any $n \geq
M_1$.

It remains to convert this into the assertion that $\mc{L}_{\lfloor
\delta n \rfloor}$ separates $P_n$ for $n$ large.   By
\cite[Theorem~6.1]{Ke1},  $d =d(r)= \dim_k \HB^0(X,\mc{L}_r) -1$
grows at least quadratically in $r$.  So we may choose $r$ large
enough so that $r/d < \delta/3$ and then take  $m \geq 3r/\delta$
in the argument above. Since $n=kd+1+q$, this ensures that
 $(k+q)r  <nr/d+qr \leq \frac{2}{3}\delta n$.
 Set $f(n) = \lfloor \delta n \rfloor
- (k+q)r$. Since $\lim_{n \to \infty} f(n) = \infty$ and $\mc{L}$ is
$\sigma$-ample, $\mc{L}_{f(n)}$ is very ample
 for $n \gg 0$, say for $n \geq M_2$.
 Since $\mc{L}_{(k+q)r}$ separates $P_n$ for all $n\geq M_1$,  it follows
from Lemma~\ref{separate} that the sheaf $ \mc{L}_{\lfloor \delta n
\rfloor} = \mc{L}_{(k+q)r} \otimes \mc{L}_{f(n)}$ also separates
$P_n$ for all $n \geq M=M_1+M_2$.
\end{proof}

For the proof of Theorem~\ref{general-ampleness} we will need    the
following concept of Castelnuovo-Mumford regularity.
  Let $X$ be a projective scheme with a  very ample invertible sheaf $\mc{N}$.
Then a coherent sheaf $\mc{F}$  is called  \emph{$m$-regular with
respect to $\mc{N}$} if $\HB^i(X,\, \mc{F} \otimes \mc{N}^{\otimes
m-i}) = 0$ for all $  1 \leq i \leq \dim X.$ If $\mc{F}$ is
$m$-regular with respect to $\mc{N}$ then it is also $m+1$-regular
with respect to $\mc{N}$ (see \cite[Theorem~1.8.5 and
Remark~1.8.14]{Laz1}).  The \emph{regularity
$\reg_{\mc{N}}\mc{F} $ of $\mc{F}$ with respect to
$\mc{N}$}\label{regularity-defn} is then defined to be
 the minimum $m$ for which $\mc{F}$ is $m$-regular with respect to $\mc{N}$.
We will delete reference to the sheaf $\mc{N}$ if it is understood.

We need the following three results, the first of which provides a useful sheaf for
measuring regularity.

\begin{lemma} \cite[Theorem~1, p.~520]{Fj}
\label{reg of amples}  Let $X$ be a projective scheme. Then there
exists  a very ample sheaf $\mc{N}$ on $X$ such that $\HB^i(X,\,
\mc{L} \otimes \mc{N}) = 0$ for every very ample  invertible sheaf
$\mc{L}$  and integer  $i \geq 1$. \qed
\end{lemma}

\begin{lemma} \cite[Proposition~2.7]{Ke2}
\label{regularity and products} Let $X$ be a projective scheme with
very ample invertible sheaf $\mc{N}$. Then there is a constant $C$
(depending on $X$ and $\mc{N}$) with the following property:
\begin{enumerate}
\item[]  If
$\mc{F}, \mc{G} $ are coherent sheaves on $X$ such that the closed
set where both $\mc{F}$ and $\mc{G}$ fail to be locally free has
dimension $\leq 2$, then $\reg_{\mc{N}} (\mc{F} \otimes \mc{G}) \leq
\reg_{\mc{N}}  \mc{F} + \reg_{\mc{N}}  \mc{G} + C$.\qed
\end{enumerate}
\end{lemma}

\begin{corollary}  \label{reg and ample seq} Let $X$ be a projective scheme with
very ample invertible sheaf $\mc{N}$. Let $\{ \mc{F}_n \}$  be a
sequence of coherent sheaves on $X$, such that, for each $n$,  the closed set where
  $\mc{F}_n$ is not locally free has dimension at most $2$.  Then
$\{ \mc{F}_n \}$  is an ample sequence if and only if $ \lim_{n \to
\infty} \reg_{\mc{N}}  \mc{F}_n = - \infty$.
\end{corollary}

\begin{proof} Suppose that  $\lim_{n \to \infty} \reg \mc{F}_n = - \infty$ and let
$\mc{G} $ be a  coherent sheaf.  Then Lemma~\ref{regularity and
products} implies that $\reg (\mc{G} \otimes \mc{F}_n) \leq    \reg
\mc{F}_n +(\reg \mc{G}+ C)$  for some constant $C$. Thus   $\lim_{n
\to \infty} \reg (\mc{G} \otimes \mc{F}_n) = - \infty$ and so
$\HB^i(X,\,\mc{G} \otimes \mc{F}_n) = 0$ for $n \gg 0$ and all
$i>0$. By Mumford's Theorem \cite[Theorem~1.8.5 and
Remark~1.8.14]{Laz1} $\mc{G} \otimes \mc{F}_n$ is also generated by
its sections for $n\gg 0$.  Hence  $\{\mc{F}_n\}$ is ample.

Conversely, suppose that  $\{\mc{F}_n \}$ is   an ample sequence and  pick $m
\in \mb{Z}$. Then  $\HB^i(X,\, \mc{N}^{\otimes m} \otimes \mc{F}_n)
= 0$ for all $i \geq 1$ and $n \gg 0$. Therefore,  $\reg \mc{F}_n
\leq m + \dim X$ for all $n\gg 0$ which, since $m$  is arbitrary,
implies that $\lim_{n \to \infty} \reg \mc{F}_n = -\infty$.
\end{proof}

One advantage of regularity is that it gives a convenient way to
rephrase Proposition~\ref{separating orbit points}.

\begin{corollary}\label{ampleness-sublemma}
Assume the  hypotheses and notation from
Proposition~\ref{separating orbit points}, let the closed point $x\in X$
 have ideal sheaf $\mc{J}$   and pick a  very ample
invertible sheaf  $\mc{N}$ by Lemma~\ref{reg of amples}. Then, for
any real number $\delta > 0$ and function $f:  \mb{N} \to \mb{Z}$,  we
have $\reg_{\mc{N}} (\mc{L}_{\lfloor \delta n
\rfloor}^{\sigma^{f(n)}} \otimes \mc{J}_n) \leq \dim X + 1$  for all $n\gg  0$.
\end{corollary}

\begin{proof}  Note  that   $P_n=\{\sigma^{-i}(x) : 0\leq i\leq n-1\}$ has ideal
 sheaf $\mc{J}_n = \mc{J}\mc{J}^{\sigma}\cdots \mc{J}^{\sigma^{n-1}}$.
By \cite[Proposition~3.2]{AV},  the sheaf  $\mc{L}_{\lfloor \delta n
\rfloor}$, and hence the sheaf $\mc{L}_{\lfloor \delta n
\rfloor}^{\sigma^{f(n)}}$, is very ample for $n \gg 0$. Thus
Lemma~\ref{reg of amples} implies that $\HB^1(X,\, \mc{N}^{\otimes
m} \otimes \mc{L}_{\lfloor \delta n \rfloor}^{\sigma^{f(n)}}) = 0$
for all $n\gg 0$ and $m \geq 1$. Pick $M$ by
Proposition~\ref{separating orbit points}; thus, for all $p
\in\mb{Z}$ and all $n\geq M$, $\mc{L}_{\lfloor \delta n
\rfloor}^{\sigma^p}$ separates $P_n$. This holds, in particular, for
$p=f(n)$ and so $\mc{L}_{\lfloor \delta n \rfloor}^{\sigma^{f(n)}}$
separates $P_n$ for $n\geq M$.  As $\mc{N}$ is very ample,
Lemma~\ref{separate} implies that $\mc{N}^{\otimes m} \otimes
\mc{L}_{\lfloor \delta n \rfloor}^{\sigma^{f(n)}}$ also separates
$P_n$ for all $n \gg 0$. This in turn  implies  that the canonical
map
\[
\HB^0(X,\, \mc{N}^{\otimes m} \otimes \mc{L}_{\lfloor \delta n
\rfloor}^{\sigma^{f(n)}}) \lra \HB^0(X,\, \mc{N}^{\otimes m} \otimes
\mc{L}_{ \lfloor \delta n \rfloor}^{\sigma^{f(n)}} \otimes
\mc{O}_X/\mc{J}_n)
\]
is surjective  for $n\gg 0$ (see the last paragraph of the proof of
\cite[Lemma~4.4]{KRS}).  From the long exact sequence in cohomology,
we conclude that $\HB^1(X,\, \mc{N}^{\otimes m} \otimes
\mc{L}_{\lfloor \delta n \rfloor}^{\sigma^{f(n)}} \otimes \mc{J}_n)
= 0$ for $n \gg 0$ and $m\geq 1$.

The higher cohomology groups are much easier to deal with. Indeed,
 fix $r\geq 2$ and $m \geq 1$  and write
$\mc{F} = \mc{L}_{\lfloor \delta n \rfloor}^{\sigma^{f(n)}}\otimes
\mc{N}^{\otimes m}$ for
 some $n\gg 0$.
Recall that $\mc{O}_X/\mc{J}_n$ is supported at a set of dimension
$0$  and that $\mc{L}$ is $\sigma$-ample.
Thus from the exact
sequence
$$\HB^{r-1}(X,\,  \mc{F}\otimes  \mc{O}_X/\mc{J}_n) \to
\HB^r(X,\, \mc{F} \otimes \mc{J}_n) \to \HB^r(X,\,\mc{F} ) \to
\HB^r(X, \, \mc{F}\otimes \mc{O}_X/\mc{J}_n),$$ one obtains $
\HB^r(X,\,  \mc{L}_{\lfloor \delta n \rfloor}^{\sigma^{f(n)}}\otimes
\mc{N}^{\otimes m} \otimes \mc{J}_n) \cong \HB^r(X,\,
\mc{L}_{\lfloor \delta n \rfloor}^{\sigma^{f(n)}}\otimes
\mc{N}^{\otimes m})$.   But, for $n\gg 0$, the sheaf
$\mc{L}_{\lfloor \delta n \rfloor}^{\sigma^{f(n)}}$ is very ample
and so $\HB^r(X,\, \mc{L}_{\lfloor \delta n
\rfloor}^{\sigma^{f(n)}}\otimes \mc{N}^{\otimes m}) = 0$   by the
choice of $\mc{N}$.  Altogether, $\reg (\mc{L}_{\lfloor \delta n
\rfloor}^{\sigma^{f(n)}} \otimes \mc{J}_n) \leq \dim X + 1$ for $n
\gg 0$, as claimed.
\end{proof}

We are now ready to prove    Theorem~\ref{general-ampleness}.

\noindent
{\bf Proof of Theorem~\ref{general-ampleness}.}
(1)  We need to prove that the sequence $\{\mc{R}_n=\mc{L}_n\otimes\mc{I}_n\}$ is ample.
All regularities in the proof will be taken with a respect to a very
ample sheaf $\mc{N}$ satisfying the conclusion of Lemma~\ref{reg of amples}.
Pick a sequence of closed points  $s_1, s_2, \dots, s_d $  in
$S$, in general with repeats, such that if $\mc{J}^{( i )}$ is the
ideal sheaf of the reduced point $s_i$, then $ \mc{J} = \mc{J}^{(  1
)} \mc{J}^{( 2 )} \dots \mc{J}^{( d )} \subseteq \mc{I}. $ Write
\[
\mc{G}_n = \mc{L}_n \otimes  \mc{J}^{( 1 ) }_n \mc{J}^{( 2 ) }_n
\dots \mc{J}^{( d ) }_n \qquad\mathrm{and}\qquad \mc{H}_n = \mc{L}_n
\otimes \mc{J}^{( 1 ) }_n \otimes \mc{J}^{(  2 ) }_n \otimes \dots
\otimes \mc{J}^{( d ) }_n.
\]

Now let $r = \lfloor n/2d \rfloor$ and $s = n - dr$ and decompose $\mc{H}_n$ as
\[
\mc{H}_n = \mc{L}_s \otimes (\mc{L}_r^{\sigma^s} \otimes \mc{J}^{( 1
) }_n) \otimes (\mc{L}_r^{\sigma^{r+s}} \otimes \mc{J}^{(  2 ) }_n)
\otimes \dots \otimes (\mc{L}_r^{\sigma^{(d-1)r +s}} \otimes
\mc{J}^{( d ) }_n).
\]
By Corollary~\ref{ampleness-sublemma},
$\reg (\mc{L}_r^{\sigma^{(ir + s)}} \otimes \mc{J}^{( i+1 ) }_n ) \leq \dim X + 1$
 for each $i \geq 0$ and all $n \gg 0$. By
Lemma~\ref{regularity and products} there therefore exists a
constant $C$ depending only on $X$ and $\mc{N}$ such that $\reg
\mc{H}_n \leq \reg \LL_s + (\dim X + 1 +C)d$.  As $\{ \mc{L}_n \}$ is an ample sequence,
Corollary~\ref{reg and ample seq} implies that $\lim_{n \to \infty}
\reg \mc{L}_n = - \infty$. Since $s \to \infty$ as $n \to \infty$
this implies that $\lim_{n \to \infty} \reg \mc{H}_n = - \infty$.

The multiplication map $\mu:\mc{H}_n\to \mc{G}_n$ yields   a short
exact sequence  $0\to \mc{K}_n\to \mc{H}_n\to \mc{G}_n\to 0$  such
that each sheaf $\mc{K}_n $ is supported on a finite set.
     For  $i>0$ and $m\in \mb{Z}$,
   the long exact sequence in cohomology therefore  shows that
   $\HB^i(X,\mc{H}_n\otimes \mc{N}^{\otimes m})=0$ $\iff$
    $\HB^i(X,\mc{G}_n\otimes \mc{N}^{\otimes m})=0$. Hence
$\lim_{ n \to \infty}  \reg \mc{G}_n = - \infty$.  Finally, the
inclusion $\mc{J} \subseteq \mc{I}$ induces    an exact sequence $ 0
\to \mc{G}_n \to \mc{R}_n \to \mc{C}_n \to 0 $ where $\mc{C}_n$ is
again supported on a finite set.  A similar long exact sequence
argument implies that $\lim_{n \to \infty} \reg \mc{R}_n = - \infty$
and so, by Corollary~\ref{reg and ample seq}, $\{ \mc{R}_n \}$ is an
ample sequence.

(2) By Proposition~\ref{bimod alg noeth}(1), $\mc{R}$ is right
noetherian and by part (1), the sequence $\{ \mc{R}_n \}$ is ample.
 Thus all of the hypotheses of
Theorem~\ref{VdB main theorem} are satisfied and so the \nsr\  $R =
\mathrm{H}^0(X,\, \mc{R})$ is right noetherian with $\rqgr \mc{R}
\simeq \rqgr R$.

As was noted in \eqref{rees opposite}, $\mc{R}^\op\cong
\mc{R}(X,\sigma(c),\mc{L}^{\sigma^{-1}},\sigma^{-1})$, by
\cite[Corollary~5.1]{Ke1}, $\mc{L}^{\sigma^{-1}}$ is
$\sigma^{-1}$-ample, and obviously $Z$ is saturating with respect to
$\sigma^{-1}$ as well.  Thus the claims on the left follow from
those on the right. This completes the proof.  \qed

One can modify Theorem~\ref{general-ampleness} so that it also works
for curves, but this case is rather uninteresting (see the
discussion in \cite[pp.~511-2]{KRS}),   so it will be ignored.

The results we have proved thus far also give a partial converse to
Theorem~\ref{general-ampleness}.

\begin{proposition}\label{not-noeth} Keep the assumptions of
\eqref{global-convention}. Suppose that $\dim X \geq 2$, $\mc{L}$ is
$\sigma$-ample, and every point $s \in S = \supp Z$ lies on a dense
$\sigma$-orbit, but $Z$ is not right saturating.
  Then neither
$\mc{R}=\mc{R}(X,Z,\mc{L},\sigma)$ nor $R = \mathrm{H}^0(X,\, \mc{R})$ is right
noetherian.
\end{proposition}

\begin{proof} The proof of  \cite[Proposition~4.8]{KRS} works without change,
except that references to \cite[Theorem~3.10 and
Proposition~4.6]{KRS} should be replaced by references to
Proposition~\ref{bimod alg noeth}, respectively
Theorem~\ref{general-ampleness}(1).
\end{proof}

  The aim of the   paper is to obtain a deeper understanding
of the algebras $\mc{R}(X,Z,\mc{L},\sigma)$ and $R=\mathrm{H}^0(X,\,
\mc{R})$ under the assumptions of Theorem~\ref{general-ampleness}
and so we make the following assumptions from now on:

\begin{assumptions}\label{global-convention2}
Let $X$ be a integral projective scheme of dimension $d\geq 2$. Fix
$\sigma\in \Aut(X)$ and a $\sigma$-ample invertible sheaf $\mc{L}$.
Finally assume that $Z\subseteq X$ is a saturating zero-dimensional
subscheme of $X$, or $Z=\emptyset$, with ideal sheaf $\mc{I} =
\mc{I}_Z$.  We will always write $\mc{R} =
\mc{R}(X,Z,\mc{L},\sigma)$ and
 $R=\mathrm{H}^0(X,\, \mc{R}) = R(X,Z,\mc{L},\sigma)$.
By Theorem~\ref{general-ampleness}, $R$ is noetherian with
 $\rqgr R\simeq \rqgr \mc{R}$.
\end{assumptions}

It is often useful to work with connected graded rings that are
generated in degree one and we give two
ways in which this may be achieved; either by replacing $R$ by a
large Veronese ring or by assuming that the invertible sheaf
$\mc{L}$ is ``sufficiently ample.''

\begin{proposition}
\label{veronese} Keep the hypotheses of \eqref{global-convention2}.
Then  the Veronese ring $R^{(q)} = \bigoplus_{n \geq 0} R_{qn}$ is
generated in degree $1$ for all $q\gg 0$.
\end{proposition}

\begin{proof}  The argument is quite similar to
that of \cite[Proposition 4.10]{KRS}, but  some technical adjustments
are needed and so we  will give a full proof.

Note that the Veronese ring $R^{(p)}$ is itself a na{\"\i}ve blowup
algebra, namely $R^{(p)}\cong R(X, Z_p, \mc{L}_p, \sigma^p)$, where $Z_p$ is the
$0$-dimensional subscheme defined by $\mc{I}_p$. Set
$\mc{J}_n=\mc{I}_n\otimes \mc{L}_n$ for $n\geq 1$. By
Theorem~\ref{general-ampleness} we may chose $r\geq 1$ such that
$\mc{J}_r$ is generated by its global sections. Thus, there exists a
short exact sequence
\[
0 \to \mc{V} \to \mathrm{H}^0(X,\mc{J}_r)\otimes \OO_X \to \mc{J}_r
\to 0,
\]
of $\OO_X$-modules, for some sheaf $\mc{V}$. Since  $\{\mc{J}_n\}$
and hence $\{\mc{J}^{\sigma^r}_n \}$ are ample sequences,  there
exists $n_0$ such that $\mathrm{H}^1(X, \mc{V}\otimes
\mc{J}^{\sigma^r}_{nr}) = 0$ for $n\geq n_0$.  Tensoring  the
displayed exact sequence on the right with $\mc{J}_{nr}^{\sigma^r}$
for $nr\geq n_0$ gives an exact sequence
\[
0 \to \tor_1(\mc{J}_r, \mc{J}_{nr}^{\sigma^r}) \to \mc{V} \otimes
\mc{J}_{nr}^{\sigma^r} \overset{\theta}{\lra}
\mathrm{H}^0(X,\mc{J}_r)\otimes \mc{J}_{nr}^{\sigma^r} \to \mc{J}_r
\otimes \mc{J}_{nr}^{\sigma^r} \to 0.
\]
 By \cite[Lemma 3.3]{KRS},  the sheaf $\tor_1(\mc{J}_r,
\mc{J}_{nr}^{\sigma^r})$ is supported on a finite set of points.
Therefore,   if $\mc{K}=\mathrm{Ker}(\theta)$, then  $\HB^1(X,\mc{K}) = \HB^1(X,\mc{V}
\otimes \mc{J}_{nr}^{\sigma^r}) = 0$ for $nr \geq n_0$.  Taking
global sections of the exact sequence
\[
0 \to \mc{K} \to \mathrm{H}^0(X,\mc{J}_r)\otimes
\mc{J}_{nr}^{\sigma^r} \to \mc{J}_r \otimes \mc{J}_{nr}^{\sigma^r}
\to 0,
\]
therefore gives the exact sequence
\begin{equation}
\label{surjection} \mathrm{H}^0(X,\mc{J}_r)\otimes
\mathrm{H}^0(X,\mc{J}_{nr}^{\sigma^r}) \to \mathrm{H}^0(X,\mc{J}_r
\otimes \mc{J}_{nr}^{\sigma^r}) \to 0.
\end{equation}

Consider the  sequence
$
0 \to \tor_1(\mc{L}_r/\mc{J}_r,\mc{J}_{nr}^{\sigma^r}) \to \mc{J}_r
\otimes \mc{J}_{nr}^{\sigma^r} \to \mc{J}_{(n+1)r} \to 0.
$
Since  $\tor_1(\mc{L}_r/\mc{J}_r,\mc{J}_{nr}^{\sigma^r})$ is
supported on a finite set, $\HB^1(X,
\tor_1(\mc{L}_r/\mc{J}_r,\mc{J}_{nr}^{\sigma^r})) = 0$ and so  the map
$\HB^0(X,\mc{J}_r \otimes \mc{J}_{nr}^{\sigma^r}) \to
\HB^0(X,\mc{J}_{(n+1)r})$ is a surjection.  Together with
\eqref{surjection}, this shows that the multiplication map
$R_r\otimes R_{nr} \to R_{(n+1)r}$ is surjective.
By induction,  we obtain
$R_{jr}R_{nr}=R_{(n+j)r}$ for all $j\geq 1$ and   $nr \geq  n_0$.
This implies that for such $nr$ the ring
$R^{(nr)}$ is generated in degree one; that is, by $R_{nr}$.
\end{proof}

\begin{proposition}\label{degree 1}
Keep the hypotheses of \eqref{global-convention2}, and assume that
$\mc{L}$ is also ample and generated by its global sections.  Then
there exists $M \in \NN$ such that, for $m \geq M$:
\begin{enumerate}
\item\label{deg11} $\mc{I}_n \otimes \mc{L}_n^{\otimes m}$ is
generated by its global sections for all $n \geq 1$.
\item\label{deg12} $R(X,Z,\mc{L}^{\otimes m},\sigma)$ is generated in
degree $1$.
\end{enumerate}
\end{proposition}

\begin{proof}   This is
similar to the proof of \cite[Proposition~4.12]{KRS} although, as
happened in the proof of Proposition~\ref{veronese}, one has to
contend  with some nonzero Tor groups. The details are left to the
reader and a full proof can be found in \cite{RS2}.
\end{proof}

The   following technical variant of
Theorem~\ref{general-ampleness}(1) is  needed in \cite{RS}.

\begin{corollary}\label{degree 2}  Keep the hypotheses of
\eqref{global-convention2} and assume that $\mc{L}$ is also very
ample.   Let $\mc{J}$ be an
ideal sheaf such that $\mc{O}_X/\mc{J}$ has finite support.  Then
there exists $M\in \mathbb{N}$ such that $\mathrm{H}^j(X,
\mc{J}^{\sigma^s}\otimes \mc{I}_n\otimes \mc{L}_n^{\otimes m})=0$
 for all $m\geq M$, $n\geq 1$, $j\geq 1$ and  $s\in \mb{Z}$.
 \end{corollary}

\begin{proof}
By Theorem~\ref{general-ampleness}(1), the corollary is true for any fixed value of $s$ and $m$,
so the point is to obtain a uniform bound $M$, which we do by appealing to
Castelnuovo-Mumford regularity.

It is easy to see that, for  $r\geq 0$, there is a constant $E(r)$ with the following property:
If $\mc{M}$ is a sheaf of ideals such that $\mc{O}_X/\mc{M}$ has  length $r$ then
$\reg_{\mc{L}} \mc{M}\leq E(r)$ (see \cite[Examples~1.8.29 and 1.8.30]{Laz1} for the case of $X=\mb{P}^n$).
In particular, $\reg_{\mc{L}}  \mc{J}^{\sigma^s}\leq E$, where $E$ is independent of $s$.

Pick any constant $C$ and follow
 the argument used to prove the equation at
the top of \cite[p.~514]{KRS} for  $\mc{K}_m = \mc{O}_X$. This
shows that there exists a constant $M \geq 1$ such that
 $\reg_{\mc{L}}  \mc{I}_n \otimes
\mc{L}_n^{\otimes m} \leq -E-C+1$ for $n \geq 1$ and  $m \geq M$. Together with
an  application of Lemma~\ref{regularity and products}, this
implies that $\reg_{\mc{L}}  \mc{J}^{\sigma^s}\otimes \mc{I}_n\otimes \mc{L}_n^{\otimes m}
\leq 1$ for all $n \geq 1$, $m \geq M$ and $s \in \mb{Z}$. The
result follows.
\end{proof}

%%%%%%%%%%%%%%%%%%%%%%%%%%%%%%%%

\section{$\mc{R}$-modules and equivalences of categories}
\label{R-modules}

The hypotheses from Assumptions~\ref{global-convention2} will remain
in force throughout this section. One nice consequence of critical
density is that it forces modules over
$\mc{R}=\mc{R}(X,Z,\mc{L},\sigma)$ and $R=\mathrm{H}^0(X,\, \mc{R})$ to have a
very pleasant structure; indeed in many cases they are just induced
from $\OO_X$-modules. This will be used in this section to give
various equivalences of categories, notably that the category of
coherent torsion $\OO_X$-modules is equivalent to the subcategory of
Goldie torsion modules in $\rqgr \mc{R}$, as defined below.
We also
give a natural analogue of the standard fact that, for a blowup
$\rho:\widetilde{X}\to X$ at a smooth point $x$, the schemes
$X\smallsetminus \{x\}$ and $\widetilde{X}\smallsetminus
\rho^{-1}(x)$ are isomorphic.  These results
are largely the same as those in \cite{KRS}, but some of the proofs need more care
since we do not have the identity  $\mc{I}\otimes\mc{I}^\sigma = \mc{I}\mc{I}^\sigma$
that was so useful in  \cite{KRS}.

If $A$ is a noetherian graded domain, a graded $A$-module $M$ is
called  \emph{Goldie torsion} (to distinguish this from the notion
of torsion already defined)  if every homogeneous element of $M$ is
killed by some nonzero homogeneous  element of $A$. Equivalently,
$M$ is a sum of modules of the form $(A/I)[n]$, for nonzero graded
right ideals $I$. The latter notion passes to all the categories $Q$
we consider; for example a right $\mc{R}$-module $\mc{M}$ is Goldie
torsion if it is a sum  of submodules of the form
$(\mc{R}/\mc{K})[n]$ for nonzero right ideals $\mc{K}$ of $\mc{R}$.
Of course, Goldie torsion $\OO_X$-modules are just the torsion
$\OO_X$-modules, as in \cite[Exercise~II.6.12]{Ha}. We write $\GT Q$
for the full subcategory of Goldie torsion modules in~$Q$.

We start by giving some technical results on the structure of Goldie
torsion modules.
If $\mc{N}\in \rGr \mc{R}$, recall from
  \eqref{nice-modules}  that  we may write
  $\mc{N}=\bigoplus {}_1(\mc{G}_n)_{\sigma^n}$ for some
sheaves  $\mc{G}_n$. It is
often convenient to write  ${}_1(\mc{G}_n)_{\sigma^n} = \mc{F}_n \otimes
\mc{L}_{\sigma}^{\otimes n}$, where $\mc{F}_n = {}_1(\mc{F}_n)_1$
has trivial bimodule structure and
$\mc{L}_\sigma^{\otimes n}=(_1\mc{L}_\sigma)^{\otimes n}$.

\begin{lemma}
\label{Goldie tors facts} {\rm{(1)}} Let  $\mc{N}\in \GT\rgr \mc{R}$
and write $\mc{N} = \bigoplus \mc{F}_n
\otimes \mc{L}_{\sigma}^{\otimes n}\in \GT\rgr \mc{R}$ for some
quasicoherent sheaves  $\mc{F}_n$.
Then there
exists a single module $\mc{F}\in \GT \OO_X\catmod$ such that
$\mc{F}_n = \mc{F}$ for all $n \gg 0$.

{\rm{(2)}} Conversely, if $\mc{F}\in \GT \OO_X\catmod$, then
$\bigoplus_{n=0}^{\infty} \mc{F}\otimes
 \mc{L}_{\sigma}^{\otimes n}\in \GT\rgr \mc{R}$.

\end{lemma}
\begin{proof}
The proof of \cite[Lemma~6.1]{KRS} goes through with only minor changes
(replace \cite[Lemma~3.9]{KRS} by Lemma~\ref{coherent}).
\end{proof}

\begin{lemma}\label{goldie-subfactors}
If $\mc{N}\in\rgr \mc{R}$, then there is an exact sequence  $0 \to
\mc{K} \to \mc{N} \to \mc{T} \to 0$, where $\mc{T} \in \GT \rgr
\mc{R}$ is Goldie torsion, and $\mc{K} \in \rgr \mc{R}$ is a direct
sum of shifts of $\mc{R}$.  In fact we can find a $\mc{K}$ such that
$\mc{K} \cong \bigoplus_{i = 0}^d \mc{R}[-n]$ for some $n \geq 0$.
\end{lemma}
\begin{proof}  This  is similar to the proof of the analogous result
for finitely generated modules over domains, and so the proof is left to the reader.
\end{proof}

\begin{lemma}
\label{shifts}  Suppose that $\mc{N}\in \rGr\mc{R}$ and write
 $\mc{N} = \bigoplus \mc{F}_n \otimes  \mc{L}_\sigma^{\otimes n}$ for some
 $\OO_X$-modules $\mc{F}_n$.
  If   $m>0$, then $\mc{N}[m] \cong \bigoplus \mc{G}_n \otimes
\mc{L}_\sigma^{\otimes n}$ where $\mc{G}_n = (\mc{F}_{n+m}\otimes
\mc{L}_m)^{\sigma^{-m}}$.
If $m<0$ then $\mc{N}[m] \cong \bigoplus \mc{G}_n \otimes  \mc{L}_\sigma^{\otimes n}$ where
$\mc{G}_n = (\mc{F}_{n+m})^{\sigma^{-m}} \otimes\mc{L}_{-m}^{-1}$.
\end{lemma}

\begin{proof} If $m>0$ the result   for the bimodule algebra
$\mc{B} = \mc{B}(X,\mc{L},\sigma)$   follows   from \cite[(3.1)]{SV}
and the same argument works for $\mc{R}$. A similar computation then
gives the required  formula for $m\leq 0$.
\end{proof}

In nice cases, the structure maps of an $\mc{R}$-module are
determined by  products instead of tensor products and   in order to
focus in on this property we make the following definition.  Let
   $\mc{N} =
\bigoplus \mc{F}_n \otimes \mc{L}_{\sigma}^{\otimes n}\in
\rGr\mc{R}$; thus the module structure of $\mc{N}$ is determined by
the structure morphisms $\theta_n: \mc{F}_n \otimes
\mc{I}^{\sigma^n} \to \mc{F}_{n+1}$ for   $n \in \mb{Z}$. Then
$\mc{N}$ is \emph{definable by products in degrees $\geq n_0$} if
for all $n \geq n_0$, the maps $\theta_n$ factor through the
multiplication map
  $\mu_n: \mc{F}_n \otimes
\mc{I}^{\sigma^n} \twoheadrightarrow  \mc{F}_n \mc{I}^{\sigma^n}$ in
the sense that there exist
 maps $\rho_n: \mc{F}_n
  \mc{I}^{\sigma^n} \to \mc{F}_{n+1}$ of $\mc{O}_X$-modules such that $\theta_n = \rho_n  \mu_n$.
We say $\mc{N}$ is \emph{definable by products} if some such $n_0$
exists.

As the next result shows,   any  coherent module $\mc{N}$ is
  definable by products.  For the \naive\ blowup at a single point, one has
  the stronger statement  that if $\mc{N}$
is coherent then $ \mu_n$ is itself an isomorphism for $n \gg 0$
\cite[Lemma~6.3]{KRS}. This stronger  result will not always be true
in our  setting; for instance it fails  for $\mc{N} = \mc{R}$ in the
case where $\supp Z$ contains two distinct points from the same
$\sigma$-orbit.

\begin{lemma}
\label{product versus tensor} Let $\mc{N}  = \bigoplus
\mc{F}_n\otimes \mc{L}_{\sigma}^{\otimes n} \in \rgr \mc{R}$.  Keep
the notation from the discussion above.  Then
\begin{enumerate}
\item There exists $n_0\geq 0$ such that
 $\mc{N}$ is definable by products in degrees $\geq n_0$.
\item  There exists $n_1\geq n_0$ such that  $\rho_n$ is an isomorphism
 for $n \geq n_1$.
\end{enumerate}
\end{lemma}

\begin{proof}   To prove that $\theta_n$ factors through the map  $\mu_n$
  is equivalent to proving  that $\ker  \mu_n \subseteq \ker
\theta_n$. Tensoring the   sequence $0 \to
\mc{I}^{\sigma^n} \to \mc{O}_X \to \mc{O}_X/\mc{I}^{\sigma^n} \to 0$
with $\mc{F}_n$ gives the exact sequence
\[
0 \to  \calTor_1^{\mc{O}_X}(\mc{F}_n, \mc{I}^{\sigma^n})
\overset{\alpha_n}{\lra} \mc{F}_n \otimes \mc{I}^{\sigma^n} \to
\mc{F}_n \mc{I}^{\sigma^n} \to 0,
\]
so we can and will  identify $\ker  \mu_n$ with $\calTor_1^{\mc{O}_X}(\mc{F}_n,
\mc{I}^{\sigma^n})$.

Suppose first  that $\mc{N} = \mc{R}[m]$
 is a shift of $\mc{R}$.  Lemma~\ref{shifts} implies that
 $\mc{F}_n = \mc{I}^{\sigma^{-m}}\cdots\mc{I}^{\sigma^{n-1}}
 \otimes \mc{L}_{|m|}^\alpha$
 for  $n>|m|$ and the appropriate  $\alpha$.  Then for large $n$, the map
$\theta_n: \mc{F}_n \otimes \mc{I}^{\sigma^n} \to \mc{F}_{n+1}$ is
then none other than the natural map from the tensor product
$\mc{I}^{\sigma^{-m}}\cdots\mc{I}^{\sigma^{n-1}} \otimes
\mc{I}^{\sigma^n}$ to the product
$\mc{I}^{\sigma^{-m}}\cdots\mc{I}^{\sigma^{n-1}} \mc{I}^{\sigma^n}$,
tensored by $\mc{L}_{|m|}^\alpha$.  So (1) holds   when
$\mc{N}= \mc{R}[m]$.

Now let $0 \to \mc{N}' \to \mc{N} \to \mc{N}''\to 0$ be an exact
sequence in $\rgr \mc{R}$, where $\mc{N}' = \bigoplus \mc{F}'_n
\otimes \mc{L}_{\sigma}^{\otimes n}$ and $\mc{N}'' = \bigoplus
\mc{F}''_n \otimes \mc{L}_{\sigma}^{\otimes n}$, with structure maps
$\theta_n'$ and $\theta''_n$, respectively.
 Consider the  commutative diagram:
 \begin{equation}\label{tensor1}
 \begin{CD}
 \calTor_1^{\mc{O}_X}(\mc{F}'_n, \mc{I}^{\sigma^n})   @>>>
  \calTor_1^{\mc{O}_X}(\mc{F}_n, \mc{I}^{\sigma^n}) @>>>
 \calTor_1^{\mc{O}_X}(\mc{F}''_n, \mc{I}^{\sigma^n})    \\
@VV \alpha' V  @VV \alpha V  @VV \alpha'' V \\
 \mc{F}'_n \otimes \mc{I}^{\sigma^n}  @>>>
     \mc{F}_n \otimes \mc{I}^{\sigma^n} @>>>
     \mc{F}''_n \otimes \mc{I}^{\sigma^n} \\
    @VV \theta' V @ VV \theta V @VV \theta'' V \\
        \mc{F}'_{n+1}  @>>>
     \mc{F}_{n+1} @>>>
     \mc{F}''_{n+1}
 \end{CD}
\end{equation}
We first suppose that $\mc{N}=\mc{N}'\oplus \mc{N}''$ and that (1) holds for
$\mc{N}'$ and $\mc{N}''$. In this case, the outside columns of \eqref{tensor1} are
complexes and the rows split. Thus the middle column is also a complex and so
(1) holds for $\mc{N}$. By induction, we have therefore proved that (1) holds when
 $\mc{N}$ is a direct sum of shifts of $\mc{R}$.

By
Lemma~\ref{goldie-subfactors},
a general coherent module  $\mc{N}$ fits into an exact sequence $0 \to
\mc{N}' \to \mc{N} \to \mc{N}'' \to 0$ where $\mc{N}'$ is a sum of
shifts of $\mc{R}$ and $\mc{N}''$ is Goldie torsion.
  By the last paragraph, the first column of
     \eqref{tensor1} is now  a complex.
By Lemma~\ref{Goldie tors facts}  there exists a coherent   torsion
 sheaf $\mc{F}''$ such that $\mc{F}''_n=\mc{F}''$, say
for all $n\gg n_0$. Since $Z$ is saturating,
 $\supp \mc{F}'' \cap \supp \mc{O}_X/\mc{I}^{\sigma^n}=\emptyset$ for $n\gg n_0$
and thus $\calTor_1^{\mc{O}_X}(\mc{F}''_n, \mc{I}^{\sigma^n}) = 0$,
for such $n$. A simple diagram chase then shows that the middle
column of \eqref{tensor1} is a complex, proving  part~(1).

(2) Lemma~\ref{coherent} implies that
  $\theta_n$ is surjective for $n \gg 0$, while  $\rho_n: \mc{F}_n \mc{I}^{\sigma^n} \to
\mc{F}_{n+1}$ is defined for all $n\gg 0$ by part~(1). Thus there exists $n_0$,
 such that $\rho_n$ is defined and
surjective for all $n \geq n_0$ and it remains to prove that $\rho_n$
is an isomorphism for $n \gg n_0$.

Pulling back to $\mc{F}_{n_0} $, for $n \geq n_0$ we may write
$\mc{F}_{n} = \mc{A}_n/\mc{B}_n$, for subsheaves $\mc{B}_n\subseteq
\mc{A}_n\subseteq \mc{F}_{n_0}$. Since $\mc{F}_{n+1}$ is a
homomorphic image of $\mc{F}_n\mc{I}^{\sigma^n} = (\mc{A}_n
\mc{I}^{\sigma^n} + \mc{B}_n)/\mc{B}_n$, we find that $\mc{B}_{n+1}
\supseteq \mc{B}_n$ for each $n\geq n_0$.
 Since $\mc{F}_{n_0}$ is noetherian, $\mc{B}_n=\mc{B}_{n+1} $ for all
 $n\gg n_0$, and hence
$\mc{F}_n \mc{I}^{\sigma^n} \cong \mc{F}_{n+1}$ for such $n$.
\end{proof}

\begin{examples}
\label{definable} (1)
Let $\mc{F} \in \mc{O}_X\Mod$ be  any
quasi-coherent sheaf on $X$ and  take $\mc{M} =
\bigoplus_{n \geq 0} \mc{F} \otimes \mc{L}_{\sigma}^{\otimes n} $ in $
\rGr \mc{R}$. Then $\mc{M}$ is easily seen to be definable by
products in degrees $\geq 0$, even though  $\mc{M}$ need  not be
coherent as an $\mc{R}$-module.

(2) In general, however,
a non-coherent $\mc{R}$-module $\mc{M}\in \rGr \mc{R}$
 need not be   definable by products. For example,
take $\mc{R} =\mc{R}(\mathbb{P}^2, Z, \mc{L}, \sigma)$ where $Z= p$
is a single reduced  saturating point with ideal sheaf $\mc{I}$ and
set $\mc{N}=\mc{O}_X/\mc{I} \oplus (\mc{I}/\mc{I}^2)\otimes
\mc{L}_\sigma\oplus 0\oplus 0\cdots$, where the structure map
$\theta_0: (\mc{O}_X/\mc{I}) \otimes \mc{I} \to \mc{I}/\mc{I}^2$ is
the natural isomorphism. Then $\mc{N}$  is clearly not
 definable by products in degrees $\geq 0$ and so
  $\mc{N}[-n]$   is not
 definable by products in degrees $\geq n$. Thus  $\bigoplus_{n\geq 0} \mc{N}[-n]$
   is not definable by products in any degree.
  \end{examples}

For   $\mc{R}$-modules  that are  definable by products, the
homomorphism groups also have a nice form and the next lemma
collects the relevant facts.  Recall from  Section~\ref{definitions}
 that  the  map from
 $\rGr\mc{R}$ to $\rQgr\mc{R}$ is denoted~$\pi$.

\begin{lemma}
\label{Hom in Qgr} Suppose that $\mc{M} = \bigoplus \mc{F}_n
\otimes\mc{L}_{\sigma}^{\otimes n} \in \rgr \mc{R}$ is coherent, and
that $\mc{N} = \bigoplus \mc{G}_n \otimes\mc{L}_{\sigma}^{\otimes
n}\in \rGr \calR$ is   definable by products.
\begin{enumerate}
\item  For some $n_1$, there is a natural isomorphism
$
\Hom_{\rQgr \mc{R}}(\pi(\mc{M}), \pi(\mc{N})) \cong \lim_{n \geq
n_1} \Hom_{\mc{O}_X}(\mc{F}_n, \mc{G}_n).$
\item Suppose that $\mc{M}$ and $\mc{N}$ are coherent
Goldie torsion modules  and, by Lemma~\ref{Goldie tors facts}, write
$\mc{F}_n = \mc{F}$ and $\mc{G}_n = \mc{G}$   for $n \gg  0$. Then
there is a natural isomorphism
$$
\Hom_{\rQgr \mc{R}}(\pi(\mc{M}), \pi(\mc{N})) \cong
\Hom_{\mc{O}_X}(\mc{F}, \mc{G}).
$$
\end{enumerate}
\end{lemma}

Before beginning the proof we need to explain the direct limit
appearing in part~(1). By Lemma~\ref{product versus tensor}, pick
$n_1$ so that $\rho_n: \mc{F}_n \mc{I}^{\sigma^n} \to \mc{F}_{n+1}$
is an isomorphism for all $n \geq n_1$.  We can also assume
that $\mc{N}$ is definable by products in degrees $\geq n_1$, so the
 structure maps $\theta'_n: \mc{G}_n \otimes \mc{I}^{\sigma^n} \to \mc{G}_{n+1}$ factor through
 maps $\rho'_n: \mc{G}_n \mc{I}^{\sigma^n} \to \mc{G}_{n+1}$ for $n \geq n_1$.
 Given $\alpha\in \Hom_{\mc{O}_X}(\mc{F}_n, \mc{G}_n)$ where $n \geq n_1$, then
$\alpha$ restricts to  a morphism $\alpha'\in
\Hom_{\mc{O}_X}(\mc{F}_n \mc{I}^{\sigma^n}, \mc{G}_n
\mc{I}^{\sigma^n})$, and so  $\alpha'$  induces a map
$\alpha''=\rho'_n\circ \alpha'\circ  \rho_n^{-1} \in
\Hom_{\mc{O}_X}(\mc{F}_{n+1}, \mc{G}_{n+1})$.

\begin{proof}[Proof of Lemma~\ref{Hom in Qgr}]
(1)  The definition of homomorphisms in quotient categories implies
that
\begin{equation}\label{hom-in-quotient}
\Hom_{\rQgr \mc{R}}(\pi(\mc{M}),\pi(\mc{N})) = \lim_{n\to \infty}
\Hom_{\rGr \mc{R}}(\mc{M}_{\geq n},\mc{N}),
\end{equation}
whenever $\mc{M} = \bigoplus \mc{F}_n
\otimes\mc{L}_{\sigma}^{\otimes n}$ is coherent (see, for example,
\cite[p.~31]{VB2}).
On the other hand, we claim that there are natural vector space maps
\begin{equation}\label{hom-in-quotient2}
\hom_{\rGr \calR}(\mc{M}_{\geq n}, \mc{N}) \overset{\phi_n}{\lra}
\hom_{\OO_X} (\mc{M}_n, \mc{N}_n) \overset{\psi_n}{\lra}
\hom_{\mc{O}_X}(\mc{F}_n, \mc{G}_n).
\end{equation}
Indeed, if $f \in \hom_{\rGr \calR}(\mc{M}_{\geq n}, \mc{N})$, then
$f$ is a morphism of right $\mc{O}_X$-modules, so define
$\phi_n(f)$ to be the restriction of $f$ to $\mc{M}_n$. The map
$\psi_n$ is the natural isomorphism obtained by tensoring with  $(\mc{L}_{\sigma}^{\otimes n})^{-1}.$

 We will prove that  $\psi_n\circ\phi_n$ is an isomorphism
for $n\geq n_1$.  Given $g\in \Hom_{\mc{O}_X}(\mc{F}_n, \mc{G}_n)$
for $n \geq n_1$, then $g$ induces a unique map
$$\mc{M}_{n+r} \cong
  \mc{F}_n \mc{I}_r^{\sigma^n} \otimes\mc{L}_\sigma^{\otimes (n+r)}
 \to \mc{G}_n \mc{I}_r^{\sigma^{n}}\otimes\mc{L}_\sigma^{\otimes (n+r)}
 \to \mc{G}_{n+r}  \otimes\mc{L}_\sigma^{\otimes (n+r)}=\mc{N}_{n+r},$$
for any $r\geq 0$. These piece together to give an
 $\mc{R}$-module map $f \in \hom_{\OO_X}
(\mc{M}_{\geq n}, \mc{N})$ and in this way we define a morphism
$\tau_n: \Hom_{\mc{O}_X}(\mc{F}_n, \mc{G}_n) \to \hom_{\OO_X}
(\mc{M}_{\geq n}, \mc{N})$ by setting $\tau_n(g) = f$. Obviously
$\psi_n \phi_n \tau_n(g) = g$ and so $\psi_n \phi_n$ is surjective.
Moreover, for $n \geq n_1$ any element in
 $\hom_{\rGr \calR}(\mc{M}_{\geq n}, \mc{N})$ is easily seen to be determined by its
 restriction to degree $n$, because $\mc{M}_{n+r} = \mc{M}_n
 \mc{R}_r$ for all $r \geq 0$.  It follows that $\psi_n
 \phi_n$ is injective and so $\psi_n \phi_n$ is an isomorphism
 for $n \geq n_1$
 as claimed.

The isomorphisms $\psi_n \phi_n$ are
compatible with the maps in the direct limits, and so they induce an
isomorphism $\lim_{n\to \infty} \Hom_{\rGr \mc{R}}(\mc{M}_{\geq
n},\mc{N}) \to \lim_{n \geq n_1} \Hom_{\mc{O}_X}(\mc{F}_n,
\mc{G}_n)$.  So we are done by \eqref{hom-in-quotient}.

(2) The $n^{\mathrm{th}}$ map in the
direct limit $\lim_{n \geq n_1} \Hom_{\mc{O}_X}(\mc{F}_n, \mc{G}_n)$
is an isomorphism  for all large $n$  and so the direct limit stabilizes
at $\Hom_{\mc{O}_X}(\mc{F}, \mc{G})$. Thus
part~(2)  is a special case of part~(1).
\end{proof}

One of the most  important cases  of  Lemma~\ref{Hom in Qgr}(1) occurs when
$\mc{M}=\mc{R}$. In this case,  \cite[Lemma~6.4]{KRS} shows that part~(1)
holds for all modules $\mc{N}$
provided  one \naive ly blows up a single point. However, when one blows up more
 than one point at once, then  Lemma~\ref{Hom in Qgr}(1) can fail for a  general
 module $\mc{N}$. Since the example is a little technical we will omit it, although it  can be found in   \cite{RS2}.

It is now easy to define an
 equivalence of categories between $\GT \rqgr R$ and $\GT
 \OO_X\catmod$, thereby  proving Theorem~\ref{mainthm}(4).

\begin{theorem}
\label{GT equiv} Keep the hypotheses from \eqref{global-convention2}.
Then there are equivalences of categories
$$\GT \rQgr R\ \simeq\ \GT \rQgr \mc{R} \ \simeq\ \GT \OO_X\catMod,$$
which restrict to equivalences
$\GT \rqgr R\simeq \GT \rqgr \mc{R} \simeq \GT \OO_X\catmod.$
This equivalence is given by mapping $\mc{F}\in \GT \OO_X\catMod$
to $\pi\left(\bigoplus \mc{F}\otimes \mc{L}^{\otimes n}_\sigma\right)
\in \rQgr \mc{R}$.
\end{theorem}

\begin{proof}  The proof of \cite[Theorem~6.7]{KRS} goes through with the following minor changes.
Specifically, the references to Theorem~4.1,  Lemma~6.1 and Lemma~6.4 of  \cite{KRS} given in that
proof  should be replaced by references to
Theorem~\ref{general-ampleness}(2), Lemma~\ref{Goldie tors facts}
and Lemma~\ref{Hom in Qgr}, respectively.  (In the proof, Lemma~\ref{Hom in Qgr}
is only applied to coherent $\mc{R}$-modules and so, by Lemma~\ref{product versus tensor},
the hypotheses of Lemma~\ref{Hom in Qgr} are satisfied.)
\end{proof}

   If $k(x)$ is the
skyscraper sheaf at a closed point   $x\in X$,  set $\overline{x}=
\bigoplus_{n\geq 0} (k(x) \otimes \mc{L}_\sigma^{\otimes n}) \in
\rGr \mc{R}$ and write $\widetilde{x}=\pi(\overline{x}) \in
\rQgr\mc{R}$. By Lemma~\ref{Goldie tors facts}(2) $\overline{x}\in
\GT\rgr\mc{R}$ and so $\widetilde{x} \in \rqgr\mc{R}$. Combined with
Proposition~\ref{non-represent2}, the next result proves
Theorem~\ref{mainthm}(5).

\begin{corollary}\label{GT equiv1}    Keep the hypotheses from
\eqref{global-convention2}. Then:
\begin{enumerate}
\item There is  a (1-1) correspondence between the closed points $x\in X$ and
isomorphism classes of simple
objects in  $ \rqgr \mc{R}$ given by $x\mapsto \widetilde{x}$.
\item The simple objects in $\rqgr R$  are the images of finitely
generated $R$-modules $M\in \rgr R$ with Hilbert series
$(1-t)^{-1}$.
\item  If $R$ is generated in degree one then the simple objects in $\rqgr
R$   are the images  of shifts of point modules.
\end{enumerate}
\end{corollary}

\begin{proof} (1)
 Clearly  $\widetilde{x}$
   is the image of   $k(x)$  under the equivalence  from Theorem~\ref{GT equiv}
 and so the simple objects
in $\rqgr\mc{R}$ are exactly these   $\widetilde{x}$.

(2) Given $\widetilde{x}$ for a closed point $x\in X$,
Theorem~\ref{VdB main theorem}(2) shows that $M(x) =
\mathrm{H}^0(X,\widetilde{x})\in \rgr R$ and so $M(x)_n=
\mathrm{H}^0(X, k(x)\otimes \mc{L}^{\otimes n}_\sigma)\cong
\mathrm{H}^0(X, k(x))$ is one-dimensional for all $n\geq 0$. By
Theorem~\ref{VdB main theorem}(2), again, the image of $M(x)$ in
$\rqgr R$ corresponds to $\widetilde{x}$ under the equivalence
$\rqgr R \simeq \rqgr \mc{R}$.

(3) Take $M(x)$ as in part~(2) and suppose that it is generated in degrees
$\leq r_0$. Since $R$ is generated in degree one it follows that
 $M(x)_n = M(x)_{r_0}R_{n-r_0}$, for all $n\geq r_0$. Thus
$M(x)_{\geq r_0}$ is a shifted point module.
\end{proof}

Unfortunately the module $M(x)$ constructed in the  proof of
Corollary~\ref{GT equiv1}(2)  will not be cyclic  when $x\in
\bigcup_{n\geq 0} \supp\sigma^{-n}( Z)$ and so one needs a more
subtle argument to define the simple objects in $\rqgr R$ in terms
of modules generated in degree zero. The details are given in
Section~\ref{further-section}.

If $\rho: \widetilde{X}\to X$ is the (classical) blowup of $X$ at a
smooth point $x$, then it is standard that  $X\smallsetminus\{x\}
\cong \widetilde{X}\smallsetminus \rho^{-1}(x)$.
The final result of this section gives the
analogous result for $\rqgr\mc{R}$, although we have to remove whole
$\sigma$-orbits rather than just isolated points. Thus, define $C_X$
to be the smallest localizing subcategory of $\mc{O}_X\catMod$
containing all
 the modules $\{k(c) | c \in \bigcup_{i \in \mb{Z}} \sigma^i(S)
\}$, where $S = \supp \mc{O}_X/\mc{I}$.
Similarly, write $C_\mc{R}$ for the localizing subcategory of
$\rQgr\mc{R}$ generated by the modules $\widetilde{c}$ for $c \in
\bigcup_{i \in \mb{Z}} \sigma^i(S)$.

\begin{proposition}
\label{factor-equiv} Assuming  \eqref{global-convention2},
there is an equivalence of  categories
$\OO_X \catMod/C_X  \simeq \rQgr \mc{R}/C_{\mc{R}}.$
\end{proposition}

\begin{proof}  The proof of this result is similar to that of
  \cite[Preposition~6.9]{KRS}, but since the result is peripheral, we leave the details to the
   interested  reader. A full proof can be found in \cite{RS2}.
 \end{proof}

%%%%%%%%%%%%%%%%%%%%%%%%%%%%%%%%

\section{Generalized na{\"\i}ve blowups and torsion extensions}
\label{gen naive blowups}

Throughout the section,  the hypotheses from
Assumptions~\ref{global-convention2} will be maintained.
If one  \naive ly blows up a single point then the corresponding  \nsr\  $R$
automatically satisfies $\chi_1$ (see \cite[Theorem~1.1(8)]{KRS}).
The $\chi$  conditions are defined in Section~\ref{chi}, but in this
section  we will just be interested in the following weaker version:
A cg Goldie domain  $R$ satisfies \emph{the weak (right)
$\chi_1$-condition} if,  given any cg algebra $R\subseteq S\subseteq
Q(R)$ such that $S/R$ is (right) torsion then $S/R$ is finite
dimensional. Remarkably, this can fail when one blows up at more
than one point. In order to analyse this situation we need to
understand the maximal torsion extensions of a \nsr\ and this leads
to a variant of na{\"\i}ve blowups, called generalized na{\"\i}ve
blowups. These will be studied in this section and applied to the
study of the chi conditions in    Section~\ref{chi}.

Here is  a simple example of this phenomenon. More examples
 will appear at the end of the section.

 \begin{example}\label{easy-eg1}
Fix $X=\mathbb{P}^2$ and
 $\sigma\in \mathrm{Aut}(X)$ for which there exists a closed point $c=c_0\in X$
 with a critically dense $\sigma$-orbit.
 Write $\lcl{\mf{m}}{0} $ for the sheaf of maximal  ideals corresponding to $c_0$.
If $\lcl{\mf{m}}{0}=(x,y)$ locally at $c_0$, let $\lcl{\mc{M}}{0}$ be the sheaf of
ideals such that $\mc{O}_X/\lcl{\mc{M}}{0}$ is supported at $c_0$ but such
that $\lcl{\mc{M}}{0}=(x^2,y^2)$ locally at $c_0$. For $i\in \mathbb{Z}$
set $c_i=\sigma^{-i}(c_0)$ and write $\lcl{\mathfrak{m}}{i} =
\lcl{\mf{m}}{0}^{\sigma^i}$ and
  $\lcl{\mc{M}}{i} = \lcl{\mc{M}}{0}^{\sigma^i}$.  The key property  of
 the $\lcl{\mc{M}}{i}$ is that $\lcl{\mc{M}}{i} \subsetneq \lcl{\mathfrak{m}}{i}^2$
 but $\lcl{\mc{M}}{i}  \lcl{\mathfrak{m}}{i}= \lcl{\mathfrak{m}}{i}^3$.

Let $\mc{I}=\lcl{\mathfrak{m}}{0} \lcl{\mc{M}}{1}$ and
$\mc{H}=\lcl{\mathfrak{m}}{0} \lcl{\mathfrak{m}}{1}^2$. A routine
computation shows that
$$\mc{I}_n=\mc{I}\mc{I}^{\sigma}\cdots\mc{I}^{\sigma^{n-1}}
= \lcl{\mathfrak{m}}{0} \lcl{\mathfrak{m}}{1}^3 \cdots
\lcl{\mathfrak{m}}{n-1}^3 \lcl{\mc{M}}{n} \qquad \mathrm{but} \qquad
\mc{H}_n =\mc{H}\mc{H}^{\sigma}\cdots\mc{H}^{\sigma^{n-1}} =
\lcl{\mathfrak{m}}{0}\lcl{\mathfrak{m}}{1}^3\cdots
\lcl{\mathfrak{m}}{n-1}^3 \lcl{\mathfrak{m}}{n}^2,$$ and
 so $\mc{H}_n\mc{I}_r^{\sigma^n}=\mc{I}_{n+r}$ for all $n \geq 0$ and $r\geq
 1$.   Thus $\mc{R}=\mc{R}(\mathbb{P}^2,Z_{\mc{I}},\mc{L},\sigma) \subset
\mc{T}=\mc{R}(\mathbb{P}^2,Z_{\mc{H}},\mc{L},\sigma)$  satisfy
$\mc{T}\mc{R}_{\geq 1} \subseteq \mc{R} $, despite the fact that
$\mc{R}_n\not=\mc{T}_n$   for all $n\geq 1$.

 Now take $\mc{L}=\mc{O}_{\mathbb P^2}(m)$ where $m$ is large enough so that
$R=R(\mathbb{P}^2,Z_{\mc{I}},\mc{L},\sigma)$ and
$T=R(\mathbb{P}^2,Z_{\mc{H}},\mc{L},\sigma)$ are both generated in
degree one and each $\mc{I}_n\otimes \mc{L}_n$ and $\mc{H}_n\otimes
\mc{L}_n$ is generated by its sections (this is possible by
Proposition~\ref{degree 1}).
Then the conclusion of the last paragraph translates into the statement that
 $R\subset T$ with $TR_{\geq 1} \subset R$ despite the fact that $T_n\not=R_n$ for each $n\geq 1$.
 Clearly, this shows that $R$ does not satisfy the weak $\chi_1$ condition on the right.  \qed
\end{example}

The rings $R\subset T$  from Example~\ref{easy-eg1} have a number of other interesting
properties that will become more  evident  as we develop the appropriate theory.
For example, $R$ is quite asymmetric and does
satisfy weak $\chi_1$ on the left (see the discussion immediately before Example~\ref{easy-eg2}).
  Examples like this are intimately
connected to the theory of idealizer rings; in fact, the ring $R$ is
the idealizer in $T$ of the left ideal $TR_{\geq 1}$ (see
Lemma~\ref{ideal-lem} and the discussion thereafter).
 Although in this example $T$ is itself a \nsr, this
does not always happen (see Example~\ref{not-naive}), and to cater for examples like that
  we will need to work with the following more general objects.

\begin{definition}
\label{gen naive def}  Keep the hypotheses from Assumptions~\ref{global-convention2}.
A \emph{\gns} is a
sequence $\{  \mcII_n \}_{n \geq 0}$ of ideal sheaves on $X$
satisfying the following properties:
\begin{enumerate}
\item $ \mcII_0 = \mc{O}_X$ and $ \mcII_m  \mcII_n^{\sigma^m} \subseteq  \mcII_{n+m}$ for all
$m, n \geq 0$.
\item There exists a constant $t \geq 1$ such that
$ \mcII_m  \mcII_n^{\sigma^m} =   \mcII_{n+m}$ for all $m, n \geq t$.
\item For $n \geq 0$, the subscheme $\supp \mc{O}_X/ \mcII_n$
 is either zero-dimensional and saturating, or empty.  
\end{enumerate}
 If (2) holds with $ t = 1$, then $\{  \mcII_n \}$ is called a
\emph{na{\"\i}ve sequence}.
\end{definition}

Given this data, we write $\mc{S}= \mc{S}(X, \{ \mcII_n\},\mc{L},
\sigma) = \bigoplus_{n\geq 0}( \mcII_n\otimes \mc{L}_n)_{\sigma^n}$.
This  is easily seen to be a bimodule algebra, which we call a
 \emph{\gnbba}.
 This notation is justified since, if
 $\{ \mcII_n\}$ is  a na{\"\i}ve sequence then
 $\mcII_n=\mcII_1 \cdots\mcII_1^{\sigma^{n-1}}$
 for all $n\geq 0$ and so  $\mc{S}$  is just  the
 bimodule algebra $\mc{R} = \mc{R}(X, Z_{{ \mcII}_1},\mc{L},\sigma)$.
The algebra of sections $\mathrm{H}^0(X,\mc{S}) =S(X, \{
\mcII_n\},\mc{L},\sigma)$ will be called \emph{a \gnba}.  We call
$\mc{S}$ or $S$ \emph{nontrivial} if $\mcII_n\not=\mc{O}_X$ for some
(and hence all) $n\gg0$.

 As we next show, many of the basic properties of $\mc{R}$
generalize easily to $\mc{S}$.

\begin{lemma}\label{gen-naive1}
Let $\{  \mcII_n \}_{n \geq 0}$ be a \gns\ and set $\mc{S}=
\mc{S}(X, \{ \mcII_n\},\mc{L}, \sigma)$.  Take $t$ as in
Definition~\ref{gen naive def}(2) and pick any $p \geq t$.
Then the   Veronese bimodule algebra $\mc{S}^{(p)} =
\bigoplus_{n\geq 0} ( \mcII_{np} \otimes \mc{L}_{np})_{\sigma^{np}}
$ equals  $\mc{R}(X, Z_{ \mcII_p}, \mc{L}_p, \sigma^p)$.
Moreover, the sequence $\{( \mcII_{np} \otimes \mc{L}_{np})_{\sigma^{np}} \}$ is ample and
both  $\mc{S}^{(p)}$ and its section ring
  $S^{(p)} = \bigoplus_{n\geq 0} \mathrm{H}^0(X,\, \mc{S}_{np})$ are noetherian.
 \end{lemma}

\begin{proof} Note that  $ \mcII_{np} =  \mcII_p
 \mcII_p^{\sigma^p} \cdots  \mcII_p^{\sigma^{np-p}}$ for all $n \geq
 1 $ and so
  $\mc{S}^{(p)} = \mc{R}(X, Z_{ \mcII_p}, \mc{L}_p,
\sigma^p)$ holds by definition. The hypotheses from
Assumptions~\ref{global-convention2} (and hence those from
\ref{global-convention}) pass to Veronese subsequences and so the result follows from
Theorem~\ref{general-ampleness}.
\end{proof}

\begin{corollary}\label{gen-naive2} Let $\{  \mcII_n \}_{n \geq 0}$ be a
\gns\  and pick $p$ as in Lemma~\ref{gen-naive1}.
Then:
\begin{enumerate}
\item   $\{  \mcII_n \otimes \mc{L}_n \}_{n \geq 0}$ is an ample sequence.
\item $\mc{S}=
\mc{S}(X, \{ \mcII_n\},\mc{L}, \sigma) $ is coherent as a left and right $\mc{S}^{(p)}$ module
and is a noetherian bimodule algebra.
\item $S = \HB^0(X, \mc{S})$ is a noetherian ring that is finitely generated as a right or left
module over $S^{(p)}$.
\item There exists a  constant  $t'$ such that $S_m S_n = S_{m+n}$ for
all $m, n \geq t'$.
\end{enumerate}
\end{corollary}

\begin{proof}
(1) By  Definition~\ref{gen naive def}(2), we can choose $m_0$ such
that $\mcII_m \mcII_{np}^{\sigma^m} = \mcII_{m+np}$ for all $m_0
\leq m \leq m_0 +p$ and $n \geq 0$.   By Lemma~\ref{gen-naive1}, the
sequence $\{  \mcII_{np} \otimes \mc{L}_{np} \}_{n \geq 0}$ is
ample and hence, for each  $m_0 \leq m \leq m_0 +p$,  so is the
sequence $\{  (\mcII_{np} \otimes \mc{L}_{np})^{\sigma^m} \}_{n \geq
0}$.  For  such $m$ and coherent sheaf $\mc{F}$, the natural
surjection
\[
\mc{F}\otimes (\mcII_m \otimes \mc{L}_m) \otimes (\mcII_{np} \otimes
\mc{L}_{np})^{\sigma^m} \twoheadrightarrow
\mc{F} \otimes \left(\mcII_{n + mp} \otimes \mc{L}_{n+mp}\right)
\]
has a  finitely supported kernel, from which it follows that  $\{ \mcII_{np+m} \otimes
\mc{L}_{np+m} \}_{n \geq 0}$   is also ample.
  Therefore $\{ \mcII_n \otimes \mc{L}_n \}_{n \geq 0}$ is  ample.

(2) We consider $\mc{S}$ as a right  $\mc{S}^{(p)}$-module via the ungraded inclusion of
bimodule algebras $\mc{S}^{(p)} \subseteq \mc{S}$.
By Lemma~\ref{coherent}, in order to show
that $\mc{S}$ is a coherent right $\mc{S}^{(p)}$-module,  it suffices  to show that $ \mcII_m
 \mcII_{np}^{\sigma^m} =  \mcII_{m+ np}$ for all $m \gg 0$ and $n
\geq 1$.  This holds by
Definition~\ref{gen naive def}(2).
It then follows  from Lemma~\ref{gen-naive1} and
 \cite[Proposition~2.10]{KRS}    that $\mc{S}$ is right noetherian.
The same argument works on the left.

(3)   By   Lemma~\ref{gen-naive1}, the hypotheses of Theorem~\ref{VdB main theorem}
are satisfied by $\mc{S}^{(p)}$. Thus, by part~(2) 
and Theorem~\ref{VdB main theorem}(2),
$S=\mathrm{H}^0(X,\, \mc{S})$ is  noetherian as both  a right and a left  $S^{(p)}$-module.

(4)  By Lemma~\ref{gen-naive1} $S^{(p)}$ is a \nsr\ for all $p\geq
t$ and so, by
 Proposition~\ref{veronese}, the Veronese ring
$S^{(q)}$ is generated in degree one for all large multiples $q$ of $p$.
 As $S$ is a finitely generated right $S^{(q)}$-module, this
implies that, for some $m_0$, $S_mS_{qr}=S_{m+qr}$ for all $r\geq 0$
and $m \geq m_0$.  Varying $p$, we can find two relatively prime
integers $q_1, q_2$ and a single $m_0$ such that the conclusion of
the previous sentence holds for both $q = q_1$ and $q = q_2$. Since
any $n \gg 0$ can be written as $n = aq_1 + b q_2$ for some $a, b
\geq 0$,
  it follows that $S_m S_n = S_{m+n}$ for $m \geq m_0$ and all $n
\gg 0$.
\end{proof}

We want to understand the asymptotic behaviour of a
\gns\ $\{\mcII_n\}$,   for which we need the following  notation.

\begin{notation}\label{cosupp-defn}
If $\mc{J}\subseteq
\mc{O}_X$ is a sheaf of ideals, we   define the
\emph{cosupport} of $\mc{J}$ to be
 $\cosupp \mc{J}=\supp \mc{O}_X/\mc{J}$.  Now consider $\mc{J}=\mcII_1$
 with cosupport  $W=W_1$.  Subdivide $W =\bigcup_{a =1}^d W(a)$  so that
 each $W(a)$ consists of the points
in $W$ contained in a single $\sigma$-orbit and write
$\mc{J}=\prod_a ({}_a\mc{J})$ for the corresponding decomposition of
$\mc{J}$. Let $c=c_0(a)\in W(a)$ be the unique element for which
$W(a)=\{c_j=\sigma^{-j}(c)\}$,   for some
 \emph{positive} set of  integers $j$.  We can  then  (uniquely) write
 ${}_a\mc{J}=\prod{}_{a}^{\hskip 3pt \ell}\hskip -2pt\mc{J}$, where
 ${}_{a}^{\hskip 3pt \ell}\hskip -2pt\mc{J}$ is supported at $c_\ell$,
 and ${}_{a}^{\hskip 3pt \ell}\hskip -2pt\mc{J}=\mc{O}_X$ if $c_\ell\not\in \cosupp \mc{J}$.
The \emph{width} of ${}_a\mc{J}$ is defined to be the maximal $j$
such that $c_j(a)$ appears in $\cosupp {}_a\mc{J}$, and the
\emph{width} of $\mc{J}$ is defined to be $\max \{
\operatorname{width} {}_a\mc{J} \, | \, 1 \leq a \leq d \}$.

 Now take $n\geq 1$ and $r\geq 0$.
 We now repeat the process of the last paragraph
 for $\mcII_n^{\sigma^r}$ and $W_{n,r}=\cosupp \mcII_n^{\sigma^r}$,
 except that we use the elements $c_0(a)$ defined for $\mcII_1$.
 By induction and the equation $\mcII_1\mcII_{n-1}^\sigma\subseteq \mcII_n$
 from Definition~\ref{gen naive def}, the $\sigma$-orbits
 defined by $W_{n,r}$ are contained in those coming from $W_1$.
 Hence each $W_{n,r}=\bigcup_{a=1}^d W_{n,r}(a)$ and, again,
 $W_{n,r}(a)=\{c_j=\sigma^{-j}(c)\}$,   for some
 positive set of  integers $j$. Of course, it is quite possible that
 $c_0(a)\not\in W_{n,r}(a)$ or even  that some $W_{n,r}(a)=\emptyset$;  the extreme case
 occurs  when $\mcII_1\not=\mc{O}_X$ but $\mcII_n=\mc{O}_X$ for all $n>1$.
  A useful observation is that
  \begin{equation}\label{cosupp-trick}
  {}_a^j\bigl(\mcII_n^{\sigma^{r+1}}\bigr)=
  \left({}_{\hskip 6pt a}^{j-1}
  \hskip -1pt\mcII_n^{\sigma^r}\right)^\sigma.
  \end{equation}
\end{notation}

 \begin{lemma}\label{cosupp-lemma}
  Let $\{\mcII_n \}_{n \geq 0}$   be a \gns, fix some    $1 \leq a \leq d$
  as in Notation~\ref{cosupp-defn} and define $t$ as in Definition~\ref{gen naive def}.
  Let $w=\mathrm{width}\, \mcII_1$.
  Then there exist   sheaves of ideals
$\mc{A},\mc{B} ,\mc{C}$, independent of $n$, such that
\begin{equation}\label{cosupp-equ}
   {}_a\mcII_{n} = \mc{A}\, \mc{B}^{\sigma^{w}}\mc{B}^{\sigma^{w+1}}\cdots
 \mc{B}^{\sigma^{n-1}}  \mc{C}^{\sigma^{n}}\qquad\mathrm{for\ all\ } n\geq M=  \max\{w,t\}.
 \end{equation}
Moreover,
  $\mc{A}= {}^0\hskip -2pt\mc{A}\cdot{}^1\hskip -2pt\mc{A}\cdots {}^{w-1}\hskip -2pt\mc{A}$
 and $\mc{C}= {}^0\mc{C}\cdots {}^{w-1}\mc{C}$ but
   $\mc{B}={}^0\mc{B}$.
    If $w=0$,  the sheaves $\mc{A}$ and $\mc{C}$ do not appear.
\end{lemma}

 \begin{proof}   Clearly we can replace
 $\{\mcII_n\}$ by $\{{}_a\mcII_n\}$ and so we will drop the subscript $a$ and put $c_i=c_i(a)$.
Since $\mcII_n \supseteq \mcII_1 \mcII_1^{\sigma} \dots
\mcII_1^{\sigma^{n-1}}$ it follows that $\cosupp \mcII_n \subseteq
\{ c_0, c_1, \dots, c_{n + w -1} \}$.  For $n,r \geq t$,
Definition~\ref{gen naive def}(2) ensures that $^j \mcII_{r+1} {}
^j (\mcII_n^{\sigma^{r+1}}) = {} ^j \mcII_{r+n+1} = {} ^j \mcII_{r}
{} ^j (\mcII_{n+1}^{\sigma^r})$. If $j > r + w$ then $^j\mcII_{r+1}
= {} ^j \mcII_{r} = \mc{O}_X$ and so $^j (\mcII_n^{\sigma^{r+1}}) =
{} ^j (\mcII_{n+1}^{\sigma^r})$. By \eqref{cosupp-trick} this is
equivalent to  $(^{k-1}\mcII_n)^{\sigma} = {} ^k \mcII_{n+1}$ for $k
>w$.

Now consider the equation $^j \mcII_{n+1} {}
^j (\mcII_r^{\sigma^{n+1}})  = {} ^j \mcII_{r+n+1}   =   {} ^j \mcII_{n} {} ^j
(\mcII_{r+1}^{\sigma^n})$, for $n,r\geq t$.  If $j < n$ then
  \eqref{cosupp-trick} implies that
   $^j  (\mcII_{r+1}^{\sigma^n}) = {} ^j (\mcII_{r}^{\sigma^{n+1}}) =
\mc{O}_X$ and so $^j \mcII_{n+1}  = {} ^j \mcII_{n}$.  Altogether,
if $w < j < n $ then $(^{j-1}\mcII_n)^{\sigma} =
{} ^j \mcII_{n+1} = {} ^j \mcII_n$.

Finally, take $n\geq\max\{w,t\}$ so that \eqref{cosupp-equ} makes sense.
The previous paragraph  certainly implies that
$ \mcII_{n} = \mc{A}\, \mc{B}^{\sigma^{w}}\mc{B}^{\sigma^{w+2}}\cdots
 \mc{B}^{\sigma^{n-1}}  \mc{C}^{\sigma^{n}}$,
 where  $\mc{A}=\mc{A}(n)$ and $\mc{C}=\mc{C}(n)$ are supported on
 $\{c_0,\dots, c_{w-1}\}$
 but  $\mc{B}= ( {}^k\mcII_n)^{\sigma^{-k}}$
for any  $w\leq k \leq n-1 $.
 Thus $\mc{B}$ is independent of $n$.  We can certainly  write
$\mc{A}= {}^0\hskip -2pt\mc{A}\cdot{}^1\hskip -2pt\mc{A}\cdots {}^{w-1}\hskip -2pt\mc{A}$
and   this decomposition is independent of $n$ simply because
 $^j \mcII_{n+1}  = {} ^j \mcII_{n}$ for $j<n$. Similarly, the fact that
  $\mc{C}= {}^0\mc{C}\cdots {}^{w-1}\mc{C}$ independently of $n$  follows from the equation
   $(^{k-1}\mcII_n)^{\sigma} = {} ^k \mcII_{n+1}$ for $k
>w$.
 \end{proof}

We now turn to  torsion extensions of graded algebras and bimodule
algebras.     Given  a cg Ore domain $A$ with homogeneous quotient
ring $Q = Q(A)$, the \emph{maximal right torsion extension of $A$}
is the ring
\begin{equation}\label{extn-def}
T(A) = \{x \in Q\ |\ x A_{\geq n} \subseteq A\ \text{for some}\ n
\geq 0 \}  \ \cong \  \lim_{n \to \infty} \Hom_A( (A_{\geq n})_A,
A_A).
\end{equation}
The ring  $T(A)$  is   again a $\mb{Z}$-graded Ore domain with quotient ring
$Q$.   The algebra $A$ is called \emph{right
torsion closed} if $T(A) = A$.  The \emph{maximal  left  torsion extension
of $A$} is defined  analogously and written $T^\ell(A)$.

We also need  the analogues of these definitions for bimodule  algebras.
 Let $\mc{K}$ be the constant sheaf of rational functions on $X$, with the induced
  action of $\sigma$,  and fix once and for all an injection $\mc{L}
\hookrightarrow \mc{K}$; thereby  giving inclusions $\mc{L}_n \subseteq
\mc{K}$ for each $n$.  Given a generalized na{\"\i}ve sequence $\{  \mcII_m \}_{m \geq 0}$,
and integers $n,m\geq 0$,  we  define $
 \mc{H}_n(m)$ to be   the unique largest subsheaf $ \mc{H}\subset
\mc{K} $  such that $\mc{H}\, \mcII_m^{\sigma^n} \subseteq
\mcII_{n+m}.$    Given  subsheaves $\mc{F}, \mc{G} \subseteq \mc{K}$,  we may
identify $\mc{H} = \calHom_{\mc{O}_X}(\mc{F}, \mc{G})$ with  the
unique largest subsheaf of $\mc{K}$ such that $\mc{H} \mc{F}
\subseteq \mc{G}$;  in this way, we have $\mc{H}_n(m) = \calHom(
\mcII_m^{\sigma^n}, \mcII_{n+m})$.

\begin{corollary}\label{cosupp-cor}  Let $\{\mcII_n \}_{n \geq 0}$
be a \gns\ and  keep the notation from
Lemma~\ref{cosupp-lemma}.
\begin{enumerate}
\item For any $n\geq 0$, the $ \mc{H}_n(m)$ are sheaves of ideals
that are equal for $m \geq M=\max\{t,w\}$.
\item For $n,m\geq M$,  one has
$ {}_a\mc{H}_{n}(m) = \mc{A}\, \mc{B}^{\sigma^w}\mc{B}^{\sigma^{w+1}}\cdots
 \mc{B}^{\sigma^{n-1}}  \mc{D}^{\sigma^{n}},$
where $\mc{A}$ and $\mc{B}$ are defined by Lemma~\ref{cosupp-lemma},
while $\mc{D}$ is independent of $n$ and $m$  and satisfies $\mc{C}\subseteq
\mc{D}\subseteq \mc{O}_X$.

\item $\{\mc{H}_n\}_{n \geq 0}$  is a  \gns.
\end{enumerate}
\end{corollary}

\begin{remark}\label{Hn-defn}  Using   Corollary~\ref{cosupp-cor}(1),
we define $\mc{H}_n=\mc{H}_n(m)$ for  any $m\geq \max\{t,w\}$.
\end{remark}

\begin{proof}  As with Lemma~\ref{cosupp-lemma} we can replace $\{\mcII_n\}$ by
$\{ {}_a\mcII_n\}$ and so we can drop the subscript $a$.

(1)   Take  $n \geq 0$ and $m\geq M=\max\{w,t\}$ and identify
$\mc{H}_n(m)=\calHom(\mcII_m^{\sigma^n},\,\mcII_{n+m})$.  By Definition~\ref{gen
naive def}(3), each $c_i$ is a smooth point on a scheme of dimension
$\geq 2$, from which  it follows that $\calHom( \mcII_m^{\sigma^n},
\mcII_{n+m}) \subseteq \calHom( \mcII_m^{\sigma^n}, \mc{O}_X) =
\mc{O}_X$ and so $\mc{H}_n(m)$ is an ideal sheaf.   It is automatic
that $\mcII_n\subseteq \mc{H}_n(m)$ and so
$\cosupp \mc{H}_n(m)\subseteq \cosupp \mcII_n\subseteq \{c_0,\dots,c_{n+w-1}\}$, whence
$^j \mc{H}_n(m) = \mc{O}_X$ unless
$0 \leq j \leq n + w- 1$.  For such $j$  \eqref{cosupp-equ} shows that $^j
\mcII_m^{\sigma^n}$ and $^j \mcII_{n+m}$ are independent of the
choice of $m \geq M$. The result follows.

(2) Throughout the proof we assume that  $n,m \geq M$ and write
$\mc{H}_n=\mc{H}_n(m)$. Thus $\mcII_{n}\mcII_m^{\sigma^{n}} =
\mcII_{n+m}$ and so $\mc{H}_{n}\mcII_m^{\sigma^{n}} = \mcII_{n+m}$.
If $0 \leq j \leq n-1$ then  $^j (\mcII_m^{\sigma^n}) = \mc{O}_X$ by
definition and hence ${}^j\mc{H}_{n}={}^j\mcII_{n+m}$. Combined with
Lemma~\ref{cosupp-lemma} and the fact that $\mc{H}_n$ is an ideal sheaf,
 this implies that $\mc{H}_n=\mc{A}\,
\mc{B}^{\sigma^w}\mc{B}^{\sigma^{w+1}}\cdots
 \mc{B}^{\sigma^{n-1}}  \mc{D}_n^{\sigma^{n}}$, where $\mc{D}_n$ is an ideal sheaf
cosupported on $\{c_0,\dots,c_{w-1}\}$.  Since $\mcII_n \subseteq \mc{H}_n$
 it is clear that $\mc{C}\subseteq \mc{D}_n$.
On the other hand, for $j  >w$, Lemma~\ref{cosupp-lemma} shows
that $^j (\mcII_{m+n}^{\sigma}) = {} ^j \mcII_{m+n+1}$ and so
\[
^j(\mc{H}_n^{\sigma}) = \calHom( {} ^j (\mcII_m^{\sigma^{n+1}}),\,
{}^j (\mcII_{n+m}^{\sigma})) =
 \calHom( {} ^j(\mcII_m^{\sigma^{n+1}}), \, {} ^j \mcII_{n+m +1}) = {} ^j\mc{H}_{n+1}.
\]
Thus $D_n$ is independent of $n \geq M$.

(3)  It follows from  the definition of $\mc{H}_n$    that
$\mc{H}_n\mc{H}_r^{\sigma^n}\subseteq \mc{H}_{n+r}$ for all $n,r\geq
0$. For $n,r\gg 0$,  equality   follows easily by combining parts
(1) and (2).  Since $\mcII_n \subseteq \mc{H}_n \subseteq \mc{O}_X$,
clearly $\cosupp \mc{H}_n$ consists of points lying on critically
dense $\sigma$-orbits and so $\{\mc{H}_n\}_{n \geq 0}$ is indeed a
\gns.
\end{proof}

 It follows from Corollary~\ref{cosupp-cor} that the sequences
  $\{\mc{H}_n\}$ are well-behaved. Indeed we have:

\begin{corollary}\label{cosupp-cor2}
  Let   $\{  \mcII_n \}_{n \geq 0}$ be a generalized na{\"\i}ve  sequence,
  set $\mc{S}=  \mc{S}(X, \{ \mcII_n\},\mc{L}, \sigma) $ with section algebra $S=\HB(X,\mc{S})$
 and define $\mc{H}_n$   by   Remark~\ref{Hn-defn}. Then:
\begin{enumerate}
\item If $\mc{T} = \bigoplus_{n \geq 0} (\mc{H}_n \otimes \mc{L}_n)_{\sigma^n}$,
then
$T = \HB^0(X, \mc{T})$  is  the maximal  right   torsion extension $T(S)$ of $S$.
\item $T(S) = S$ in high degree if and only if $\mc{H}_n=\mcII_n$ for  $n\gg 0$.
\item $T$ is a finitely generated left $S$-module such that
$TS_{\geq n}\subseteq S$ for some $n\geq 1$.
\end{enumerate}
\end{corollary}

\begin{proof}  (1)
By Corollary~\ref{cosupp-cor}(3),  $\mc{T}$ is a $(\mc{O}_X,
\sigma)$-bimodule algebra containing $\mc{S}$.  The given embedding
of $\mc{L}$ in $\mc{K}$ induces an embedding of $\mc{T}$ in the
$(\mc{O}_X, \sigma)$-bimodule algebra $\widetilde{\mc{K}}
 = \bigoplus_{n\in \mathbb{Z}} \mc{K}_{\sigma^n}$
and it follows that $Q(S)\subseteq Q(T) \subseteq K[t,t^{-1};\sigma] \cong
 \mathrm{H}^0(X,\widetilde{\mc{K}})$ (although we will not need it, these inclusions are actually equalities).

Consider the maximal right torsion extension $  T'=T(S)\subset Q(S)$
of $S$. For $n \geq 0$,   $T'_n $ generates a sheaf of
$\mc{O}_X$-modules
 $\mc{T}'_n\subset \mc{K}$, which we  write  as $\mc{T}'_n = \mc{J}_n \mc{L}_n$ for some sheaf
$\mc{J}_n$. The fact that $(T'/S)_S$ is torsion means that, for any given
$n$, one has $T'_n S_m \subset S$ for all $m \gg 0$.  By
Corollary~\ref{gen-naive2}(1),
 $\mc{S}_m$ is generated by its global sections
for all $m\gg 0$  and $\mc{T}'_n$ is generated by
its global sections by definition. Consequently,  $\mc{T}'_n
\mc{S}_m \subseteq \mc{S}$  for $m \gg 0$;
equivalently,  $\mc{J}_n  \mcII_m^{\sigma^n} \subseteq
 \mcII_{n+m}$ for all $m \gg 0$.  Thus
  $\mc{J}_n \subseteq
\mc{H}_n$ and  $\mc{T}' \subseteq \mc{T}$, from which it follows that  $T' \subseteq
\HB^0(X, \mc{T}') \subseteq \HB^0(X, \mc{T}) = T$. Conversely, it is easy
to see  that $(T/S)_S$ is right torsion, and so $T
\subseteq T'$.

(2)  Suppose that $T_n =S_n$ for all $n \gg 0$. By Corollary~\ref{cosupp-cor}(3) and
Corollary~\ref{gen-naive2}(1), $\mc{T}_n = \mc{H}_n \otimes
\mc{L}_n$ and $\mc{S}_n = \mcII_n \otimes \mc{L}_n$ are generated by
their respective global sections $T_n$ and $S_n$ for all $n \gg 0$.
It follows that $\mc{H}_n \otimes \mc{L}_n = \mcII_n \otimes
\mc{L}_n$ and hence $\mc{H}_n = \mc{I}_n$ for all $n \gg 0$.  The
converse follows immediately from part~(1).

(3) As in the proof of  Corollary~\ref{cosupp-cor}(2), for $n ,r
\geq M$ we have $\mc{H}_n \mc{I}_r^{\sigma^n} = \mc{I}_{r+n}$. Thus
$\mc{T}_n \mc{S}_r \subseteq \mc{S}$, and hence $T_n S_r \subseteq
S$ by taking sections. Since $(T/S)_S$ is torsion, for each  $0<m
<M$ we have $T_m S_r \subseteq S$ for all $r \gg 0$. Thus we can
pick a single $r$ such that $T S_{\geq r} \subseteq S$.
\end{proof}

Next, we study  the maximal right torsion extensions  of Veronese
rings.

\begin{lemma}
\label{finite change}
Keep the hypotheses  from Corollary~\ref{cosupp-cor2} and set $T=T(S)$.
  Then, for  $q\geq 1$ one has
\begin{enumerate}
\item If $T=T(S)$, then
   $T^{(q)} = T(S^{(q)})$.
\item $S$ is equal to $T=T(S)$ in large degree if and
only if $S^{(q)}$ is equal to $T^{(q)}$ in large degree.
\end{enumerate}
\end{lemma}

\begin{proof}  (1)  Since $S^{(q)}$ is also a generalized na{\"\i}ve blowup
algebra,  we may apply
Corollary~\ref{cosupp-cor2}(1) to find its maximal right torsion
extension $T(S^{(q)})$. But it is clear from Corollary~\ref{cosupp-cor}(2) that,
for any $n \geq 0$,
$$\mc{H}_{nq} =
\calHom( \mcII_m^{\sigma^{nq}},  \mcII_{m + nq}) = \calHom(
\mcII_{mq}^{\sigma^{nq}},  \mcII_{mq + nq}) \qquad\mathrm{for\ all
}\ m \gg 0.
$$
Thus  $T^{(q)}$ must be the maximal right torsion extension of
$S^{(q)}$.

(2)
Suppose that $S_r\not=T_r$ for infinitely many $r$. Then $\mc{H}_r
\neq \mcII_r$ for infinitely many $r$ and so
Corollary~\ref{cosupp-cor}(2) implies that $\mc{C} \neq \mc{D}$ and
hence $\mc{H}_r \neq \mcII_r$ for all $r \gg 0$.
 In particular,  $\mc{H}_{qr}\not=\mcII_{qr}$ for all $r\gg 0$  and so
   Corollary~\ref{cosupp-cor2}(2)  implies that $S^{(q)}_u=
   S_{qu}\not=T_{qu}=T^{(q)}_u$ for all $u\gg 0$.
 The other direction is trivial.
\end{proof}

 The previous results all have analogs for left torsion
extensions.   Indeed, let   $\{  \mcII_n \}_{n \geq 0}$ be a generalized na{\"\i}ve  sequence,
and set
\begin{equation}
\label{Hdef-left} \mc{H}^{\ell}_n = \text{the unique largest
subsheaf}\ \mc{H}\ \text{ of}\ \mc{K}\ \text{such that}\  \mcII_m
\mc{H}^{\sigma^m} \subseteq \mcII_{n+m},\  \text{for some}\ m \gg 0,
\end{equation}
which is again   independent of the choice of $m \gg 0$.
Now set $\mc{S} = \mc{S}(X, \{
\mcII_n \}, \mc{L}, \sigma)$  with section algebra $S$.
Then \emph{the maximal
left torsion extension $T^{\ell}=T^\ell(S)$} of $S$ may be calculated as the section
algebra  of $\mc{T}^{\ell} = \bigoplus_{n \geq 0}
(\mc{H}^{\ell}_n \otimes \mc{L}_n)_{\sigma^n}$.
These facts will be used without particular comment in the next few examples and their proof is  left  to the
reader.

We end the section with a few more examples to illustrate what can
happen in the passage from a \nsr\ to its maximal right or left
  torsion extension and we keep the notation  developed in Example~\ref{easy-eg1}.
We first note a few more properties of that example:

\noindent
 \begin{example}\label{easy-eg1'} This is a continuation of   Example~\ref{easy-eg1},
 and we use the notation set up there.
Then it is easy to see that $\mc{H}^{\ell}_n= \mcII_n$ for all $n\geq 1$ (use the fact that
$\lcl{\mf{m}}{0}$ is a  maximal  sheaf of  ideals of $\mc{O}_X$).
Hence  $T^{\ell} = R$ by Corollary~\ref{cosupp-cor2}.
Using Theorem~\ref{chi_1, not chi_2}, below,  this   also implies that
 $R$ satisfies $\chi_1$ on the left but not on the right. \qed
\end{example}

It is easy to modify this example so that the $\chi_1$ condition fails on both sides:

\begin{example}\label{easy-eg2}
\emph{In the notation from  Example~\ref{easy-eg1}, Let
$\mc{I}=\lcl{\mc{M}}{0} \lcl{\mathfrak{m}}{1} \lcl{\mc{M}}{2}$.
  Then
$R=R(\mathbb{P}^2,Z_{\mc{I}},\mc{L},\sigma)$ does not satisfy weak
$\chi_1$ on either side.}

 Set $\widehat{\mc{H}}= \lcl{\mc{M}}{0} \lcl{\mathfrak{m}}{1}
\lcl{\mathfrak{m}}{2}^2$. Then a simple computation shows that
$\widehat{\mc{H}}\mcII^\sigma=\lcl{\mc{M}}{0}
\lcl{\mathfrak{m}}{1}^3 \lcl{\mathfrak{m}}{2}^3\lcl{\mc{M}}{3}
=\mcII\mcII^\sigma.$ Thus  $\mc{H} \supseteq \widehat{\mc{H}}$
 (in fact one has $\mc{H} = \widehat{\mc{H}}$).
 On the other hand ${}^{n+1}\widehat{\mc{H}}_n =
\lcl{\mathfrak{m}}{n+1}^2 \supsetneq \lcl{\mc{M}}{n+1} =
{}^{n+1}\mcII_{n+1}$. By Corollary~\ref{cosupp-cor2}, this implies
that  $T(R)_n\not= R_n$ for all $n\gg 0$. By symmetry,
$T^{\ell}(R)_n\not= R_n$ for  $n\gg 0$.  \qed
  \end{example}

In Example~\ref{easy-eg2},  it happens that passing to the maximal
torsion extension on one side and then the other  leads to the same
ring; indeed
$T(T^\ell(R))=T^\ell(T(R))=R(\mathbb{P}^2,Z_{\mc{J}},\mc{L},\sigma)$,
  where  $\mc{J}=\lcl{\mathfrak{m}}{0}^2 \lcl{\mathfrak{m}}{1}
\lcl{\mathfrak{m}}{2}^2$. However, as the next example shows, this
does not always happen.

\begin{example}\label{easy-eg3} \emph{There exists a \nsr\ $R$ such that
 $T(R)$ and $T^\ell(R)$ are distinct  infinite dimensional extensions of $R$,
yet both $T(R)$ and $T^\ell(R)$ are left and right torsion closed.}

In the proof we again use the notation from Example~\ref{easy-eg1}.
We first seek ideals primary to $(x,y)$   with the following properties:
$P\supsetneq I\supsetneq K\subsetneq J\subsetneq Q$ and $PJ=IJ=IQ=K$, but
$$P=(K:J)=\{  r\in k[x,y] : rJ\subseteq K\},$$ $J = (K:P)$, $Q = (K:I)$, and $I = (K:Q)$.
 An easy calculation shows that the following ideals
have the desired properties.
 $$J= (x^6, x^5y, xy^5, y^6) + (x,y)^7\ \subset\  Q=J+(x^3y^3),\qquad
I = (x^4y^2, x^2y^4) + (x,y)^7\ \subset  \ P=I+(x^3y^3)  $$ and
$$
K= (x^{10}y^2, x^9y^3, x^8 y^4, x^7 y^5, x^5 y^7, x^4y^8, x^3 y^9,
x^2 y^{10}) + (x,y)^{13}.$$ (A crucial property of $K$ is that
$x^6y^6\not\in K$).  As before, for $N=P,Q,$ etc, write $\mc{N}$ for
the sheaf of ideals on $\mathbb{P}^2$ that equals $N$  locally at
$c_0$ and is cosupported at $c_0$, and put $\lcl{\mc{N}}{i} =
\mc{N}^{\sigma^i}$.

Now take $R=R(\mathbb{P}^2, Z_{\mc{F}},\mc{L},\sigma)$, where
$\mc{L}=\mc{O}_{\mathbb{P}^2}(m)$ for suitably large $m$   and
$\mc{F}= \lcl{\mc{I}}{0} \lcl{\mc{J}}{1}$. Either by direct
computation or using Corollary~\ref{cosupp-cor}, one shows using $Q
= (K:I)$ that $T=R(\mathbb{P}^2, Z_{\mc{G}},\mc{L},\sigma)$ for
$\mc{G}= \lcl{\mc{I}}{0} \lcl{\mc{Q}}{1}$, while a similar
left-sided computation using $P = (K:J)$ gives
$T^\ell=R(\mathbb{P}^2, Z_{\mc{H}},\mc{L},\sigma)$ for $\mc{H}=
\lcl{\mc{P}}{0} \lcl{\mc{J}}{1}$. By Corollary~\ref{cosupp-cor}, the
fact that $J=(K:P)$ shows that $T(T^\ell(R))=T^\ell$, while
similarly $I=(K:Q)$ implies that $T^{\ell}(T(R))=T$. \qed
\end{example}

To end the section,  we give the promised  example where the    sequence
$\{\mc{H}_n\}$ arising from a \naive\ sequence $\{\mcII_n\}$   is not  a
   \naive\ sequence. Thus   one does   need  the
theory of generalized \naive\ sequences.

\begin{example}\label{not-naive}
Keep the notation introduced in Example~\ref{easy-eg1}.    Define
the ideals $M = (x^2, y^2)$, $N = (x^3, y^3)$, $F = N + (x,y)^4$, and $G
= (x^4, x^3y, xy^3, y^4)$.
As usual given $P=M,N,$ etcetera,  let  $ \mc{P}$  be the sheaf of ideals
with cosupport $c_0$ that equals $P$ locally at $c_0$.
Now take $\mc{I}=\lcl{\mc{M}}{0} \lcl{\mc{N}}{1}
\lcl{\mathfrak{m}}{2}$ and consider
$\mc{R}=\mc{R}(\mathbb{P}^2,Z_{\mc{I}},\mc{L},\sigma)$ with maximal
right torsion extension $\mc{T}=\bigoplus
(\mc{H}_n\otimes\mc{L}_n)_{\sigma^n}$.

One computes that $\lcl{\mc{M}}{0}\lcl{\mc{N}}{0}\lcl{\mathfrak{m}}{0}
=\lcl{\mathfrak{m}}{0}^5 =\lcl{\mc{M}}{0}\lcl{\mc{F}}{0}\lcl{\mathfrak{m}}{0}$,
and that  this sheaf of ideals equals $\mc{B}$ in the notation of
Lemma~\ref{cosupp-lemma}. With the help of Lemma~\ref{cosupp-lemma}
and Corollary~\ref{cosupp-cor},  one then calculates   that
${}\mc{H}_1 = \lcl{\mc{M}}{0} \lcl{\mc{F}}{1}
\lcl{\mathfrak{m}}{2}$.   It follows that
$^2 (\mc{H}_1\mc{H}_1^\sigma) = \lcl{\mc{G}}{2}$, whereas $^2
\mc{H}_2$ is equal to $\lcl{\mathfrak{m}}{2}^4$. Thus
$\mc{H}_2\not=\mc{H}_1\mc{H}_1^\sigma $.
\qed
 \end{example}

%%%%%%%%%%%%%%%%%%%%%%%%%%%%%%%%

\section{The $\chi$ conditions and cohomological conditions}
\label{chi}

The hypotheses from Assumptions~\ref{global-convention2} remain in
force in this section and we first define the $\chi$ conditions.
Let $A$ be a
 cg $k$-algebra and identify $k$ with the factor ring $k=A/A_{\geq 1}$.
For $n\geq 1$, we say that  \emph{$A$ satisfies   $\chi_n$ on the right}
 if $\dim\ext^i_{\rmod A}(k,M)<\infty$ for all finitely generated graded
 right $A$-modules $M$ and all $i\leq n$.
It is immediate that a ring satisfying $\chi_1$ also satisfies weak $\chi_1$.

  As was shown in Section~\ref{gen naive blowups},
the \nsr\ $R=R(X,Z,\mc{L},\sigma)$ need not satisfy even weak
$\chi_1$ when one \naive ly blows up more than one point and this is
in marked contrast to  the case of blowing up a single point, where
$\chi_1$ always holds   \cite[Theorem~1.1(8)]{KRS}.
  In this section we continue our study of the $\chi$ conditions,
  showing in particular that the maximal right torsion extension
of $R$ will satisfy $\chi_1$ on the right.  On the other hand the
higher $\chi$ conditions behave the same way whether one  \naive ly
blows up one or more than one point---they always fail. We will also
want to consider the $\chi_1$ condition at the level of individual
modules for which we need another definition. Recall that $\pi: \rGr
R\to \rQgr R$ is the natural morphism. For  $N \in \rgr R$, we say
that \emph{the condition $\chi_1(N)$ holds} provided that the
natural map
\begin{equation}\label{chi-first}
N \to \lim_{n \to \infty} \Hom_R((R_{\geq n})_R, N_R) = \bigoplus_m
\Hom_{\rQgr R}(\pi(R), \pi(N)[m])
\end{equation}
has a right bounded  cokernel, as defined on
page~\pageref{bounded-defn}. By \cite[Propositions~3.11(2)
and~3.14(2)]{AZ1}, $R$ satisfies $\chi_1$ if and only if condition
$\chi_1(N)$ holds for all $N\in \rgr R$.

The main result of this section is the following theorem, which also
proves part~(8) of  Theorem~\ref{mainthm}.

\begin{theorem}
\label{chi_1, not chi_2} Keep the hypotheses from
\eqref{global-convention2}, let $R=R(X,Z,\mc{L},\sigma)$   and write $T=T(R)$
for the maximal right torsion extension of $R$. Then:
\begin{enumerate}
\item Condition $\chi_1(N)$ holds for all Goldie torsion modules $N\in \rgr R$.
\item
$\chi_1$ holds for $R$ on the right if and only if $R = T$ in large
degree.

\item If $R$ is a nontrivial \nsr, then $\chi_2$ fails for $R$ on   the right; indeed
$\ext_{\rcatMod R}^2(k,R)$ is infinite dimensional.

\item If $R$ is nontrivial, then ${\coH}^1(\pi(R))=\ext^1_{\rQgr R}(\pi(R),\pi(R))$ is infinite
dimensional.  \end{enumerate}
\end{theorem}

\begin{proof}
(1) Fix  $N \in \GT \rgr R$.  We convert \eqref{chi-first} into a
statement about the $\mc{R}$-module $\mc{N}=N\otimes_R\mc{R}$. By
the equivalence of categories, Theorem~\ref{general-ampleness},
$\xi^{-1}\circ\pi_R(N)=\pi_{\mc{R}}(\mc{N})$ and
$N=\mathrm{H}^0(X,\, \mc{N})$ in high degree.   Thus
\eqref{chi-first} has a right bounded cokernel if and only if the
morphism
\begin{equation}
\label{max tors ext} \mathrm{H}^0(X,\, {\mc{N}}) \to \bigoplus_m \Hom_{\rQgr
\mc{R}}(\pi(\mc{R}), \pi(\mc{N})[m])
\end{equation}
does. Since $N\otimes_R\mc{R}$ is Goldie torsion,
it will suffice to show that \eqref{max tors ext} has right bounded
cokernel for an aribtrary Goldie torsion module $\mc{N}\in \GT\rgr
\mc{R}$.

Write $\mc{N} = \bigoplus  \mc{F}_n\otimes  \mc{L}_{\sigma}^{\otimes
n}$  and $\mc{R} = \bigoplus \mc{I}_n \otimes
\mc{L}_{\sigma}^{\otimes n}$. For $n \geq n_0\gg 0$,
Lemma~\ref{Goldie tors facts} implies that $\mc{F}_n =
\mc{F}_{n+1}=\mc{F}$, say. Now fix some $m\geq n_0$ and write
$\mc{N}[m] = \bigoplus \mc{G}_n \otimes \mc{L}_{\sigma}^{\otimes
n}$; thus
$\mc{G}_n=\mc{G}_{n+1}=\mc{G}=(\mc{F}\otimes\mc{L}_m)^{\sigma^{-m}}$
for all $n\geq 0$, by Lemma~\ref{shifts}.  By Lemma~\ref{Hom in
Qgr}(1) we have an isomorphism
$$
\Hom_{\rQgr \mc{R}}(\pi(\mc{R}), \pi(\mc{N})[m]) \cong \lim_{n \to
\infty} \Hom_{\mc{O}_X}(  \mc{I}_n, \mc{G}_n) = \lim_{n \to \infty}
\Hom_{\mc{O}_X}(  \mc{I}_n, \mc{G}).
$$
If we can show, for $m \gg 0$, that the maps in the direct limit
$\lim_{n \to \infty} \Hom_{\mc{O}_X}(  \mc{I}_n, \mc{G})$ are
isomorphisms for all $n \geq 0$, then we are done, since the zeroth
term of the limit is nothing more than $\hom_{\mc{O}_X}(  \mc{I}_{0},
\mc{G}_{0})=\coH^0(X,\mc{G}_0) \cong\coH^0(X,\mc{F}_m\otimes
\mc{L}_m)= \mathrm{H}^0(X,\, \mc{N})_m$.

Since $\mc{F}$ is a Goldie torsion sheaf, its support $\supp \mc{F}$  is a
proper closed subset of $X$.
As $Z$ is saturating,  each point in $S=\supp Z$ lies on a
 critically dense $\sigma$-orbit and so
 $\sigma^{-m}(S)\cap \supp \mc{F}=\emptyset$ for all $m\gg 0$.
 Since $\supp \mc{G}=\sigma^m(\supp \mc{F})$,
  we can therefore choose $m\gg 0$ such that
 $\supp \mc{G}\cap \sigma^{-j}(S)=\emptyset$ for all $j\geq 0$.
Now since $\supp   \mc{I}_n/  \mc{I}_{n+1} \subseteq \supp
\mc{O}_X/\mc{I}^{\sigma^n} \subseteq \sigma^{-n}(S)$,
\cite[Lemma~7.2(1)]{KRS}    implies that
$\Hom(\mc{I}_n/\mc{I}_{n+1}, \mc{G}) = 0 =
\ext^1(\mc{I}_n/\mc{I}_{n+1}, \mc{G})$ for all $n \geq 0$ for such
large $m$.  This implies that the $n$th map of the direct limit
$\lim_{n \to \infty} \Hom_{\mc{O}_X}(\mc{I}_n, \mc{G})$ is an
isomorphism for all $n \geq 0$ as we needed.

(2) By \cite[Proposition~3.14]{AZ1} we need to prove that $\chi_1(N)$
 holds for all modules $N\in \rgr R$.
This condition clearly holds for $N$ if and only if it holds for a
shift $N[r]$. Since $N$ has a filtration by shifts of $R$ and Goldie
torsion modules, it suffices to prove the condition in those two
cases. When $N$ is Goldie torsion, the result is given by part~(1), so assume that
 $N = R$.  Then $\lim_{n \to \infty} \Hom_R(R_{\geq n}, N)$ is
simply the maximal right torsion extension of $R$, namely $T$.  Thus
in this case the condition demanded by \eqref{chi-first} is precisely that $T$ and
$R$ are equal in large degree.

(3)  By \cite[$(\dagger)$, p.~274]{AZ1}, it suffices to prove that
 $\dim_k\ext^1_{\rQgr R}(\pi(R),\pi(R))=\infty$.  Thus (3)  follows from~(4).

(4) By Theorem~\ref{general-ampleness}, again,
 this is equivalent to showing  that
$
\dim_k\ext^1_{\rQgr \calR}(\pi(\calR) , \pi(\calR))=\infty.
$
Write $\mc{B} = \bigoplus_{n \geq 0} \mc{L}_{\sigma}^{\otimes n}$,
which we think of as a right $\mc{R}$-module. The long exact
sequence in $\Hom$ induced from the inclusion $\pi(\mc{R})\subset
\pi(\mc{B})$ provides the following exact sequence:
\begin{equation}
\label{chitwo}
 \Hom_{\rQgr \mc{R}}(\pi(\mc{R}), \pi(\mc{B}))
 \overset{\phi_1}{\longrightarrow}
\Hom_{\rQgr \mc{R}}(\pi(\mc{R}), \pi(\mc{B}/\mc{R}))
\overset{\phi_2}{\longrightarrow} \ext^1_{\rQgr \mc{R}}(\pi(\mc{R}),
\pi(\mc{R})).
\end{equation}
We need to understand the first two terms in \eqref{chitwo}. As in
Example~\ref{definable}(1),  we see that $\mc{B}$
is definable by products in all degrees $\geq 0$, and thus
Lemma~\ref{Hom in Qgr}(1) implies that $\Hom_{\rQgr
\mc{R}}(\pi(\mc{R}), \pi(\mc{B})) \cong {\displaystyle \lim_{n\to
\infty}} \Hom_{\mc{O}_X}(\mc{I}_n, \mc{O}_X)$. By
\cite[Lemma~7.2(2)]{KRS}, $\ext^1(\mc{O}_X/\mc{I}_n, \mc{O}_X) = 0$
and so
 the natural map $\Hom(\mc{O}_X, \mc{O}_X) \to \Hom(\mc{I}_n, \mc{O}_X)$ is an
isomorphism for all $n \geq 0$. Thus $\Hom_{\rQgr
\mc{R}}(\pi(\mc{R}), \pi(\mc{B})) \cong \Hom_{\mc{O}_X}(\mc{O}_X,
\mc{O}_X) = k$.

On the other hand, $\mc{B}/\mc{R}$ will generally not be definable
by products  and so we cannot use Lemma~\ref{Hom in Qgr}(1) to
calculate $\Hom_{\rQgr \mc{R}}(\pi(\mc{R}), \pi(\mc{B}/\mc{R}))$
directly.  Instead we will examine a simpler submodule of
$\mc{B}/\mc{R}$.   By   assumption,  $Z \neq \emptyset$  and so   we
may choose some $c \in \supp Z$ and write $c_i = \sigma^{-i}(c)$ for
$i \in \mb{Z}$. More than one point of $\supp Z$ might lie on the
same $\sigma$-orbit as $c$, but we can choose the point $c\in \supp
Z$ so that $c_i \not \in \supp Z$ for all $i < 0$.
 Since  $c_{n-1} \in \supp \mc{O}_X/\mc{I}_n$  for   $n\geq 1$,   there exists
 an ideal sheaf $\mc{J}_n$ such that $\mc{J}_n/\mc{I}_n \cong
k(c_{n-1})$.  Set $\mc{M}_n = \sum_{1 \leq i \leq n}
\mc{J}_i\mc{I}_{n-i}^{\sigma^i}$, and notice that for each $i$ we
 have $\mc{J}_i \mc{I}_{n-i}^{\sigma^i}/\mc{I}_n \cong
k(c_{i-1})$, simply because $c_{i-1} \not \in \supp
\mc{O}_X/\mc{I}_{n-i}^{\sigma^i}$.  Consequently,
 for each $n$ we have $\mc{M}_n/\mc{I}_n\cong \bigoplus_{i = 1}^n
\mc{J}_i\mc{I}_{n-i}^{\sigma^i}/\mc{I}_n \cong  \bigoplus_{i = 1}^n  k(c_{i-1})$.
 Now for  each $i \geq 1$ we can define a right $\mc{R}$-module $\mc{N}^{(i)}
= \bigoplus_{n \geq i} \mc{J}_i\mc{I}_{n-i}^{\sigma^i}/\mc{I}_n
\otimes \mc{L}_{\sigma}^{\otimes n}$ in the obvious way, and it is
clear that $\mc{N}^{(i)} \cong \overline{c}_{i-1}\,[-i]$ in the notation
from before Corollary~\ref{GT equiv1}.  Moreover, $\mc{N}^{(i)}$ is definable by
products in degrees $\geq i$.  Then $\bigoplus_{n \geq 1}
\mc{M}_n/\mc{I}_n \otimes \mc{L}_{\sigma}^{\otimes n}$ is a
submodule of $\mc{B}/\mc{R}$ which is isomorphic to $\bigoplus_{i
\geq 1} \mc{N}^{(i)}$.

Applying  Lemma~\ref{Hom in Qgr}(1) to $\mc{N}^{(i)}$ shows that
$$\Hom_{\rQgr \mc{R}}(\pi(\mc{R}), \pi(\mc{N}^{(i)})) = \lim_{n \to
\infty} \Hom_{\mc{O}_X}(\mc{I}_n, k(c_{i-1})) = \lim_{n \to \infty}
\Hom_{\mc{O}_{X,c_{i-1}}}( (\mc{I}_n)_{c_{i-1}}, k(c_{i-1})).$$
 The final
direct limit clearly stabilizes for $n \gg 0$ to something nonzero.
Then $\Hom_{\rqgr \mc{R}}(\pi(\mc{R}), \bigoplus_{i \geq 1}
\pi(\mc{N}^{(i)}))$ is infinite-dimensional over $k$.  By
left-exactness of $\Hom$, we see that $\dim_k \Hom_{\rQgr
\mc{R}}(\pi(\mc{R}), \pi(\mc{B}/\mc{R})) = \infty$ as well.  Then
the map $\phi_2$ in \eqref{chitwo} has an infinite dimensional
cokernel, and we are done.
\end{proof}

The characterization of $\chi_1$ given in the theorem easily
extends to the case of generalized na{\"\i}ve blowups, as the following corollary shows.

\begin{corollary}
\label{fixing chi} Let  $S=S(X,\{\mcII_n\},\mc{L},\sigma)$ be a  \gnba\ and write $T=T(S)$
 for the maximal right  torsion extension of $S$.  Then:
  \begin{enumerate}
\item $T$ is right torsion closed and satisfies right $\chi_1$.
\item $S$ satisfies right $\chi_1$ if and only if $S$ equals $T$ in
large degree.
\item If $S$ is left torsion closed, then so is $T$.
  If $S$ satisfies left $\chi_1$ then so does  $T$.
\item If $S$ is a nontrivial \gnba\ then
$\chi_2$ fails for $S$.
 \end{enumerate}
\end{corollary}

\begin{proof}
(1)   Let $x \in Q(T)=Q(S)$ be a homogeneous element with  $x T_{\geq n} \subseteq
T$ for some $n$. For  any homogeneous element  $z \in T_{\geq n}$, we have
  $xz S_{\geq m} \subseteq S$ for some $m$
  depending on $z$.  By Corollary~\ref{gen-naive2},  $S$ is a cg
noetherian algebra and so   it is   finitely generated  as an  algebra, say
 in degrees $\leq d$. Set  $V =(\bigoplus_{i=n}^{n+d} S_i)$.
 Then there exists a single $m$ such that $(x V) S_{\geq m} \subseteq S$. Since
  $V S_{\geq m} \supseteq S_{\geq (m+n+d)}$, this implies that $x S_{\geq
(m+n+d)} \subseteq S$. Thus  $x \in T$ and  $T$ is right torsion
closed.
  Combining  Corollary~\ref{cosupp-cor2}(1), Corollary~\ref{cosupp-cor}(3) and
   Lemma~\ref{gen-naive1} shows  that $T^{(q)}$ is a \nsr\ for $q \gg 0$.
 By  Lemma~\ref{finite change},   $T^{(q)}$ is also right torsion closed and
 so    Theorem~\ref{chi_1, not chi_2}(2)
shows  that $T^{(q)}$  satisfies right $\chi_1$.
By   Corollary~\ref{gen-naive2}(3) and
\cite[Theorem 8.3(1)]{AZ1}  $T$ satisfies right $\chi_1$.

(2) If $S$ does not equal $T$ in large degree, then $T$ is an
infinite dimensional right torsion extension of $S$ and so  $S$ fails right
$\chi_1$.  If  $S$ does equal $T$ in large degree, then    $S$   satisfies right
$\chi_1$ by part~(1) combined with  \cite[Lemma~8.2(5)]{AZ1}.

(3)  Suppose that $S$ is left torsion closed and let
  $x\in Q(S)=Q(T)$ be a  homogeneous element  such that
 $T_{\geq m} x  \subseteq T$ for some $m\geq 0 $.  By
Corollary~\ref{cosupp-cor2}(3), $T_{\geq m} x S_{\geq n} \subseteq
TS_{\geq n}\subseteq S$  for some $n$, so in
particular $S_{\geq m} x S_{\geq n} \subseteq S$. Since $S$ is left
torsion closed, this implies that $x S_{\geq n} \subseteq S$ and
hence that $x \in T$. So $T$ is also left torsion closed.

Assume that $S$ satisfies left $\chi_1$.  Then
$\dim_kT^\ell(S)/S<\infty$, whence $T^\ell(S)\subseteq T$.   By
\cite[Lemma~8.2(5)]{AZ1} $T^\ell(S)$ satisfies left
 $\chi_1$ and so, replacing $S$ by   $T^\ell(S)$, we can
assume that $S$  is left torsion closed. Then
     $T$ is  left torsion closed by the   last paragraph and satisfies left $\chi_1$
  by the left hand analogue of part~(2).

(4)  By Lemma~\ref{gen-naive1}, some Veronese ring $S^{(p)}$ is a
\nsr\ and it is obviously still nontrivial.   By the theorem, $\chi_2$ fails for $S^{(p)}$ and
so, by the proof of \cite[Proposition~8.7]{AZ1}, it also fails for
$S$.
\end{proof}

\begin{remark}\label{artin-zhang3}  The significance of the $\chi_1$ condition is that
it allows one to recover the ring $R$ from $\rqgr R$ and to apply
the results from \cite{AZ1}. For  example, suppose that $R = R(X, Z,
\mc{L}, \sigma)$ satisfies
Assumptions~\ref{global-convention2} and that    $ R= T(R)$ in large
degree.  Then it follows from Theorem~\ref{general-ampleness} and
\cite[Theorem~4.5(2)]{AZ1} that \emph{$R$ is equal in large degree
to $  \bigoplus_{n\geq 0} \hom_{\rqgr R}( \pi(R), \pi(R)[n])$.}
See \cite[pp.~528-9]{KRS} for a further discussion.
\end{remark}

Let  $R = R(X, Z, \mc{L}, \sigma)$ satisfy
Assumptions~\ref{global-convention2}. As Example \ref{easy-eg2} and
Example \ref{easy-eg3} show,  $R$ may well fail $\chi_1$ on both
sides. However Corollary~\ref{fixing chi} shows that we can repair
this failure without greatly changing the  properties of $R$.
 Specifically,   apply Corollary~\ref{fixing chi} to $R$   on the right to give
the algebra  $T=T(R)$   that satisfies  right $\chi_1$.
  Then   apply the left-sided analogue of this construction  to
$T$, giving an algebra $U=T^\ell(T)$ which,  by Corollary~\ref{fixing
chi},  will satisfy $\chi_1$ on both sides. In terms of noncommutative geometry these
 operations are fairly innocuous.
Indeed, recall that  the
noncommutative projective scheme $\rproj S$ for a cg $k$-algebra $S$ is defined to be the
  the pair $(\rqgr S, \pi(S))$.   By  \cite[Proposition 2.7]{SZ2},
  one has $\rproj R \simeq \rproj T$. The same is not quite true in passing from
  $T $ to $U$, although by mimicking the proof of
  \cite[Lemma~3.2]{Ro2}, one can show that
 $(\rqgr U, \pi(I_T)) \simeq \rproj T$ for an appropriate module $I$;  thus   the underlying category
  will be the same,   although the distinguished object may change.

Such phenomena occur
elsewhere in noncommutative geometry; for example, by
\cite[Lemma~2.2(iii)]{SZ2} they occur for the idealizer ring  $R$
from \cite[Theorem~2.3]{SZ2}. In fact,  using  the following observation,
we can interpret a number of our examples as idealizer rings.

\begin{lemma}
\label{ideal-lem} Suppose that the cg Ore domain $S$ is left torsion
closed in $Q(S)$ and that its maximal right torsion extension $T = T(S)$,
satisfies $T S_{\geq n} \subseteq S$ for some $n$.  Then $S$ is the
\emph{idealizer}  $S=\{\theta\in Q(S) : I\theta\subseteq I\}$
   of the left ideal $I = TS_{\geq n}$ of $T$.
\end{lemma}

\begin{proof}
Obviously $IS \subseteq I$.  Conversely, suppose that  $x \in Q(S)$ is
homogeneous and $Ix \subseteq I$.  Then $S_{\geq n}x \subseteq
TS_{\geq n} x \subseteq T S_{\geq n}$, and so $x$ is in the maximal
left torsion extension of $S$, namely $S$.
\end{proof}

Given a generalized \naive\ blowup algebra $S$, then $T^\ell (S)$ is
left torsion closed by the left-hand analogue
 Corollary~\ref{fixing chi} and so, by Lemma~\ref{ideal-lem},
$T^\ell (S)$ is an idealizer subring of $T(T^\ell(S))$.
For example,     the ring $R$ from
Example~\ref{easy-eg1} satisfies $T^\ell(R)=R$ (see Example~\ref{easy-eg1'}) and so $R$
is an idealizer ring inside $T=T(R)$; indeed, we  showed that $\mc{T}\mc{R}_{\geq 1}
\subseteq \mc{R}$ from which it follows that $R$ is actually the idealizer of
$I=TR_{\geq 1}=R_{\geq 1}$.

There is a curious contrast between these examples and
 earlier appearances of idealizer domains in  noncommutative geometry in
\cite{AS, Ro2, SZ2}.   Those  earlier  examples all
 have the property that no Veronese ring is generated in
degree $1$ (see \cite[Proposition~6.6]{AS},
\cite[Theorem~8.2(6)]{Ro2} and \cite[Corollary~3.2]{SZ2}). In
contrast, by Lemma~\ref{gen-naive1} and Proposition~\ref{veronese},
any idealizer ring $S$ which is also a generalized \naive\ blowup
algebra  will always have some  Veronese ring $S^{(q)}$ that is
generated in degree $1$.  In particular, by replacing some such
idealizer  $S$ by $S^{(q)}$ one obtains an example of an idealizer
which is a cg domain generated in degree~$1$.

Let $S=S(X, \{\mcII_n\},\mc{L},\sigma)$ be a nontrivial \gnba.
 We end the section by
studying the homological and cohomological dimensions of
 $\rQgr S\simeq \rQgr\mc{S}$.  Here,
the \emph{global dimension} of $\rQgr \mc{S}$ (or $\rQgr S$) is
 defined to be
$ \gld( \rQgr \mc{S}) = \sup \{i \mid \ext^i_{\rQgr \mc{S}} (\mc{M},
\mc{N}) \neq 0\ \text{for some}\ \mc{M}, \mc{N} \in \rQgr \mc{S} \}.
$ The \emph{cohomological dimension} of $\rQgr \mc{S}$ (and $\rQgr
S$) is $\cd(\rQgr\mc{S}) = \sup \{\cd(\mc{N})\mid  \mc{N}\in \rQgr
\mc{S}\},$ where $\cd(\mc{N}) =
 \sup \{i \mid  \ext^i_{\rQgr \mc{S}}
(\mc{S}, \mc{N})  \neq 0.\} $

Before stating the theorem, we need the following lemma.

\begin{lemma}\label{nonrep-lemma} {\rm (1)}  Let $S$ be a cg noetherian domain
 such that, for some $t\geq 0$, one has $S_nS_m=S_{n+m}$ for all $n,m\geq t$.
 Then the
Veronese ring $S^{(t)}$ is noetherian and $\rqgr S \simeq \rqgr S^{(t)}$
via the functor $M\mapsto M^{(t)}$.

{\rm (2)} Let $S=S(X,\{\mcII_n\},\mc{L},\sigma)$ be a \gnba.  Then
there exists a \nsr\ $R$
 such that $R=T(R)$ is generated in degree one, satisfies $\chi_1$,
and has $\rqgr R \simeq \rqgr S$.
\end{lemma}

\begin{proof}  (1)
This is similar to the proof of \cite[Proposition~6.1]{AS}. The ring
$S^{(t)}$ is noetherian by \cite[Proposition~5.10]{AZ1}.  Thus,
given $M\in \rgr S$, then $M^{(t)}\in \rgr S^{(t)}$ and if $N\in
\rgr S^{(t)}$ then $N\otimes_{S^{(t)}}S\in \rgr S$. Clearly
$(N\otimes S)^{(t)}=N$, so consider the kernel and cokernel of the
natural map from $M^{(t)}\otimes S \to M$. Either module $L$ is a
noetherian $S$-module  satisfying $L^{(t)}=0.$ We claim that $L$ is
bounded.

 If $L$ is not bounded, pick some $a\in L_r$ such that $ aS$ is infinite dimensional.
 Pick $u\in \mathbb{N}$ with $ut>r$.
For all $m\geq 0$ we have $aS_{mt+(ut-r)} \in L^{(t)}= 0 $ whence
$0= aS_{mt+(ut-r)}S_v = aS_w$ for all $v\geq t$ and $w=v+(mt+ut-r)$.
But all integers $w\gg 0$ can be so written, implying that $L_w=0$
for all $w\gg 0$. Thus $L$ is indeed bounded. It follows routinely
that the maps    $M\mapsto
M^{(t)}$ and $N\mapsto N\otimes_{S^{(t)}} S$ define the equivalence
between $\rqgr S$ and $ \rqgr S^{(t)}$.

(2)  By Lemma~\ref{cosupp-cor2} and Corollary~\ref{cosupp-cor}(3),
$T=T(S)$ is a \gnba\ and, as mentioned after
Remark~\ref{artin-zhang3}, $\rqgr S\simeq \rqgr T$ follows from
\cite[Proposition~2.7]{SZ2}. By Lemma~\ref{gen-naive1}  and
Proposition~\ref{veronese}, for some $q\gg 0$ the ring $T^{(q)}$ is a \nsr\ that is
generated in degree one. By part~(1) and Lemma~\ref{gen-naive2}(4),
$\rqgr S \simeq \rqgr T^{(q)}$ for such $q$. Finally, by
Lemma~\ref{finite change}, $T^{(q)}$ is right torsion-closed and so
Theorem~\ref{chi_1, not chi_2} implies that $R=T^{(q)} $ satisfies
$\chi_1$.
 \end{proof}

\begin{theorem}\label{fin-cohom-dim} Let $S=S(X, \{\mcII_n\},\mc{L},\sigma)$ be a
nontrivial \gnba. Then one has
 $ \cd(\rQgr S)\leq \dim X$.
 If  $X$ is smooth, then
   $  \gld(\rQgr S)\leq 1+\dim X$.
 \end{theorem}

\begin{remark} This result proves
  Theorem~\ref{mainthm}(7)   from the introduction.
\end{remark}

 \begin{proof} By Corollary~\ref{gen-naive2}(4) and Lemma~\ref{nonrep-lemma}
we can replace $S$ by some large Veronese ring $S^{(p)}$ and so, by
Lemma~\ref{gen-naive1}, assume that
$S=R(X,Z_{\mc{I}},\mc{L},\sigma)$ is a \nsr.

The proof of  the corresponding assertions in  \cite[Theorem~8.2 and
Corollary~8.3]{KRS} now  go through with the following minor
changes.   First, the statement and proof of
  \cite[Lemma~7.2]{KRS} go through unchanged using the
  definition of $\mc{I}_n$ from this paper.  Then one should   replace,
  in order of their appearance,
\cite[Theorem~4.1]{KRS} by  Theorem~\ref{general-ampleness};
\cite[Lemma~6.1]{KRS} by Lemma~\ref{Goldie tors facts};
\cite[Lemma~6.4]{KRS} by   Lemma~\ref{Hom in Qgr}(1) and
  Example~\ref{definable};
\cite[Theorem~6.7]{KRS} by Theorem~\ref{GT equiv};
 finally,   \cite[Lemma~6.2]{KRS} by Lemma~\ref{goldie-subfactors}.
\end{proof}

In fact, one  can prove that $\dim X-1 \leq \cd(\rQgr S)\leq \dim X$
and $\dim X\leq \gld(\rQgr S)\leq 1+\dim X$  in
Theorem~\ref{fin-cohom-dim}.   A
 detailed proof of this assertion  can be found in \cite{RS2},
but we will not give it here, in part
 because   we conjecture that the correct dimension is $\dim X$ in both cases.
In the commutative case, and in contrast to  Theorem~\ref{fin-cohom-dim}(2),
  one can easily blow up a nonsingular integral scheme
at a zero-dimensional subscheme and obtain a  scheme that is
singular (see, for example, \cite[Section~IV.2.3]{EH}).

\begin{remark}\label{dualizing-remark}
As a final application of Theorem~\ref{chi_1, not chi_2}  note
that, by \cite[Theorem~4.2]{YZ}, Theorem~\ref{chi_1, not chi_2}
implies that a nontrivial \nsr\   $R$
does not have a balanced dualizing complex, in the sense of
Yekutieli \cite{Ye}. By \cite{Jg} and Theorem~\ref{fin-cohom-dim},
it does however have a dualizing complex in the weaker sense of~\cite{Jg}.
\end{remark}

%%%%%%%%%%%%%%%%%%%%%%%%%%%%%%%%

\section{Generic flatness and parametrization}\label{further-section}

The hypotheses from Assumptions~\ref{global-convention2} will be
assumed throughout this section. In this final section, we give
several further results about the structure of nontrivial \nsr s
$R=R(X,Z,\mc{L},\sigma)$, and, more generally, for nontrivial
  \gnba s  $R=S(X,\{\mcII_n\},\mc{L},\sigma)$.
We prove in particular that
generic flatness always fails for $R$ and that both the point
modules in $ \rgr R$ and their analogues in $\rqgr R$ fail to be
parametrized by any scheme of locally finite type. This is in marked
contrast to Corollary~\ref{GT equiv1}(1) which shows that the latter
are in 1-1 correspondence with the closed points of $X$.

We first consider generic flatness, which is  defined as follows. If
$M$ is a module over a commutative domain $C$, then $M$ is
\emph{generically flat} over $C$ if there exists $f\in
C\smallsetminus \{0\}$ such that $M[f^{-1}]$ is a flat
$C[f^{-1}]$-module.    If $A$ is a cg $k$-algebra and $C$ is
  a commutative $k$-algebra, set $A_C = A\otimes_kC$,
 regarded as a graded $C$-algebra.

\begin{lemma}\label{non-represent1}
Let $\mc{R}=\mc{S}(X,\{\mcII_n\},\mc{L},\sigma)$  and $R=S(X,\{\mcII_n\},\mc{L},\sigma)$
 be as in  Definition~\ref{gen naive def} and suppose that $R'\subseteq R$ is a cg
 subalgebra such that  $\dim_kR/R'<\infty$. Then:
\begin{enumerate}
\item There exists $n_0\geq 0$ such that,
for any open affine subset $U \subset X$, the
$R'_{\mc{O}_X(U)}$-module $\mc{R}(U)$ is equal in degrees $\geq n_0$ to the
submodule $1\cdot R'_{\mc{O}_X(U)}$ generated by $1\in
\mc{R}(U)_0=\mc{O}_X(U)$.

\item If $R'_1\not=0$ then $\mc{R}(U) = 1\cdot R'_{\mc{O}_X(U)}$ for some
open affine  set $U\subset X$.
\end{enumerate}
\end{lemma}

\begin{proof} (1) As $\{\mc{L}_n\otimes\mcII_n\}$ is ample, there exists
 $n_0\geq 0$ such that $\mc{R}_n=\mc{L}_n\otimes\mcII_n$ is generated by its
sections $R_n$ for all $n\geq n_0$.  By hypothesis, $R_n=R'_n$ for all $n\gg 0$ so, after
possibly increasing $n_0$, we can assume that $R_n=R'_n$ for all $n\geq n_0$, as well.
Thus, for any open affine set
$U\subset X$, the element $1\in \mc{R}(U)_0$ generates $\mc{R}_n(U)
= R_n\mc{O}_X(U)=1\cdot\left(R'_{\mc{O}_X(U)}\right)_n$.

(2) In this case, for each $1\leq m< n_0$ we pick $\alpha_m\in
R'_m\smallsetminus\{0\}$,  and then we can find an open affine
subset $U_m\subset X$   such that
$\alpha_m\mc{O}_X(U_m)=\mc{L}_m(U_m)=\left(\mc{L}_m\otimes\mc{I}_m\right)(U_m)$.
So, replace $U$ by $U\cap U_1\cap\cdots \cap U_{n_0-1}$.  \end{proof}

\begin{comment}
\begin{lemma}\label{non-represent1}
Let $\mc{R}=\mc{S}(X,\{\mcII_n\},\mc{L},\sigma)$  and $R=S(X,\{\mcII_n\},\mc{L},\sigma)$
 be as in  Definition~\ref{gen naive def}. Then:
\begin{enumerate}
\item There exists $n_0\geq 0$ such that,
for any open affine subset $U \subset X$, the
$R_{\mc{O}_X(U)}$-module $\mc{R}(U)$ is equal in high degree to the
submodule $1\cdot R_{\mc{O}_X(U)}$  generated by $1\in
\mc{R}(U)_0=\mc{O}_X(U)$.

\item If $R_1\not=0$ then $\mc{R}(U) = 1\cdot R_{\mc{O}_X(U)}$ for some
open affine  set $U\subset X$.
\end{enumerate}
\end{lemma}

\begin{proof} (1) As $\{\mc{L}_n\otimes\mcII_n\}$ is ample, there exists
 $n_0\geq 0$ such that $\mc{R}_n=\mc{L}_n\otimes\mcII_n$ is generated by its
sections $R_n$ for all $n\geq n_0$. Thus, for any open affine set
$U\subset X$, the element $1\in \mc{R}(U)_0$ generates $\mc{R}_n(U)
= R_n\mc{O}_X(U)=1\cdot\left(R_{\mc{O}_X(U)}\right)_n$.

(2) In this case, for each $1\leq m< n_0$ we pick $\alpha_m\in
R_m\smallsetminus\{0\}$,  and then we can find an open affine
subset $U_m\subset X$   such that
$\alpha_m\mc{O}_X(U_m)=\mc{L}_m(U_m)=\left(\mc{L}_m\otimes\mc{I}_m\right)(U_m)$.
So, replace $U$ by $U\cap U_1\cap\cdots \cap U_{n_0-1}$.  \end{proof}

\end{comment}

We can now show that generic flatness fails for some very natural
 $R$-modules, thereby proving parts~(3) and (9) of Theorem~\ref{mainthm}.

\begin{theorem}
\label{not strong noeth} Let
$\mc{R}=\mc{S}(X,\{\mcII_n\},\mc{L},\sigma)$ be a nontrivial \gnbba\
with  $R=\mathrm{H}^0(X,\mc{R})$.
 Let $V$ be any open affine subset of $X$
 and write $C = \mc{O}_X(V)$ and
 $M = \mc{R}(V)$. We regard $M$ as a right $R_C$-module,
 with $R$ acting from the right and $C$ from the left.

Then $M$ is a finitely generated right $R_C$-module
which is not  generically flat over $C$.
It follows that  $R$ is neither strongly right noetherian nor strongly left noetherian.
\end{theorem}

\begin{proof} The proof is  similar to that
 of \cite[Theorem~9.2]{KRS}.
By Lemma~\ref{non-represent1}(1), $M$ is finitely generated as a
 right $R_C$-module.
Any localization $C[f^{-1}]$ of $C$ equals $\mc{O}_X(U)$ for an open
subset $U\subseteq V$, and so we can always replace $V$ by $U$ in
the statement of the result. In particular,  in order to prove that
$M$ is not generically flat over $C$,    it
suffices to prove that $M_n$ is not flat for  $n\gg 0$.

Consider the short exact sequence
$$0\to M_n \to \mc{L}_n(V) \to  \left(\mc{O}_X/\mcII_n \right)(V)\to 0.    $$
By nontriviality and the saturation property,  the final
term is nonzero for $n \gg 0$. Thus,  as $\mc{O}_X/\mcII_n$ is
zero-dimensional and supported at nonsingular points of $X$,
 the $C$-module $\left(\mc{O}_X/\mcII_n\right)(V)$ has
projective dimension equal to $\dim X$.  Thus, for $n\gg 0$,
  the $C$-module  $M_n$ has projective
dimension equal to $\dim X-1 \geq 1$, as required.

The second assertion  of the theorem follows from the first
  combined with  \cite[Theorem~0.1]{ASZ}.  \end{proof}

We next turn to the representability of functors, for which we need
some notation.
Let $S=\bigoplus_{n\geq 0} S_n$ be a cg $k$-algebra
and write $\ptfn(S,C)$   for
the set of isomorphism classes of graded factors $V$ of $S_C$ with
the property that  each $V_n$ is a flat $C$-module of constant rank
$h(n)=1$.  Note that  $\ptfn(S,k)$  denotes the
 point modules for $S$, as defined in the introduction.
 Moreover,   $ \ptfn(S, -)$ defines a functor from  the category
of commutative  $k$-algebras to the category of sets. Following
\cite[Section~E5]{AZ2}, we  also have an analogue of point modules
in $\rqgr S$. Specifically, let $\mathcal{P}'(S,C)$ denote the set
of isomorphism classes of  graded factors $V$ of $S_C$ with the
property that, for $n\gg 0$, the $C$-module  $V_n$ is  flat of
constant rank
 $h(n)=1$. Then write  $\qptfn(S,C)$   for the image of
$\mathcal{P}'(S,C)$   in $\rqgr S_C$.

  When $S$ is a \gnba,
it is clear from Corollary~\ref{gen-naive2}(4) that $\qptfn(S,k)$ consists of
simple objects in $\rqgr S$ but, as we show next, the converse is
also true. Combined with Corollary~\ref{GT equiv1}, this completes
the proof of Theorem~\ref{mainthm}(5).

\begin{proposition}\label{non-represent2}
Let $R=S(X,\{\mcII_n\},\mc{L},\sigma)$ be a  \gnba.
Then  $\qptfn(R,k)$ is  the set of isomorphism classes
of simple objects in $\rqgr R$.
\end{proposition}

\begin{remark}\label{non-rep-remark}
 A couple of comments about the proof are in order. One would like to claim
that the module $M(x)$ constructed in the proof of Corollary~\ref{GT
equiv1} is a point module, as this would essentially prove the
proposition. Unfortunately this is not always true. Instead, we will
use the more subtle global sections functor from \cite{AZ1}.
Unfortunately, again, this does not behave well for rings that do
not satisfy $\chi_1$, as may happen for $R$ (see
Example~\ref{easy-eg1}). So we will have to
 work simultaneously with  $R$ and its maximal right torsion extension~$T(R)$.
\end{remark}

\begin{proof}  As suggested in the remark, we
first study modules for $T=T(R)= S(X,\{\mc{H}_n\},\mc{L},\sigma)$,
as defined in \eqref{extn-def}, and its associated bimodule algebra
$\mc{T} = \mc{S}(X,\{\mc{H}_n\},\mc{L},\sigma)$.    Note that, by
Corollaries~\ref{cosupp-cor}(3) and \ref{cosupp-cor2}, $T$ is a
\gnba\ so the earlier results of the paper are available to us.

We now follow the proof of \cite[Proposition~10.7]{KRS}. Explicitly,
for a module  $\mc{N}\in \rqgr T$, write
$$\Gamma_{\mathrm{AZ}}(\mc{N})
 = \bigoplus_{m\geq 0} \hom_{\rqgr T}( \pi_T(T), \mc{N}[m])$$
 for the image of $\mc{N}$ under the Artin-Zhang global section functor \cite{AZ1}.
 By Corollary~\ref{fixing chi},
$T$ satisfies $\chi_1$ and so \cite[Theorem~4.5(2)]{AZ1} implies
that $T=\Gamma_{AZ}(\pi_T(T))$. Thus, by \cite[S2, p.~252 and S5,
p.~253]{AZ1}, $\Gamma_{AZ}(\mc{N})$ is a finitely generated
$T$-module that is torsion-free in the sense that it has no finite
dimensional submodules.

Fix a closed point $x \in X$, and recall the notation $\overline{x}=
\bigoplus_{n\geq 0} (k(x) \otimes \mc{L}_\sigma^{\otimes n})$ from
before Corollary~\ref{GT equiv1}. The proof of
\cite[Lemma~6.1]{KRS}(2) shows that $\overline{x} \in \rgr \mc{T}$
and so, by Theorem~\ref{VdB main theorem}, $\HB^0(X,\overline{x})$
is finitely generated as an $T$-module. Let $\mc{N} =
\pi_T(\HB^0(X,\overline{x}))$, considered as an element of $\rqgr T$,
and set $N(x) = \Gamma_{\mathrm{AZ}}(\mc{N}) \in \rgr T$.
 As $\HB^0(X,\overline{x})$  is noetherian, its
  maximum torsion submodule  must be finite-dimensional. Thus, as  $T$
satisfies $\chi_1$, it follows from   \cite[(3.12.3) and Proposition 3.14]{AZ1}
that  the natural map
$\HB^0(X,\overline{x}) \to N(x)$ is an isomorphism in large
degree.  In particular,
$\dim_{k}N(x)_m=1$, for $m\gg 0$.

   We next   show that $N(x)_0=\hom_{\rqgr T}( \pi(T), \mc{N})$
is nonzero.   By saturation,  we can  choose $t\geq 1$ such that
$x \not\in \bigcup_{m\geq 0}\supp \mc{O}_X/\mc{H}_m^{\sigma^t}$.  We may also assume that
$\mc{H}_m\mc{H}_n^{\sigma^m} = \mc{H}_{m+n}$ for all $m,n\geq t$. Then any
surjection of sheaves
$\theta: \mc{H}_{t} \twoheadrightarrow  \mc{H}_{t} /\mc{M}  \cong k(x)$
induces a canonical  surjection
$\theta_m : \mc{H}_{t+m} = \mc{H}_t\mc{H}_m^{\sigma^t}  \twoheadrightarrow
 \mc{H}_{t+m} +\mc{M}/\mc{M} \cong k(x),$
and hence a surjection of $\mc{O}_X$-modules
$\theta_m\otimes \mathrm{Id} : \mc{T}_{t+m} =
\mc{H}_{t+m}\otimes \mc{L}^{t+m}_\sigma
  \twoheadrightarrow  k(x) \otimes \mc{L}^{t+m}_\sigma$ for all $m\geq t$.
These  are the structure maps for a surjective homomorphism $f:
\mc{T}_{\geq 2t} \to \overline{x}_{\, \geq 2t}$ in $\rgr \mc{T}$. Finally,
by taking global sections and passing to $\rqgr T$, the morphism $f$
induces a nonzero element of $N(x)_0$.

Now consider $N(x)$ as an $R$-module  and fix a nonzero element $a\in
N(x)_0$. We claim that $aR$ is not torsion. Indeed, otherwise
$aR_{\geq r}=0$ for some $r\geq 1$. But, Lemma~\ref{cosupp-cor2}(3)
implies that $T$ is a finitely generated left $R$-module and so
$T/R_{\geq r}T$ is finite dimensional as a left (and therefore
right) $k$-module. Hence $T_{\geq s}\subseteq R_{\geq r}T$ for some
$s$. Thus $aT_{\geq s}=0$, contradicting the fact that $aT\subseteq
N(x)$ is a torsion-free right $T$-module.

So, $aR$ is not torsion. By Corollary~\ref{gen-naive2}(4), there
exists $u$ such that $R_mR_n=R_{m+n}$ for all $m,n\geq u$. It
follows that $a R_n \neq 0$ for all $n \geq u$. (To see this, note
that if $a R_n  = 0$ for some $n \geq u$, then $0=a R_n R_m = a
R_{n+m}$ for all $m \geq u$, which leads to the contradiction
 $a R_{\geq (n + u)} = 0$.)
 Since $aR_n\subseteq N(x)_n$ and $\dim_kN(x)_n = 1$ for $n\gg 0$, it
follows that $aR_n=N(x)_n$ is $1$-dimensional for all $n \gg 0$. In
particular, $N(x)$ is a finitely generated right $R$-module.
Moreover, if $\pi_R$ denotes the natural map $\rgr R\to \rqgr R$,
then $\pi_R(aR)=\pi_R(N(x))\in \qptfn(R,k)$.

Finally, since $\HB^0(X,\overline{x})$ and $N(x)$ are isomorphic in
large degree, $\pi_R(aR)$ is also equal to
$\pi_R(\HB^0(X,\overline{x}))$.
 But, if we use the equivalence of categories, Theorem~\ref{VdB main theorem},
to identify $\rqgr R$ with $\rqgr \mc{R}$, then
$ \pi_R(\HB^0(X,\overline{x}))=\pi_{\mc{R}} (\overline{x})=\widetilde{x}$.
If $R$ is a \nsr\ then,  by Corollary~\ref{GT
equiv1}(1), the $\widetilde{x}$ are also just the simple
objects in $\rqgr R$. In other words, the set of isomorphism classes
of  simple objects in $\rqgr R$ is  just  $\qptfn(R,k)$
 as is required to prove the theorem.
If $R$ is not a \nsr, apply Lemma~\ref{gen-naive1} to
 pick  $t\geq 1$ such that $R^{(t)}$ is one.
Then Corollary~\ref{GT equiv1}(1) can be applied to show that  the
simple objects in $\rqgr R^{(t)}$
 are just the images of the closed points in $X$; that is the
objects $\pi \big( \bigoplus_{n\geq 0} k(x)\otimes
(\mc{L}_{t})_{\sigma^{t}}^{\otimes n} \big) =\widetilde{x}^{(t)} $,
 in the notation of this proof. However, by
Lemma~\ref{nonrep-lemma}(1),
 $\rqgr R \simeq \rqgr R^{(t)}$ via the functor $M\mapsto M^{(t)}$.
Thus, the simple objects in $\rqgr R$ are still the
 $\widetilde{x} $ for $ x\in X$,  as we needed.
\end{proof}

\begin{theorem}\label{non-represent}
Let $R=S(X,\{\mcII_n\},\mc{L},\sigma)$ be a  nontrivial \gnba\
and suppose that $R'\subseteq R$ is a cg subalgebra such that $\dim_kR/R'<\infty$.
\begin{enumerate}
\item
If $R'_1\not=0$ then  $\ptfn(R',-)$ is not represented by any scheme $Y$ of locally finite type.
\item Whether  $R'_1=0$ or not,  $\qptfn(R',-)$ is
 not represented by any scheme $Y$ of locally finite type.
\end{enumerate}
\end{theorem}

\begin{remarks} (1) This proves Theorem~\ref{mainthm}(6).

(2) If $R'_1=0$, then part (1) of the theorem will fail. Indeed, in this case, given  \emph{any}
commutative $k$-algebra $C$ and  $R'_C$-module
$M$ generated in degree zero, then   $M_1=0$. In other words, there are no
point modules for $R'$ and $\ptfn(R',-)$ is represented by the empty scheme.
\end{remarks}

\begin{proof}    The proof is similar to that of
\cite[Theorem~10.4 and Corollary~10.5]{KRS} and, as there, the idea
of the proof is that, for any open affine $U\subset X$, the module
$\mc{R}(U)$ is ``trying but failing'' to be the module corresponding
to the commutative ring $\mc{O}_X(U)$. We need to make this
assertion formal.

Assume that $\mathcal{P}(-) = \ptfn(R',-)$ is represented by the
scheme $Y$ of locally finite type. Pick an open affine set $U\subset
X$ by Lemma~\ref{non-represent1}(2), fix a closed point $p\in
U\smallsetminus \bigcup_{m\in \mathbb{Z}} \supp
\mc{O}_X/\mcII_1^{\sigma^m}$ and set $C =\mc{O}_{X,p}$.  Then
$(\mc{L}_n\otimes \mcII_n)_p \cong C$ for all $n$ and so
$\mc{R}_p=\mc{R}(U)\otimes_{\mc{O}_X(U)}C \cong\bigoplus_{n\geq
0}C$. By Lemma~\ref{non-represent1}(2), $\mc{R}_p$ is generated as
an $R'_C$-module by the element $1$ in degree zero, so $\mc{R}_p\in
\mc{P}(C)$. Thus there exists $\theta_p\in
\mathcal{P}(C)=\Morph(\spec C, Y)$ corresponding to~$\mathcal{R}_p$.

By the definition of locally finite type \cite[p.84]{Ha}, we may pick an open affine
neighbourhood $V$ of $\theta_p(p)$ in $Y$ of finite type over $k$.
Then  we get a map of algebras
$\theta_p' : \OO_Y(V)\to  \OO_{\spec C}(\theta_p^{-1}(V))$.
Since $\theta_p^{-1}(V)$ is an
open set containing $p$, it is necessarily $\spec C$ and so
 $\mathrm{Im}(\theta_p')\subseteq C$.
Since $\OO_Y(V)$ is a finitely generated $k$-algebra
  and $\OO_X(U)$ is a domain,
 $\theta_p'(\OO_Y(V))\subseteq \OO_X(U')$, for some open set $U'\subseteq U$.
 Since it does no harm to replace $U$ by a smaller open set containing $p$,
 we may as well assume that $U=U'$.
 In other words, we have extended $\theta_p$ to a map
 $\widetilde{\theta}_p \in \Morph(U,Y)$ such that $\theta_p=\widetilde{\theta}_p\circ \pi_p$,
  where $\pi_p : \spec C\to U$ is the natural
morphism.

By construction, $\widetilde{\theta}_p$ corresponds to a module
$M_U\in \mathcal{P}(\OO(U))$
with the property that $M_U\otimes_{\OO(U)}C \cong \mathcal{R}_p$.
But $\mathcal{R}(U)$ is a second finitely generated $R'_{\OO(U)}$-module
with $\mathcal{R}(U)\otimes_{\OO(U)}C \cong \mathcal{R}_p$. This local
isomorphism of $R'_C$-modules  lifts to an isomorphism
$M_W=M_U\otimes_{\OO(U)}\OO(W) \cong \mathcal{R}(W)$ of $R'_{\OO(W)}$-modules,
 for some open affine set
$W\subseteq U$.  By the definition of $\mathcal{P}$, the $\OO(W)$-module
$(M_W)_n=(M_U)_n\otimes_{\OO(U)}\OO(W)$ is flat for
all $n$. On the other hand, for $n\gg 0$, the proof of
Theorem~\ref{not strong noeth}
implies that  $(M_W)_n\cong \mathcal{R}(W)_n$ is \emph{not} flat over
$\OO(W)$.  This contradiction proves (1).

(2) To begin, assume that $R'_1\not=0$ and
consider the proof of part~(1).
In the final paragraph of that proof,
$M_W \in \mathcal{P}(\OO(W))$ and so $\pi(M_W)$ certainly lies in
$\qptfn(R',\OO(W))$.
In contrast, as $\mathcal{R}(W)_n$ is not flat
as an $\OO(W)$-module for any $n\gg 0$, no tail $\mathcal{R}(W)_{\geq n}$
of $\mathcal{R}(W)$ is a flat $\OO(W)$-module. Hence
$\pi(\mathcal{R}(W))$ cannot belong to
$\qptfn(R',\OO(W))$.
Thus, the proof of part~(1)  also proves part~(2).

If $R'_1=0$ then the same proof works, except that one now uses
Lemma~\ref{non-represent1}(1) in place of
Lemma~\ref{non-represent1}(2) and, for each module $N$ that appears
in the proof,
 one ignores the terms $N_n$ for   $0<n<n_0$.
 \end{proof}

\begin{remark}\label{intro-remark-proof}
To  end the paper we justify the comments made in Remark~\ref{intro-remark}.
 So, assume that  the hypotheses (and conclusions)
 of Theorem~\ref{mainthm}  hold for a \nsr\
$R=R(X,Z,\mc{L},\sigma)$  and let $R'\subseteq R$ be a cg subalgebra such that $\dim R/R'<\infty$.
We need to prove that the conclusions of  that theorem also hold for $R'$.

First of all,  it is routine that  $R'$ is noetherian, proving 
part~(2), while part~(3) is trivial. By \cite[Proposition~2.7]{SZ2},
$\rqgr R'=\rqgr R$ and so   parts (1,4,7,8)  immediately hold for $R'$.
  Part~(5) and hence the first part of~(6) are  also  easy exercises.
Moreover, the module $\mc{R}(U)$ is still a finitely generated
$R'\otimes_k\mc{O}_X(U)$-module, so part (9) also holds for $R'$.
 Thus, it only remains to prove that  the point modules for $R'$  are
 not parametrizable, and this  was proved directly in Theorem~\ref{non-represent}.
 \end{remark}

\providecommand{\bysame}{\leavevmode\hbox to3em{\hrulefill}\thinspace}
\providecommand{\MR}{\relax\ifhmode\unskip\space\fi MR }
\providecommand{\MRhref}[2]{
  \href{http://www.ams.org/mathscinet-getitem?mr=#1}{#2}
}
\providecommand{\href}[2]{#2}

\end{document}